\documentclass[11pt,reqno]{amsart}

\usepackage{graphicx}
\usepackage{amsmath}
\usepackage{amssymb}
\usepackage{amsfonts}
\usepackage{amsthm}
\usepackage{epstopdf}

\newtheorem{theorem}{Theorem}

\newtheorem{conjecture}{Conjecture} 
\theoremstyle{remark}
\newtheorem{remark}{Remark}[section]

\newcommand{\defeq}{\stackrel{\rm{def}}{=}}

\begin{document}
\title[Zakharov-Kuznetsov equation in 3D]
{Numerical study of soliton stability, resolution and interactions in the 3D Zakharov-Kuznetsov equation}

\author[C. Klein]{Christian Klein}
\address{Institut de Math\'ematiques de Bourgogne, UMR 5584; \\
Universit\'e de Bourgogne-Franche-Comt\'e, 9 avenue Alain Savary, 21078 Dijon
                Cedex, France} 
\email{Christian.Klein@u-bourgogne.fr}

\author[S. Roudenko]{Svetlana Roudenko}
\address{Department of Mathematics \& Statistics\\Florida International University,  
Miami, FL, 33199, USA}
\curraddr{}
\email{sroudenko@fiu.edu}
    
\author[N. Stoilov]{Nikola Stoilov}
\address{Institut de Math\'ematiques de Bourgogne, UMR 5584; \\
Universit\'e de Bourgogne-Franche-Comt\'e, 9 avenue Alain Savary, 21078 Dijon
                Cedex, France} 
\email{Nikola.Stoilov@u-bourgogne.fr}


\begin{abstract}
We present a detailed numerical study of solutions to the
Zakharov-Kuznetsov equation in three spatial dimensions. The equation is a three-dimensional generalization of the Korteweg-de Vries equation, though, not completely integrable.
This equation is $L^2$-subcritical, and thus, solutions exist globally, for example, in the $H^1$ energy space.  

We first study stability of solitons with various perturbations in 
sizes and symmetry, and show asymptotic stability and formation of 
radiation, confirming the asymptotic stability result in 
\cite{FHRY2020} for a larger class of initial data. We then 
investigate the solution behavior for different localizations and rates of decay including exponential and algebraic decays, and give positive confirmation toward the soliton resolution conjecture in this equation. Finally, we investigate soliton interactions in various settings and show that there is both  
a quasi-elastic interaction and a strong interaction when two solitons merge into one, in all cases always emitting radiation in the conic-type region of the negative $x$-direction.   
\end{abstract}

\maketitle

\section{Introduction}
We are interested in the 3D quadratic Zakharov-Kuznetsov (ZK) equation
\begin{equation}\label{ZK}
u_t + (u_{xx}+u_{yy}+u_{zz} + u^2)_x  = 0,
\end{equation}
where $u=u(x,y,z,t)$ is real-valued, $(x,y,z) \in \mathbb R^3$, and $t \in \mathbb R$. 
This equation is a three-dimensional generalization of the well-known 
Korteweg-de Vries (KdV) equation, which is a one-dimensional model 
for weakly nonlinear waves in shallow water. The 3D ZK equation was 
originally proposed by Zakharov and Kuznetsov in the description of 
weakly magnetized ion-acoustic waves in a low-pressure magnetized 
plasma \cite{ZK1974}, where they raised the  question of studying 
solitons in a higher-dimensional setting. In particular, their main 
question was about the stability of solitons, for which they argued 
that for the 3D ZK equation the Lyapunov-type functional attains its minimum on a soliton. The orbital stability of solitons was obtained by de Bouard \cite{deB} by adapting the KdV argument of Grillakis, Shatah \& Strauss \cite{GSS} to the 2D and 3D ZK equation. 
The more delicate asymptotic stability of solutions (in $H^1$) close 
to a soliton for the 3D ZK equation was recently obtained by the 
second author together with Holmer, Farah and Yang in 
\cite{FHRY2020}. It is the goal of the present work to investigate 
soliton formation, stability and interaction in the 3D ZK equation numerically. 

While originally the equation was proposed by Zakharov and Kuznetsov 
in the 3D setting, the first rigorous derivation as a long-wave 
small-amplitude limit of the Euler-Poisson system in the cold-plasma approximation was done by Lannes, Linares and Saut in \cite{LLS}, see also \cite{LS}. Other derivations exist as well, for a review see \cite{LLS, FHRY, FHR3} and references therein. 

Unlike KdV and its other generalizations such as 
Kadomtsev-Petviashvili or Benjamin-Ono equations,  the Zakharov-Kuznetsov equation is not completely integrable. Nevertheless, it has a Hamiltonian structure with three conserved quantities: 
energy (Hamiltonian), $L^2$-norm (often called mass) and the integral, defined as follows
\begin{align}\label{MC}
M[u(t)] &\defeq\int_{\mathbb{R}^3} u^2(t) 
= M[u(0)], \\
\label{EC}
E[u(t)] & \defeq \dfrac{1}{2}\int_{\mathbb{R}^3}
\big[u_{x}^{2}(t)+u_{y}^2(t)+u_z^2(t)\big]  
- \dfrac{1}{3}\int_{\mathbb{R}^3} u^{3}(t)  
= E[u(0)], \\
\label{L1-inv}
\int_{\mathbb{R}}  u(x, & \, y, z,t) \, dx = \int_{\mathbb{R}} u(x,y,z,0) \, dx.
\end{align}

The equation \eqref{ZK} has a scaling invariance: 
if $u(x,y,z,t)$ is a solution of \eqref{ZK}, then so is the rescaled version 
\begin{equation}\label{E:scaling}
u_\lambda(x,y,z,t)=\lambda^2 u(\lambda x, \lambda y, \lambda z,\lambda^3 t), ~~\lambda > 0.
\end{equation}
This symmetry makes the Sobolev norm $\dot{H}^s$ with $s=-1$ invariant, thus, making the equation \eqref{ZK} $L^2$-subcritical ($s<0$). 
The 3D ZK equation has other invariances such as translation and dilation.

The well-posedness theory for the Cauchy problem for the 3D ZK equation 
with $H^s$ initial data has attracted significant interest in the last decade. 
The local well-posedness can be established via the classical Kato method in $H^s$ for $s>\frac52$. This was remarked and improved by Linares \&  Saut in \cite{LS} to the local well-posedness in $H^s$ with $s>\frac98$ following the method of Kenig \cite{K2004}, which was then further improved by Ribaud \& Vento \cite{RV} down to $H^s$ with $s>1$.  The global well-posedness in  $H^s$, $s>1$, was established by Molinet \& Pilod \cite{MP}, and had been open for a while, until the recent work of Herr \& Kinoshita \cite{HK}, obtaining  the local well-posedness in $H^s$ for $s>-\frac12$. We are interested in studying finite energy solutions, hence, $H^1$ global well-posedness suffices for our purposes.

The equation has a family of traveling waves called solitary waves (sometimes called solitons, although the model is not integrable), moving only in the positive $x$-direction:
\begin{equation}\label{E:wave}
u(x,y,z,t) = Q_c(x - c \, t, y, z), ~~c>0,
\end{equation}
where $Q_c$ is the dilation 
\begin{equation}
	Q_c(\cdot) = c \, Q(\sqrt{c} \, \cdot).
	\label{dilation}
\end{equation}
We only consider solitary waves vanishing at infinity, thus, $Q$ is the vanishing at infinity ground state solution of the well-known nonlinear elliptic equation 
\begin{equation}\label{E:Q}
-\Delta_{\mathbb R^3} Q + Q  - Q^2 = 0,
\end{equation}  
i.e., the unique radial positive smooth solution in $H^1(\mathbb{R}^3)$. 
The properties of this ground state include $Q \in C^{\infty}(\mathbb{R}^3)$, $\partial_r Q(r) <0$ for any $r = |(x,y,z)|>0$, and for any multi-index $\alpha$
\begin{equation}\label{prop-Q}
|\partial^\alpha Q(x, y, z)| \lesssim_\alpha e^{-|(x,y,z)|} \quad \mbox{for any}\quad {(x,y,z)} \in \mathbb{R}^3.
\end{equation}

As it was mentioned in the beginning, these solitary waves are stable, both orbitally (by a result of de Bouard \cite{deB}) and also asymptotically stable by the result of the second author with Farah, Holmer and  Yang \cite{FHRY2020}:
\begin{theorem}[\cite{FHRY2020}]\label{T:1}
For $\alpha \ll 1$ and $u_0\in H^1(\mathbb R^3)$ with 
$\|u_0 - Q \|_{H^1} \leq \alpha$,
the solution $u(x,t)$ to the 3D ZK 
\eqref{ZK} is \emph{asymptotically stable}: 

$\bullet$ \emph{(orbital stability)} there exist trajectories $c(t)>0$ and $(a_1(t), a_2(t), a_3(t)) \in \mathbb{R}^3$ such that 
$$
\big\| c^2(t) \, u\big(c(t)x +a_1(t), c(t) y +a_2 (t), c(t) z + a_3(t), t\big) - 
Q(x,y,z) \big\|_{H^1} \lesssim \alpha,
$$

$\bullet$ \emph{(convergence of trajectories)} there exists $c_*$ such that $|c_*-1| \lesssim \alpha$ such that 
$$
c(t) \to c_*   \quad \text{and} \quad \big(a_1'(t), a_2'(t), a_3'(t) \big) \to (c_*^{-2},0,0) 
\quad \text{as} \quad t\to +\infty, 
$$

$\bullet$ \emph{(weak convergence)} as $t \to + \infty$ 
\begin{equation}
\label{E:main-weak}
c^2(t) \, u\big(c(t)x +a_1(t), c(t) y +a_2 (t), c(t) z + a_3(t), t \big) \rightharpoonup Q(x,y,z) ~~ \text{ (weakly) in } H^1, 
\end{equation}

$\bullet$ \emph{($L^2$ strong convergence)} for any $\delta\gtrsim \alpha$ the strong convergence in $L^2$ 
holds outside of the radiation cone (i.e., on the conic right-half space $\mathcal C$)
\begin{equation}
\label{E:window-conv-b}
\big \| c^2(t) \, u(c(t)x +a_1(t), c(t) y +a_2 (t), c(t) z + a_3(t), t) - 
Q(x,y,z) \big \|_{L^2_{\mathcal C}} \to 0 \text{ as } t\to +\infty,
\end{equation}
where 
\begin{equation}\label{E:cone}
\mathcal C \defeq \{ (x,y,z) \in \mathbb R^3: ~ x > (-1+\delta) t -\sqrt{y^2+z^2}\tan \theta \}
\end{equation}
for all $\theta$ such that 
\begin{equation}\label{E:angle}
0\leq  \theta \leq \frac{\pi}{3}-\delta.
\end{equation}
\end{theorem}
We remark that the $L^2$ convergence stated in \eqref{E:window-conv-b} is in the reference frame of the soliton (being at the origin); in the reference frame of the solution, the rightward shifting external conic region is  $x> \delta t -\sqrt{y^2+z^2} \, \tan \theta$. 

Theorem \ref{T:1} describes the asymptotic behavior of solutions very close (in $H^1$ sense) to solitons. In order to understand a more general picture of the behavior of solutions in this equation, we perform a detailed numerical study. The main goal of this study is to investigate 
stability of solitons, soliton resolution, the radiation regime and the interaction of solitons in this 3D non-integrable model via numerical approaches. 
\begin{conjecture} \label{C:1} 
Consider the 3D ZK equation \eqref{ZK}.
\begin{enumerate}
\item[I.]
The soliton solutions \eqref{E:wave}-\eqref{dilation}-\eqref{E:Q} are orbitally and asymptotically stable. The stability holds on a large set of perturbations of soliton-like initial data.
\item[II.] 
Solutions of \eqref{ZK} with sufficiently localized and smooth initial data decompose  into solitons and radiation as $t\to\infty$, i.e., the soliton resolution conjecture holds for the 3D ZK equation.
\end{enumerate}
\end{conjecture}

In this work we give positive numerical confirmations to Conjecture 1, as well 
as to the asymptotic stability result in Theorem 1 showing that 
asymptotic stability holds on a much larger class of solutions.
For the asymptotic stability we consider perturbations of solitons within $10\%$ difference in mass or $L^\infty$ norm. To study the soliton resolution, in our simulations we are able to consider not only exponentially decaying initial data, but also data with a sufficiently rapid algebraic decay. We have also investigated solutions with initial data not necessarily having a maximum peak at a single point but on a continuous set (such as an interval). However, we only studied positive valued initial data, as this is the first detailed numerical study for the 3D ZK equation. 

Besides the soliton stability and soliton resolution, we also investigate the soliton interactions. We observe two types of interactions: one that preserves the number of solitons before and after the interaction (though the soliton parameters such as amplitude or speed can change) -- we term this interaction as {\it quasi-elastic}, and another one that combines the two solitons into one single soliton (and radiation, which is present in any interaction) -- we refer to that interaction as {\it strong}.

For completeness, we mention that the 2 dimensional ZK equation with 
different powers of nonlinearity has been studied intensively in 
various aspects and some questions that have been answered in a two-dimensional setting would be interesting to investigate in three dimensions (and higher). For well-posedness results in the 2D quadratic ZK see \cite{F95}, \cite{LP2009}, \cite{MP2015}, \cite{GH2014}, \cite{K2019}, for other powers of ZK see \cite{LP2009}, \cite{LP2011}, \cite{FLP2012}, \cite{BFR}, \cite{K2019b}; uniqueness results have been studied in \cite{CFL}, propagation of regularity in \cite{LP}; orbital stability in 2D quandratic ZK was done in \cite{deB}, the asymptotic stability in \cite{CMPS}, in the same paper the authors investigate  the behavior of $N$ well-separated solitons, two soliton interaction in \cite{V2020}.  For instability of solitons in the critical and supercritical regimes (i.e., with cubic power and higher), see \cite{FHR2}, \cite{FHR4}, \cite{FHR3}, and existence of blow-up in the critical regime (2D cubic ZK) was proved in \cite{FHRY}. We have also investigated numerically soliton resolution, interaction, as well as the blow-up behavior in the 2D generalized ZK equations in \cite{KRS}. 


The paper is organized as follows: in Section \ref{S:Num} we present the numerical methods used to solve the 3D ZK equation. In Section \ref{S:Stability} we first study stability of solitons, where we consider various perturbations of the soliton itself. We then consider different types of initial data including exponential and algebraic decays as well as different localization features and show the soliton resolution and formation of radiation. Section \ref{S:Interaction} contains our study of the soliton interaction in different settings, including quasi-elastic and strong interactions.   
\smallskip

{\bf Acknowledgements.} CK and NS were partially supported by the ANR-FWF project ANuI - ANR-17-CE40-0035, the isite BFC project NAANoD, the EIPHI Graduate School (contract ANR-17-EURE-0002) and by the European Union Horizon 2020 research and innovation program under the 
Marie Sklodowska-Curie RISE 2017 grant agreement no. 778010 IPaDEGAN.
SR was partially supported by the NSF grant DMS-1927258.

\section{Numerical methods}\label{S:Num}
We start with a brief review of numerical methods used in this work.  
For the spatial discretisation we use Fourier spectral methods, i.e., 
we approximate the solutions via trigonometric polynomials in each 
spatial variable. Solitons to the 3D ZK equation are constructed 
after this discretisation via a Newton-Krylov approach.  The integration of the ZK 
equation in time is then performed with a fourth order exponential 
time differencing method. 

\subsection{Solitons}
Fourier spectral methods are known to be efficient in the 
approximation of 
sufficiently decreasing smooth functions. This means we choose spatial 
periods of sufficient size such that the studied function as well as its 
first derivatives are small (ideally, on the order of the machine 
precision, which is of the order $10^{-16}$ here) at the boundaries 
of the computational domain. It is known that the Fourier 
coefficients of such functions will be rapidly decreasing. Specifically,  
we work with $x\in L_{x}[-\pi,\pi]$, $y\in L_{y}[-\pi,\pi]$, and $z\in L_{z}[-\pi,\pi]$, 
where $L_{x}$, $L_{y}$ and $L_{z}$ are positive real numbers. We denote the Fourier variables 
dual to $x$, $y$ and $z$ by $k_{x}$, $k_{y}$ and $k_{z}$, 
respectively, and approximate a function $Q(x,y,z)$ via
\begin{equation}\label{Qhat}
Q(x,y,z)\approx \sum_{k_{x}=-N_{x}/2+1}^{N_{x}/2} \,
\sum_{k_{y}=-N_{y}/2+1}^{N_{y}/2} \,  \sum_{k_{z}=-N_{z}/2+1}^{N_{z}/2}
\hat{Q}(k_{x},k_{y},k_{z})\, e^{i(k_{x}x+k_{y}y+k_{z}z)}.
\end{equation}
The \emph{discrete Fourier transform} $\hat{Q}=\mathcal{F}Q$ can be conveniently 
computed with a \emph{fast Fourier transform} (FFT). An advantage of 
Fourier methods is that the numerical resolution can be controlled 
via the decay of the Fourier coefficients, the highest coefficients 
(the coefficients with indices $|k_{x}|\sim N_{x}/2$, $|k_{y}|\sim 
N_{y}/2$, and $|k_{z}|\sim N_{z}/2$) 
indicate the numerical error introduced by the truncation of the 
series. 

We first obtain the soliton solution for the equation \eqref{ZK} by 
solving the equation \eqref{E:Q}. This is a well-known nonlinear elliptic equation, it is the same equation for the ground state solution in the nonlinear Schr\"odinger or Klein-Gordon equations in 3D and has non-trivial solutions with radial symmetry. Here, however, we do not use this fact: since we apply Fourier methods throughout the paper, we directly construct the soliton solutions  on the given numerical grids to use that later either as the initial data for the ZK time evolution or for the appropriate fitting of a rescaled soliton.

With the Fourier discretization \eqref{Qhat}, the equation \eqref{E:Q} 
is approximated by an $N_{x} N_{y} N_{z}$ dimensional system of 
nonlinear equations for the $\hat{Q}$. The latter is iteratively 
solved by a Newton-Krylov iteration as in \cite{KRS}. This means that we invert the 
Jacobian via Krylov subspace methods as in \cite{AKS}, here, GMRES \cite{gmres}.  
We use $N_{x}=N_{y}=N_{z} = 2^{7}$, $L_{x}=L_{y}=L_{z}=3$ and $Q=2\,e^{-(x^{2}+y^{2}+z^{2})}$ as the initial iterate. The iteration is stopped when the residual is smaller than $10^{-10}$. 

\begin{figure}[!htb]
\includegraphics[width=0.49\hsize]{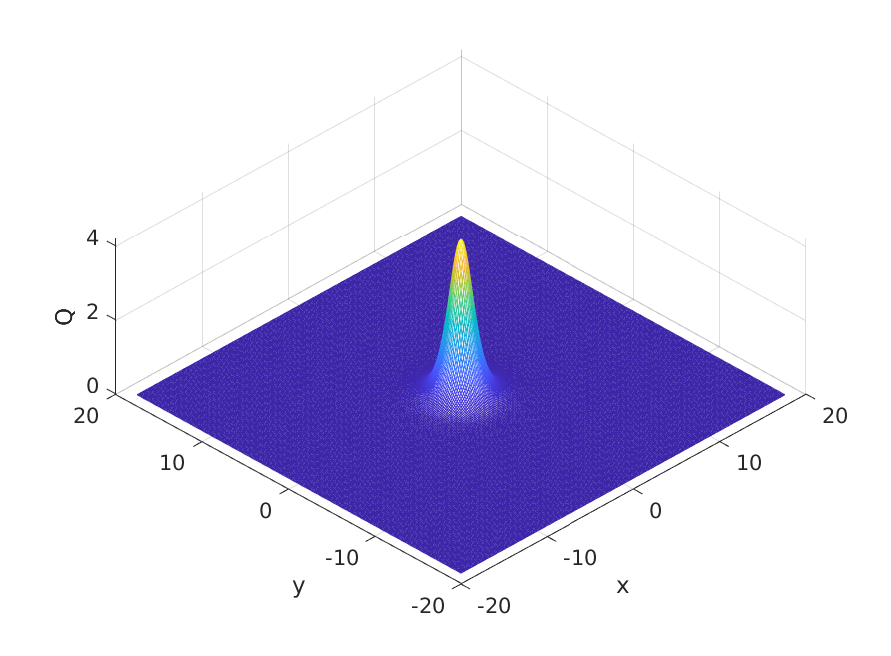}
\includegraphics[width=0.49\hsize]{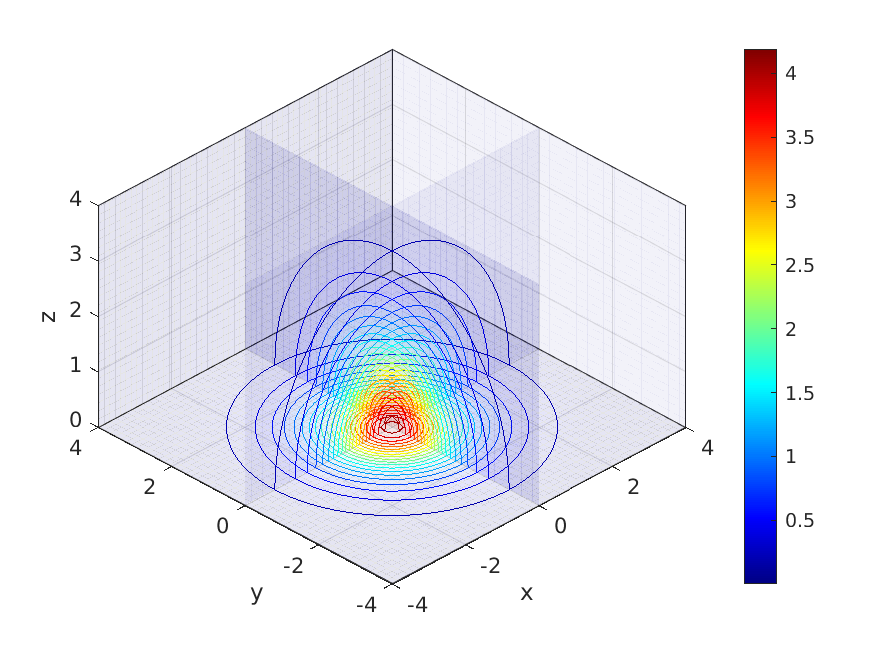}
\caption{The ground state solution to \eqref{E:Q}. 
Left: plot of $Q$ with $z = 0$. Right: 3D contour plots of $Q$ on the slices of the coordinate planes. The color bar indicates the magnitude of the solution.}
\label{solfig}
\end{figure}

The ground state solution of $-c \,Q + \Delta_{\mathbb R^3} Q + Q^2=0$ for $c=1$  
is shown in Fig.~\ref{solfig}. 

\subsection{Time evolution}
The Fourier discretization (\ref{Qhat}) is used also for 
the full ZK equation \eqref{ZK}, which is thus approximated by an 
$N_{x}N_{y}N_{z}$ dimensional system of ordinary differential equations in $t$ of the form
\begin{equation}
    \hat{u}_{t}= \mathcal{L}\hat{u}+\mathcal{N}[\hat{u}]
    \label{sys},
\end{equation}
where $\mathcal{L}=ik_{x}(k_{x}^{2}+k_{y}^{2}+k_{z}^2)$ and 
$\mathcal{N}[\hat{u}]=-ik_{x}\mathcal{F}(u^{p})$. Because of the 
appearance of third derivatives in $x$, $y$, and $z$, 
this system is \emph{stiff}, implying that explicit methods will be 
inefficient due to stability conditions as they necessitate prohibitively small times 
steps in order to stabilize the code. 
Implicit schemes are less restrictive in this sense, but are 
computationally expensive, since the resulting nonlinear equation has 
to be solved in each time step. 
In \cite{etna,KR} we compared various adapted integrators for stiff systems with a 
diagonal $\mathcal{L}$ as we have here, which 
are explicit and of fourth order. It turned out that 
\emph{exponential time differencing} (ETD) schemes, see \cite{HO} for a comprehensive 
review with many references, are most efficient  in the context 
of the KdV-type equations. There are various fourth order ETD methods, which 
all showed a similar performance in our tests. As in \cite{KRS}, we apply the method 
by Cox and Matthews \cite{CM} in the implementation described in 
\cite{etna, KR}. The accuracy of the time integration scheme can be 
controlled via the conserved energy of the equation. Due to 
limitations in the accuracy of numerical methods, the 
computed energy (again Fourier 
techniques are applied to \eqref{EC}) will not be 
exactly conserved. The quantity $\Delta E = |E(t)/E(0)-1|$ can be 
used as discussed in \cite{etna,KR} as an estimate of the numerical 
error. Typically, it overestimates the accuracy of the numerical 
solution by 1-2 orders of magnitude. 

\begin{remark}\label{remperiodic}
Note that in this paper we  approximate solutions on 
$\mathbb{R}^{3}$ by simulations on the torus. 
Within machine 
precision, this does not make a difference for  
rapidly decreasing solution, if sufficiently large periods are chosen. 
This is, for instance, possible for stationary localized 
solutions as the solitons of the ZK equation. However, if radiation 
appears in non-stationary solutions, one would have to 
choose prohibitively large computational domains to avoid the 
reappearance of emitted radiation (always emitted in the negative 
$x$-direction) for positive values of $x$, which is in practice 
impossible for 3D computations. The reappearence of radiation is acceptable as long 
as it has much smaller amplitudes than the studied bulk of 
the solution. 
\end{remark}

\subsection{Test}

To test the time evolution code and the soliton solution at the same time, 
we consider the ground state solution as the initial condition. 
For $t\in[0,1]$ we apply $N_{t}=1000$ time steps. The numerically computed energy is conserved 
to the order of $10^{-14}$. The difference between the numerically 
computed solution and the soliton 
increases with time, but it is of the same order  
$10^{-14}$ as in the energy conservation. The code is thus
able to propagate solitons on the considered time intervals,   
essentially with machine precision.

\section{Soliton Stability and Soliton Resolution}\label{S:Stability}

In this section we consider initial data with monotone decay and a 
single maximum value (though not necessarily assumed at a single 
point) and investigate how the time evolutions of such data resolve 
into coherent structures (here, rescaled and shifted solitons) and radiation. We recall that the ZK equation is not integrable, so we always expect some radiation to form unlike for integrable equations such as the KdV equation. We first consider examples of perturbed soliton initial data and then examine various other settings with different rates of decay and symmetry. 

\subsection{Soliton Stability}
We start with the investigation of the soliton stability. For that we consider the following examples of initial data with perturbed solitons: 
\begin{itemize}
\item[(a)] 
a multiple of the soliton, $u(x,y,z,0) = \lambda \, Q(x,y,z)$, $\lambda \approx 1$;
in particular, considering $\lambda \gtrsim 1$ and $\lambda \lesssim 1$; 

\item[(b)] 
a localized asymmetric perturbation of the soliton,
$u(x,y,z,0)=Q(x,y,z) +e^{-(x^2+y^2 + \alpha \,z^2)}$, $\alpha \neq 1$, $\alpha>0$.
\end{itemize}
We work in this section with \( L_{x}=L_{y}=6 \) and \( 
N_{x}=N_{y}=N_{z} = 2^{8} \) 
Fourier modes and \( N_{t}=10,000 \) time steps on the considered time 
intervals. In all studied cases the Fourier coefficients decrease at 
least to the order of $10^{-10}$, and the relative energy is conserved 
at least to the same order.  
We use a co-moving frame (with the unperturbed soliton), i.e., we solve
\begin{equation}\label{ZKcom}
u_t + (\Delta u + u^2 - v_{x}u)_x  = 0 
\end{equation}
with $v_{x}=1$.
\smallskip

\underline{\bf Case (a): multiple of a soliton with $\lambda \gtrsim 1$.} 
\smallskip

\begin{figure}[!htb]
\includegraphics[width=0.51\hsize]{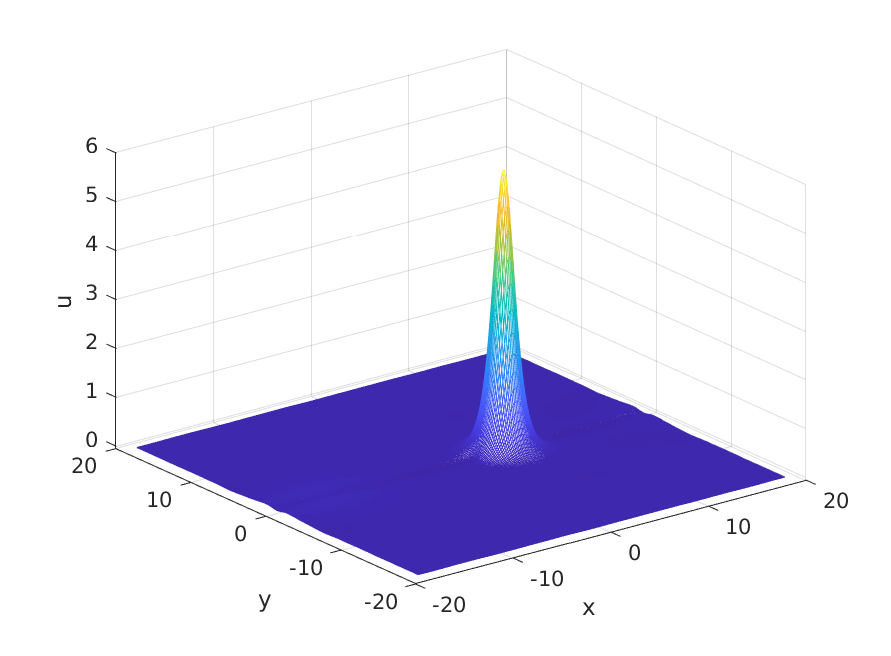}
\includegraphics[width=0.46\hsize]{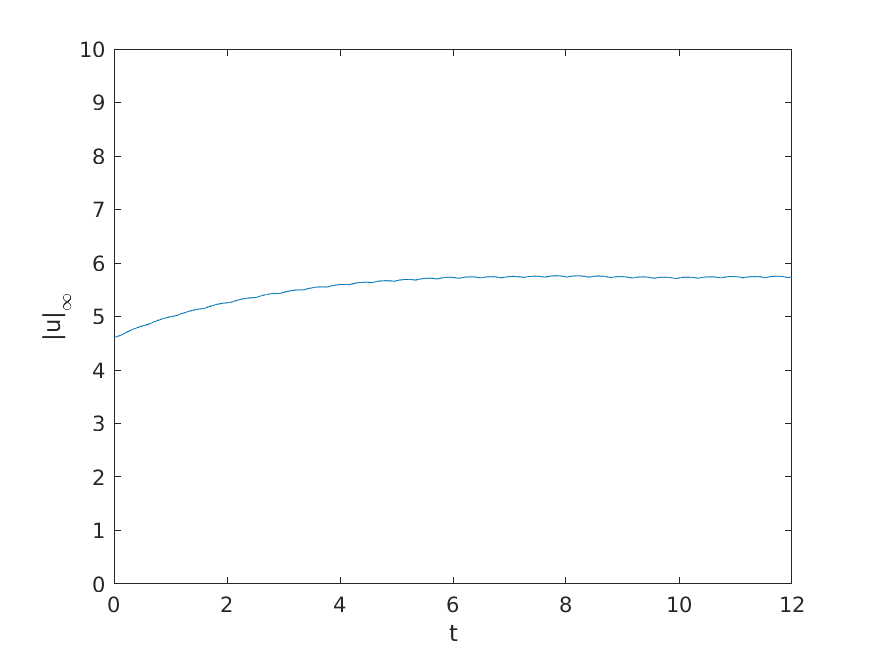} 
\includegraphics[width=0.49\hsize]{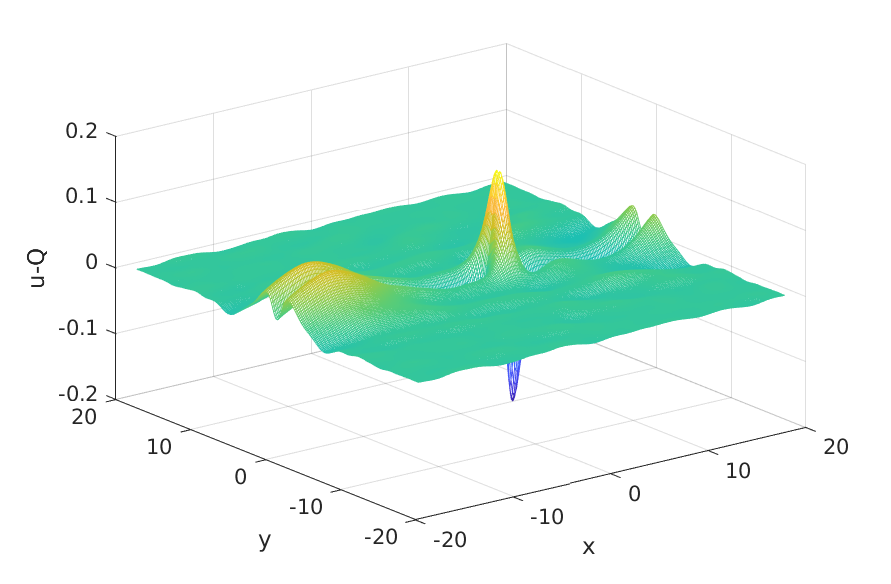} 
 \includegraphics[width=0.49\hsize]{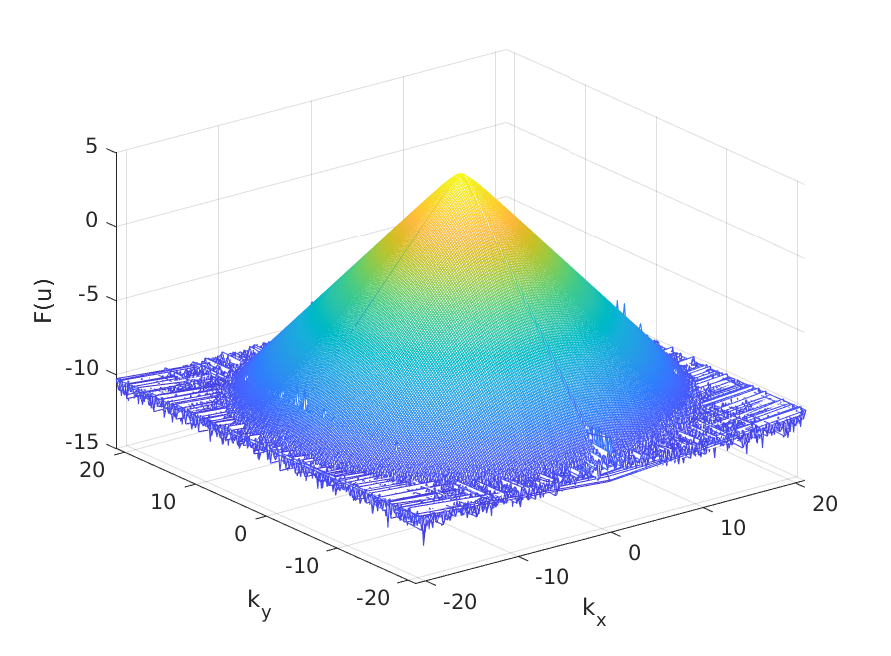}
\caption{ZK solution with $u_0 = 1.1\, Q$: 
solution (projected onto the $z=0$ plane) at $t=12$ (top left), the 
time dependence of the $L^{\infty}$ norm (top right), the difference 
between the solution and a rescaled soliton $Q_c$ at $t=12$ (bottom left), the Fourier coefficients at $t=12$ depending on $k_{x}$ and $k_{y}$ (the solution is symmetric in $k_{y}$ and $k_{z}$, thus, we project onto $k_{z}=0$) (bottom right).}
\label{Stability_11Q}
\end{figure}

We begin with investigating the ZK time evolution for the initial 
data from the part (a) above $u_0 = \lambda \, Q$, $\lambda \gtrsim 
1$. The solution travels in the $x$-direction and increases its 
amplitude till it finds a suitable rescaled soliton $Q_c$ with $c>1$, 
while shedding some radiation in the negative $x$-axis direction. An example of such an  evolution with $\lambda=1.1$ at $t=12$ is shown on the top left of Fig.~\ref{Stability_11Q}. 

Since the perturbation of the original soliton $Q$ is radially symmetric, the solution itself has a symmetry in $y$ and $z$ due to the corresponding symmetry of the ZK equation. In such cases one can plot the solution either in the coordinates $(x, \rho)$, where  
\begin{equation}
	\rho := \sqrt{y^{2}+z^{2}}
	\label{rho},
\end{equation}
or in the plane $z=0$ (which is equivalent), we will typically plot the snapshots of solutions as projections onto the plane $z=0$ as in the top left of Fig.~\ref{Stability_11Q}; nevertheless, the computation is done in the full 3D space without any assumption on symmetry.

The next important feature is tracking the $L^\infty$ norm and its leveling off or saturation in time, if any, which would indicate that the solution is asymptotically approaching the rescaled and shifted version of the soliton $Q_c$ as in \eqref{E:wave}-\eqref{dilation}. The scaling parameter $c$ is defined as 
\begin{equation}\label{E:c}
c= \frac{\| u (t^*) \|_{L^\infty}}{\|Q \|_{L^\infty}},
\end{equation}
where ideally $t^* = \infty$ but for numerical purposes $t^*$ is the 
time when the $L^\infty$ norm levels off (to acceptable precision). The center of the shifted soliton can be easily identified (e.g., in MATLAB) by getting the location coordinates of the local maximum. 
Both the scaling parameter and the shift is used then to measure the difference between the solution and the rescaled and shifted soliton at that time $t^*$ to check if that difference is on the order of the radiation or less (we can then call the time $t^*$ as the `final state' time for numerical investigations).

\begin{figure}[!htb]
 \includegraphics[width=0.51\hsize]{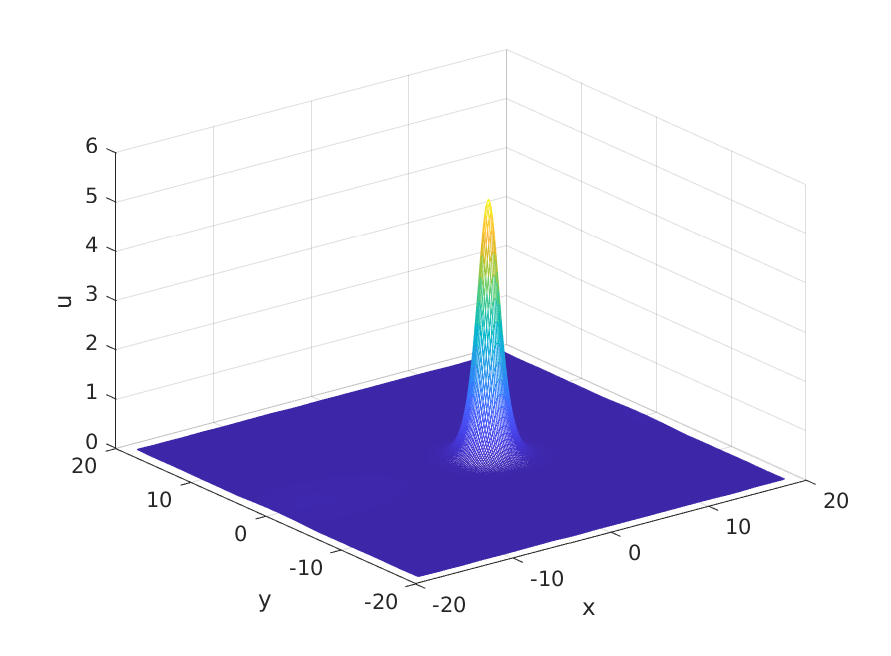}
 \includegraphics[width=0.46\hsize]{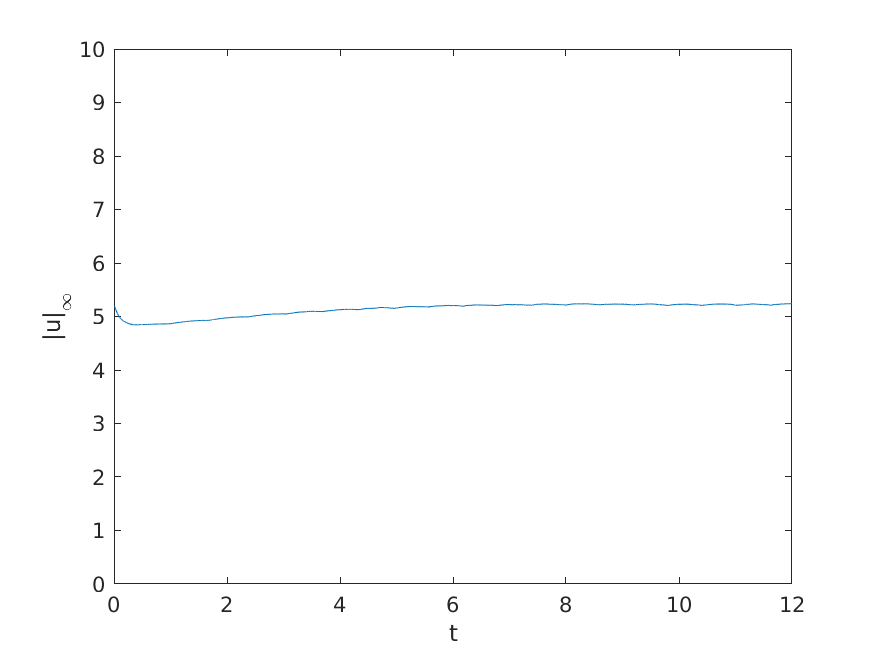}
\includegraphics[width=0.49\hsize]{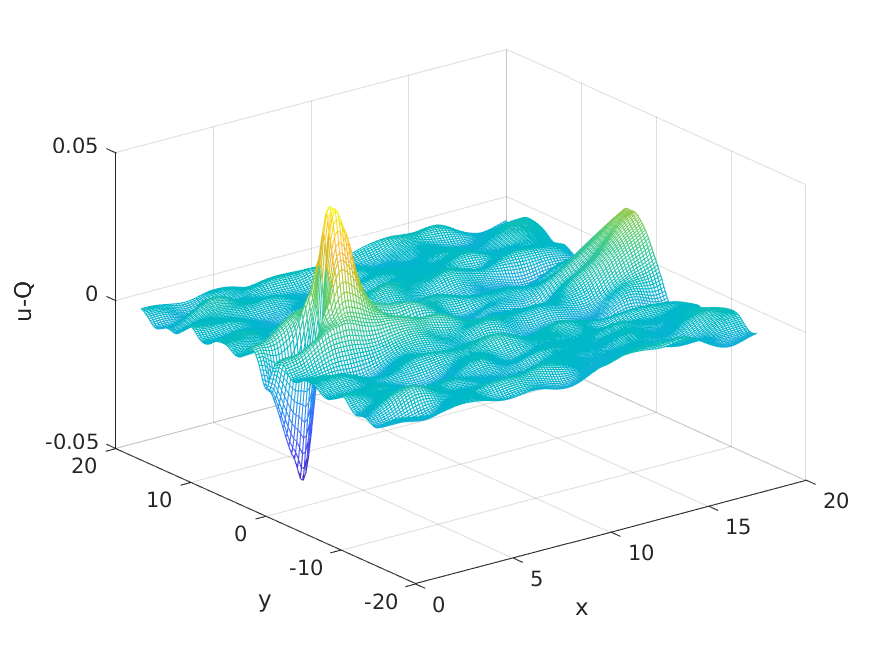}
\includegraphics[width=0.49\hsize]{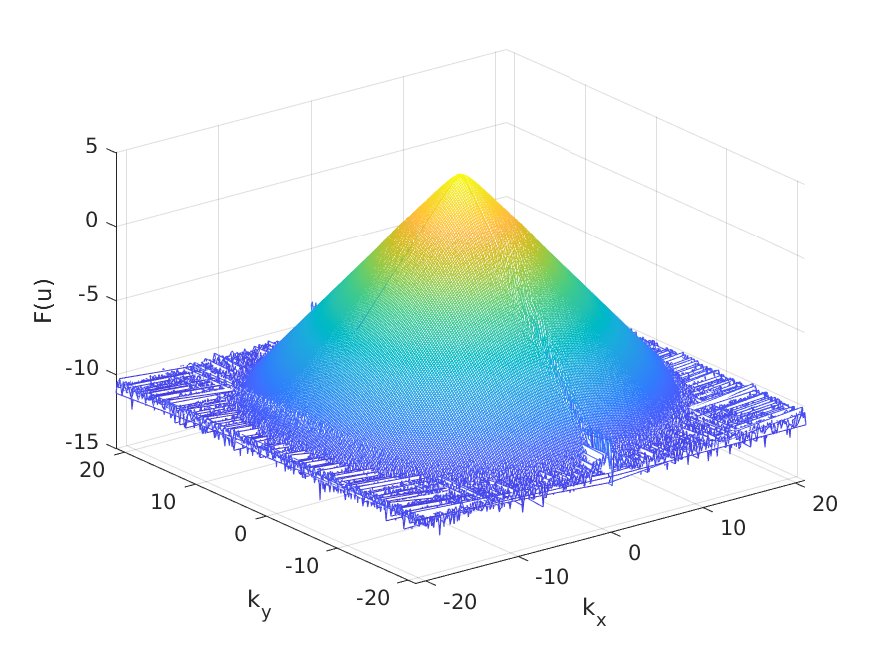}
\caption{ZK solution with asymmetric initial data \eqref{asym}, $\alpha=4$: 
solution (projection onto $z=0$) at $t=12$  (top left), 
the $L^{\infty}$ norm of the solution depending on time (top right),
the difference between this solution and a rescaled soliton $Q_c$ at $t=12$, 
the Fourier coefficients at $t=12$ with respect to $k_{x}$ and  $k_{y}$ 
(the dependence on $k_{z}$ is very similar to $k_y$, and thus, not shown 
for the ease of presentation).}
\label{Stability_Assym}
\end{figure}

The $L^{\infty}$ norm of the solution with data $u_0=1.1\,Q$ is plotted in the top right of Fig.~\ref{Stability_11Q}, saturating at a slightly higher amplitude than the initial value. This seems to be typical: an initial condition that is a multiple of the ground state with higher amplitude than the original soliton ($\lambda > 1$) leads to a rescaled version of $Q_c$, defined in \eqref{dilation}, with $c>1$, and thus, faster moving in the $x$-direction than the original soliton with speed $c=1$. This can also be observed on the top left plot of Fig.~\ref{Stability_11Q}, since the location of the maximum is slightly shifted away from the origin in the positive $x$-direction (in other words, its maximum does not stay at the origin in the co-moving frame \eqref{ZKcom} with $v_{x}=1$ and shifts to the right). 
The difference between the solution at $t=12$ and a rescaled soliton $Q_c$ is plotted on the bottom left of Fig.~\ref{Stability_11Q}. 
One can notice that the difference is on the order of $10^{-1}$, thus, confirming the asymptotic shape of the solution as a rescaled and shifted soliton (plus radiation).  
There is some radiation propagating into the negative $x$-direction. 
Since we work on a periodic domain, the radiation cannot escape to (negative) infinity, and it reenters the domain on the opposite side, which leads to a slightly noisy background (as we discussed in Remark \ref{remperiodic}). Therefore, it might be easier to see the radiation on the right side of either top left or bottom left plots of Fig.~\ref{Stability_Assym}. 
The Fourier coefficients at $t=12$ are shown in 
Fig.~\ref{Stability_11Q} on the bottom right. One can observe that an accuracy much better than plotting accuracy is achieved. 
\smallskip

\underline{\bf Case (a): multiple of a soliton with $\lambda \lesssim 1$.}
\smallskip

Solutions of the ZK equation \eqref{ZK} with initial data from the part (a) $u_0=\lambda \, Q$ with $\lambda \lesssim 1$ also asymptotically approach solitons, though with amplitudes smaller than the initial data (and hence, smaller than $Q$ itself).
\begin{figure}[!htb]
\includegraphics[width=0.51\hsize]{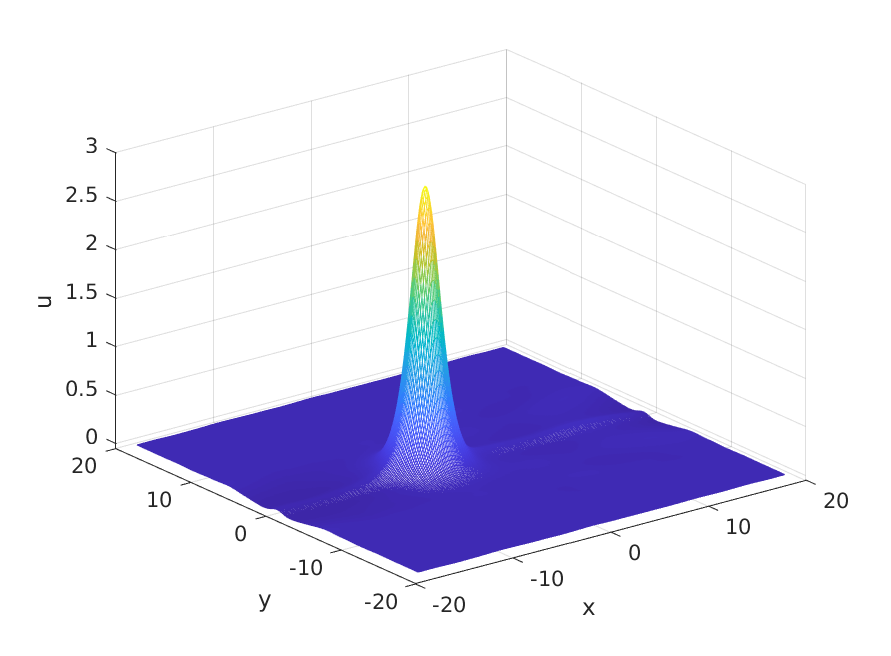}
\includegraphics[width=0.46\hsize]{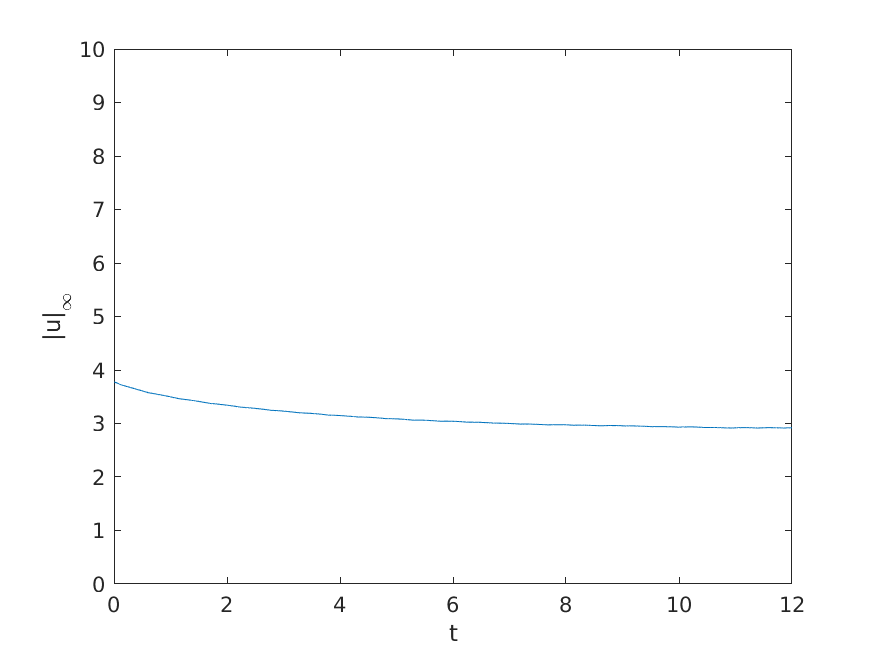} 
\includegraphics[width=0.49\hsize]{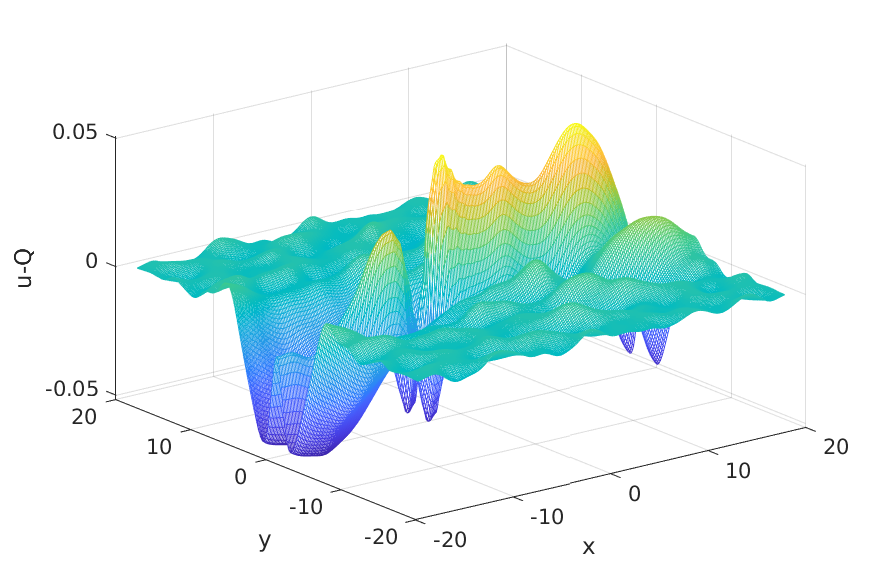}
\includegraphics[width=0.49\hsize]{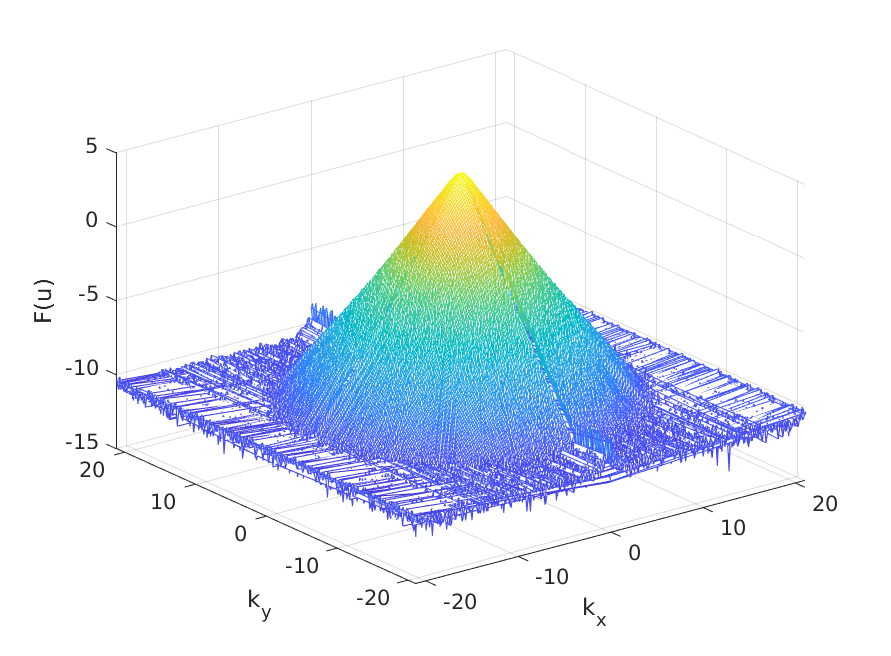}
\caption{ZK solution with $u_0=0.9\,Q$: solution (projected onto $z=0$) at $t=12$ (top left), 
the $L^{\infty}$ norm of the solution depending on time (top right), 
the difference between the solution and a rescaled $Q_c$ at $t=12$ (bottom left), 
the Fourier coefficients at $t=12$ depending on $k_{x}$ and $k_{y}$ 
(projected onto $k_{z}=0$).}
\label{Stability_09Q}
\end{figure}

An example with $\lambda = 0.9$ is shown in Fig.~\ref{Stability_09Q}. First, note that the peak moves to the left of the origin in the co-moving frame (that is, in the negative $x$-direction), indicating that the resulting soliton has a slower speed and smaller amplitude than the soliton $Q$ itself (and also, smaller than the initial amplitude). 
This is also confirmed on the top right of Fig.~\ref{Stability_09Q}, in the decay and then leveling off of the $L^{\infty}$ norm of this solution.

On the bottom left of Fig.~\ref{Stability_09Q} we show the difference between the solution at $t=12$ and a rescaled soliton $Q_c$, where $c$ is computed from \eqref{E:c}; the difference is on the order of $10^{-2}$.  
The Fourier coefficients, plotted on the bottom right of the same figure, indicate the 
numerical resolution of the solution. 
\smallskip

\underline{\bf Case (b): asymmetric perturbations of a soliton.}
\smallskip

The previous perturbations of the soliton-type initial data had the same 
symmetry as the unperturbed soliton, and thus, led to a solution only 
depending on $\rho$ defined in \eqref{rho}; we now look at the non-symmetric (in $y$ and $z$) localized perturbations. For example, one can consider
\begin{equation}\label{asym}
u(x,y,z,0)=Q(x,y,z) +e^{-(x^2+y^2+ \alpha \, z^2)}, ~\alpha \neq 1, \alpha > 0.
\end{equation}

An example of the ZK solution with the initial data \eqref{asym} and 
$\alpha = 4$ is shown in Fig.~\ref{Stability_Assym}. A snapshot of 
the solution at $t=12$ is on the top left, and the time dependence of the $L^{\infty}$ norm is plotted on the top right of 
Fig.~\ref{Stability_Assym}. Both show that this solution converges to a rescaled soliton with 
the amplitude slightly higher than the initial one. Since the amplitude is higher than the original, the rescaled soliton, to which the solution is converging asymptotically, is moving slightly faster than the original unperturbed soliton: on the top left plot in Fig. \ref{Stability_Assym} the location of the peak is slightly shifted in the positive $x$-direction; this is also confirmed by the difference of this solution and a rescaled soliton $Q_c$ on the bottom left of the same figure. 
The Fourier coefficients on the bottom right 
show that the solution is numerically well resolved.

\begin{figure}[!htb]
\includegraphics[width=0.49\hsize]{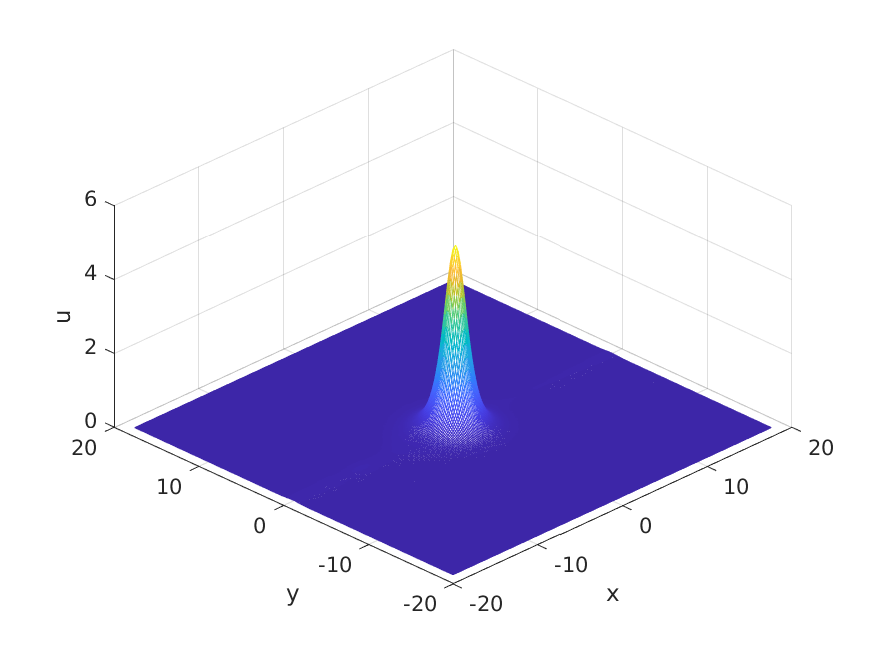}
 \includegraphics[width=0.49\hsize]{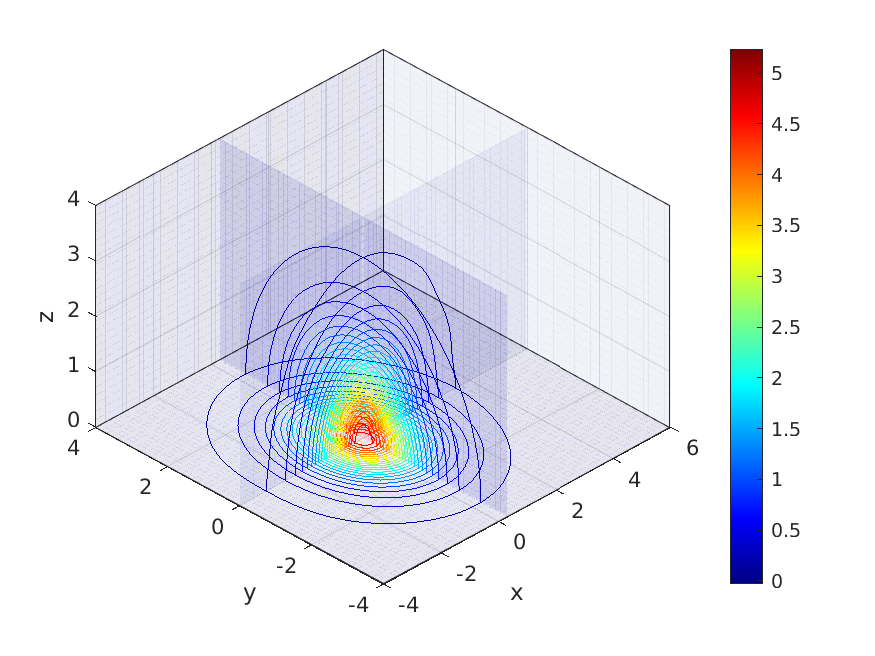}\\
\includegraphics[width=0.49\hsize]{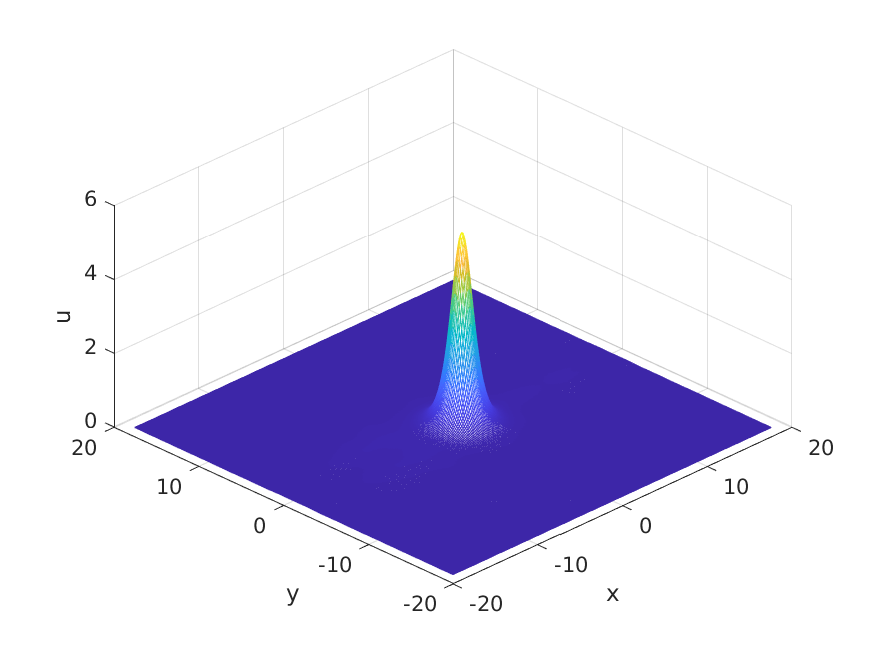}
 \includegraphics[width=0.49\hsize]{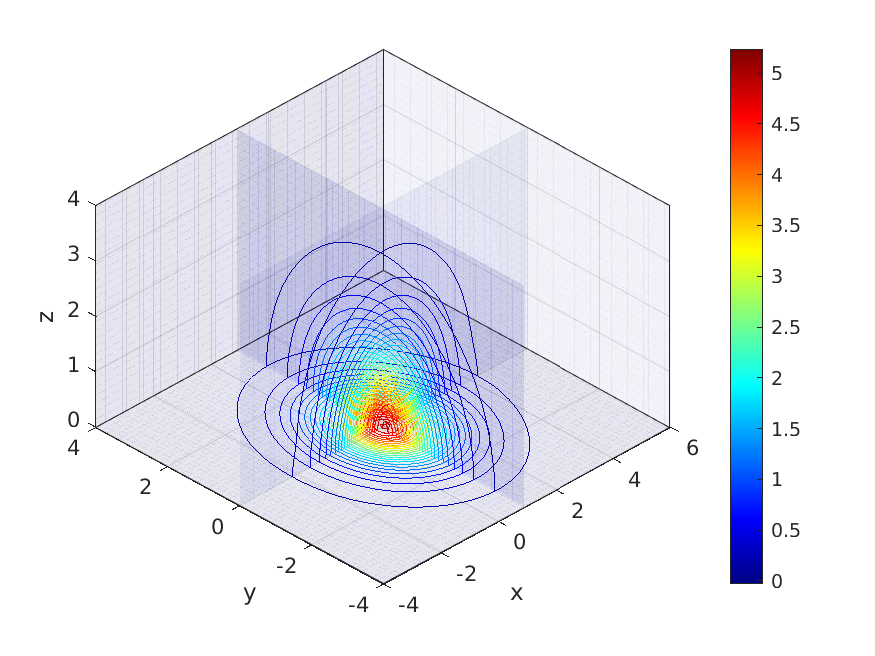}\\ 
 \includegraphics[width=0.49\hsize]{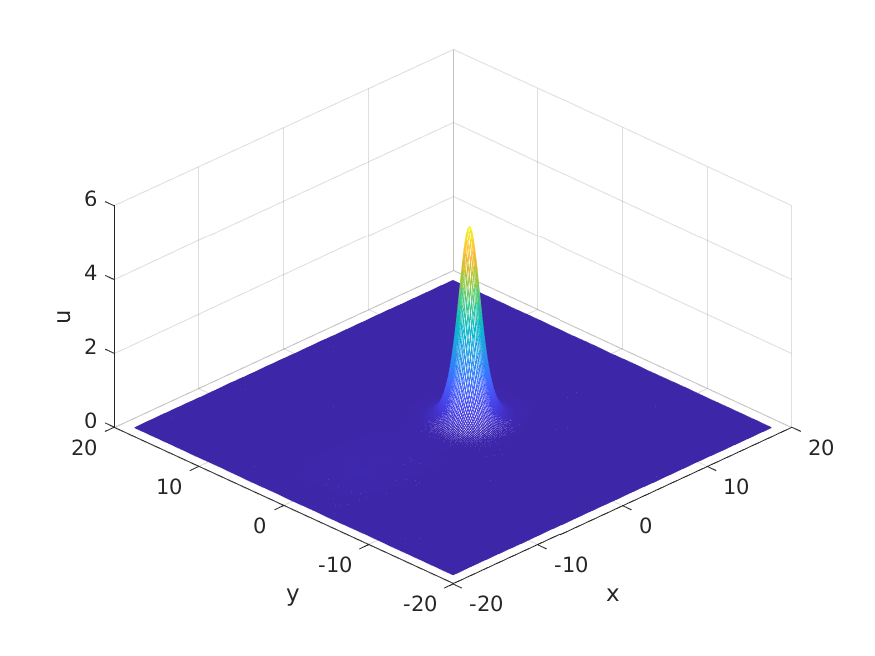}
 \includegraphics[width=0.49\hsize]{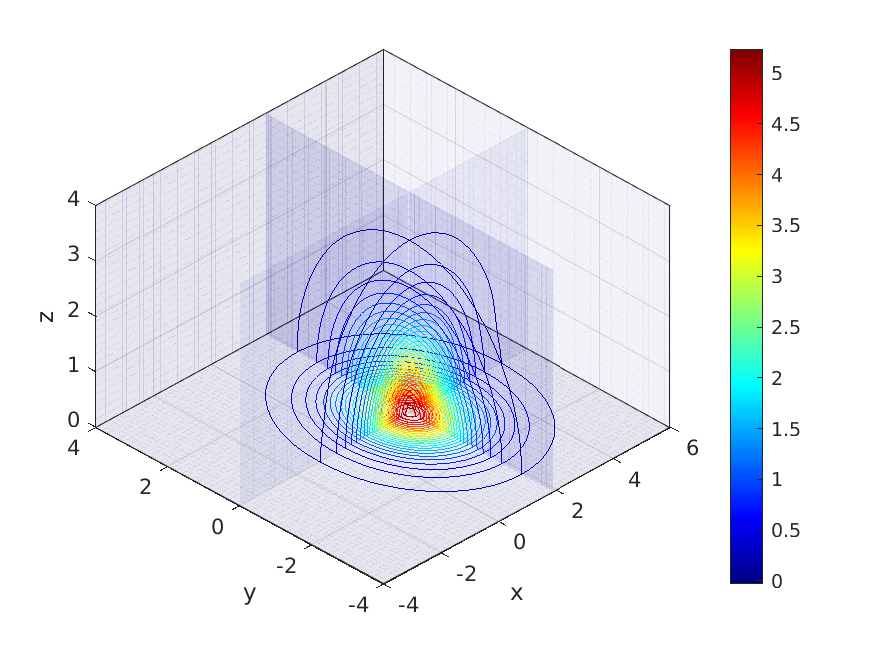}\\%
\caption{Snapshots of ZK solution with $u_0$ as in \eqref{asym}, $\alpha=4$, 
at $t=1.2$, $4.8$, $8.4$. Left: solutions projected onto $z=0$. 
Right: 3D isocurves on the slices of the coordinate planes.}
\label{Stability_Assym-3d}
\end{figure}

For better understanding of this solution, 
we provide a 3D visualization in Fig.~\ref{Stability_Assym-3d} 
by showing contour plots on each of the coordinate plane slices through the solution. 

One may notice that the initial data in this case are not spherically symmetric, nevertheless,  the time evolution leads to a more symmetric configuration. For that consider the top right 3D isocurves in Fig.~\ref{Stability_Assym-3d} at the points where $y=0$ and note that those curves are more flat compared with the corresponding isocurves at $y=0$ in the very bottom plot (of course, the solution has moved in a positive $x$-direction); from ellipses the isocurves in the $xy$-plane become more circular, indicating a more symmetric solution.

\begin{figure}[!htb]
 \includegraphics[width=0.49\hsize]{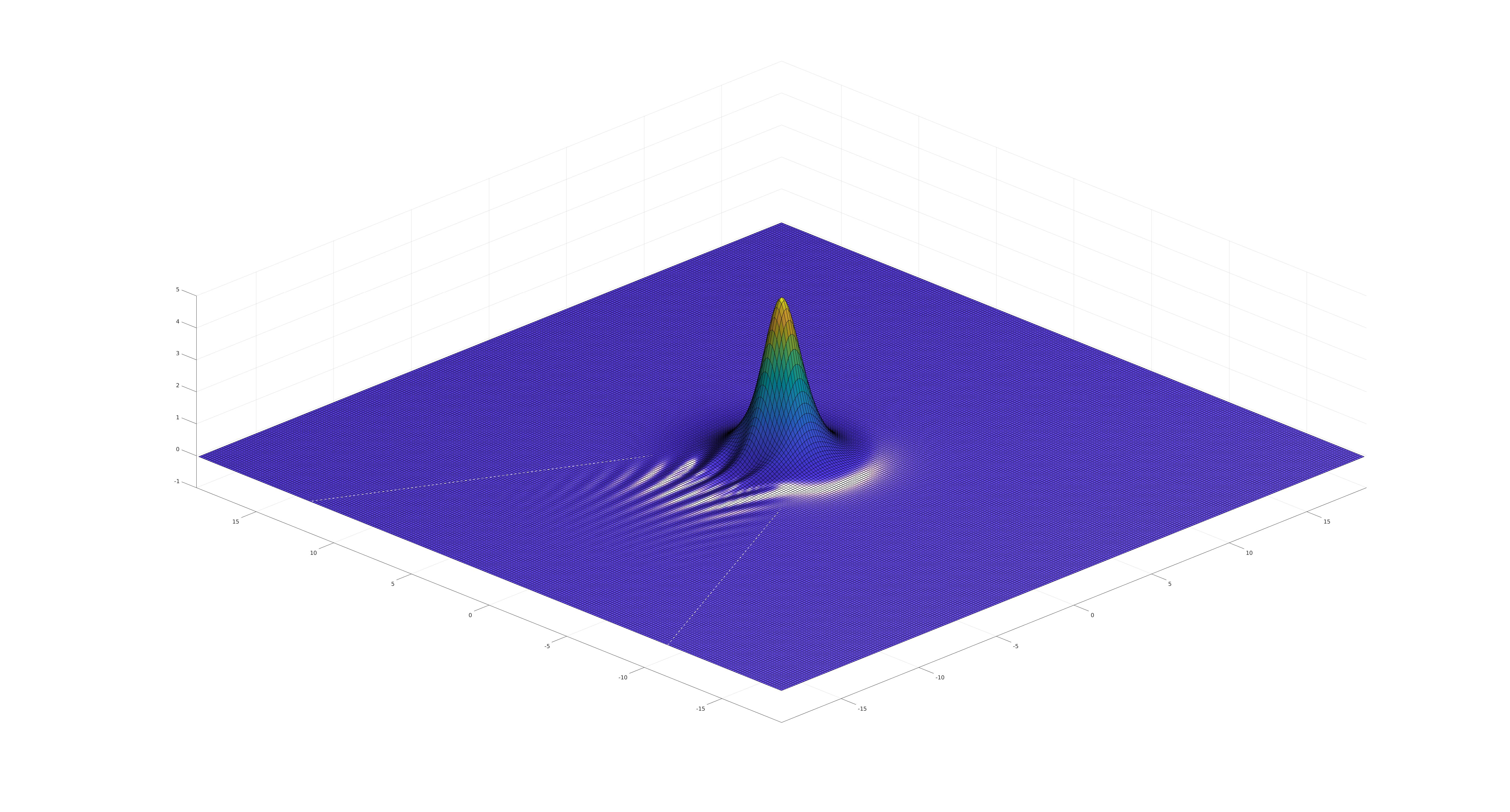}
 \includegraphics[width=0.49\hsize]{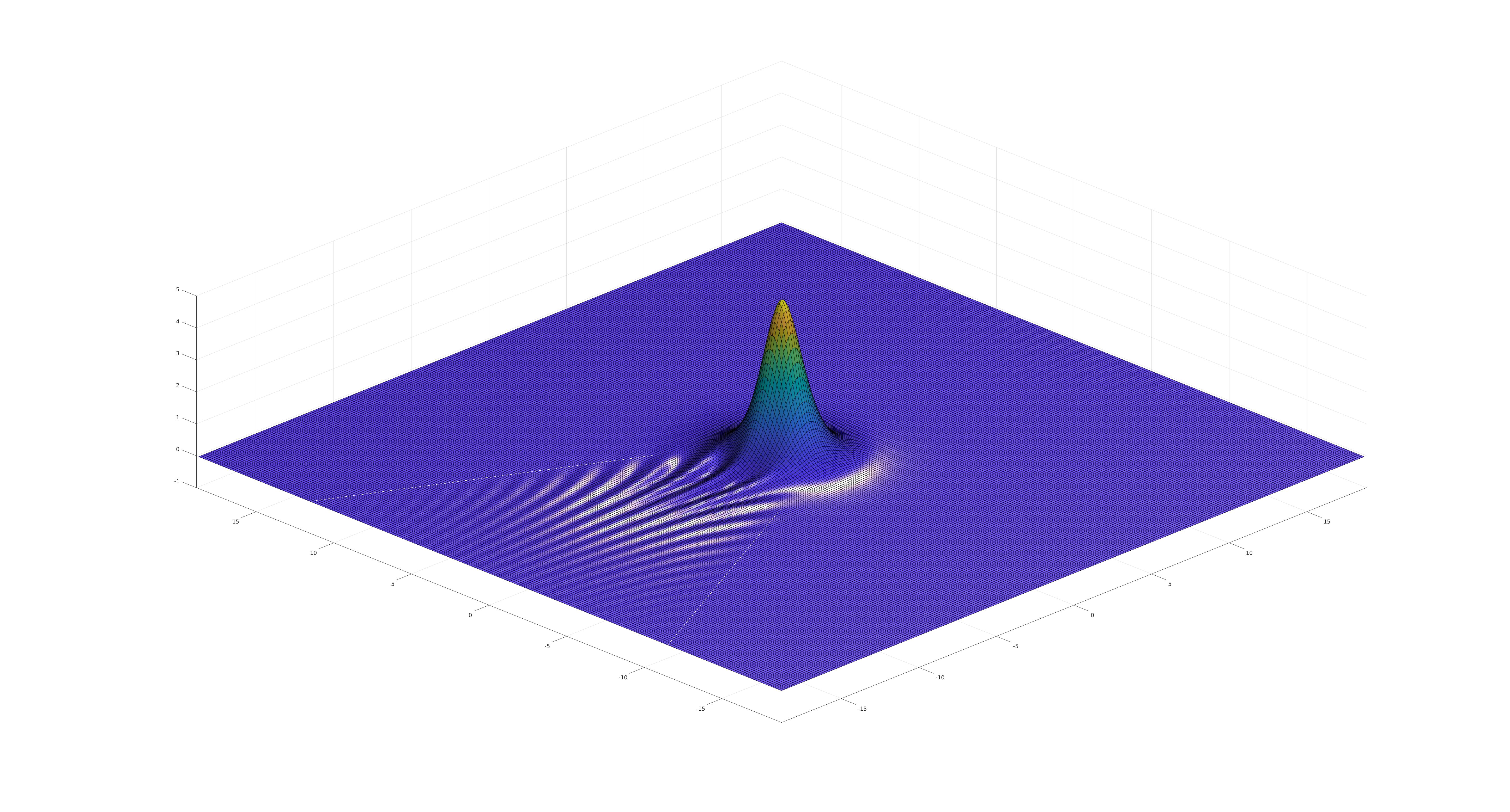} 
 \includegraphics[width=0.49\hsize]{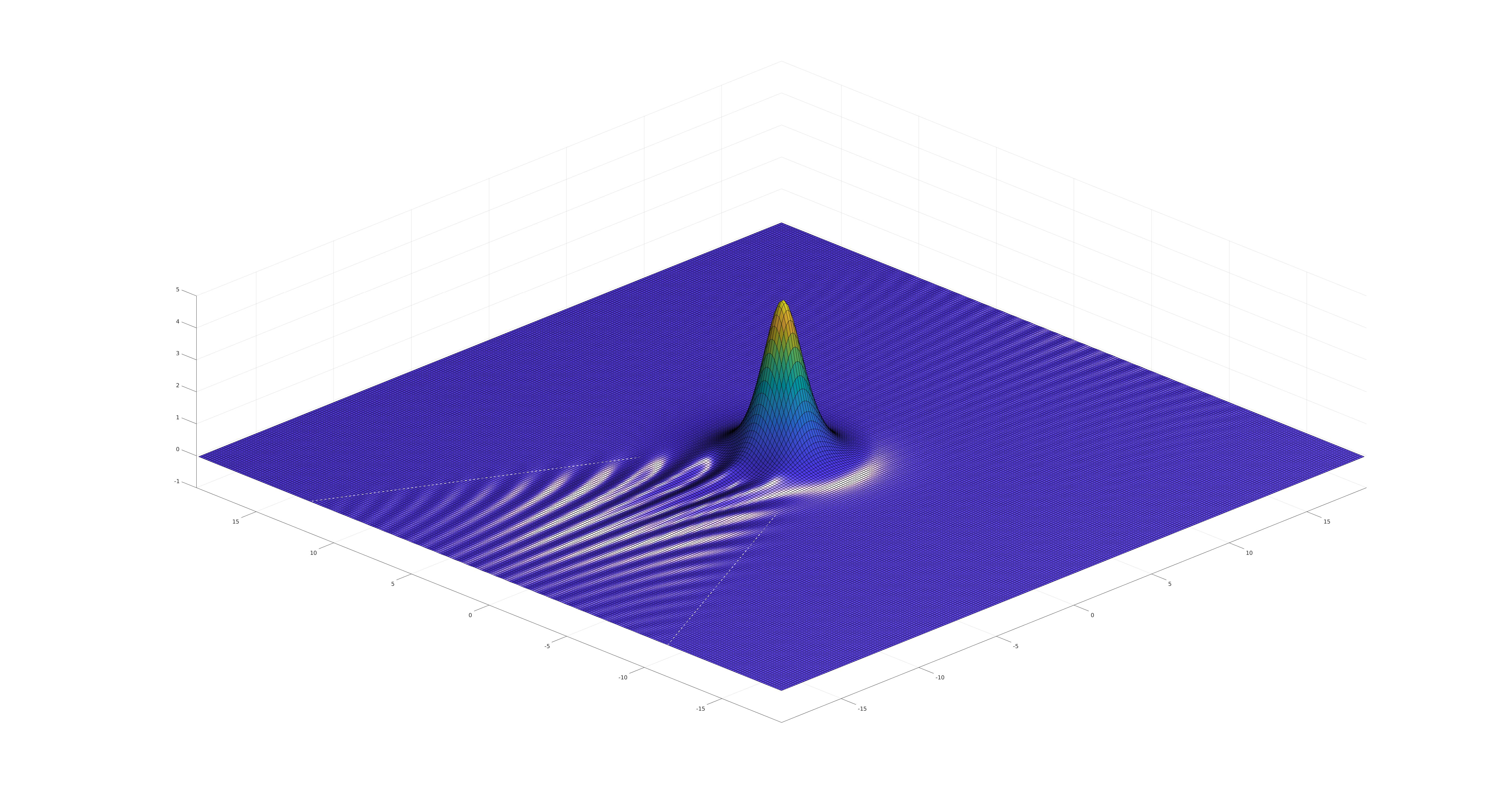}
 \includegraphics[width=0.49\hsize]{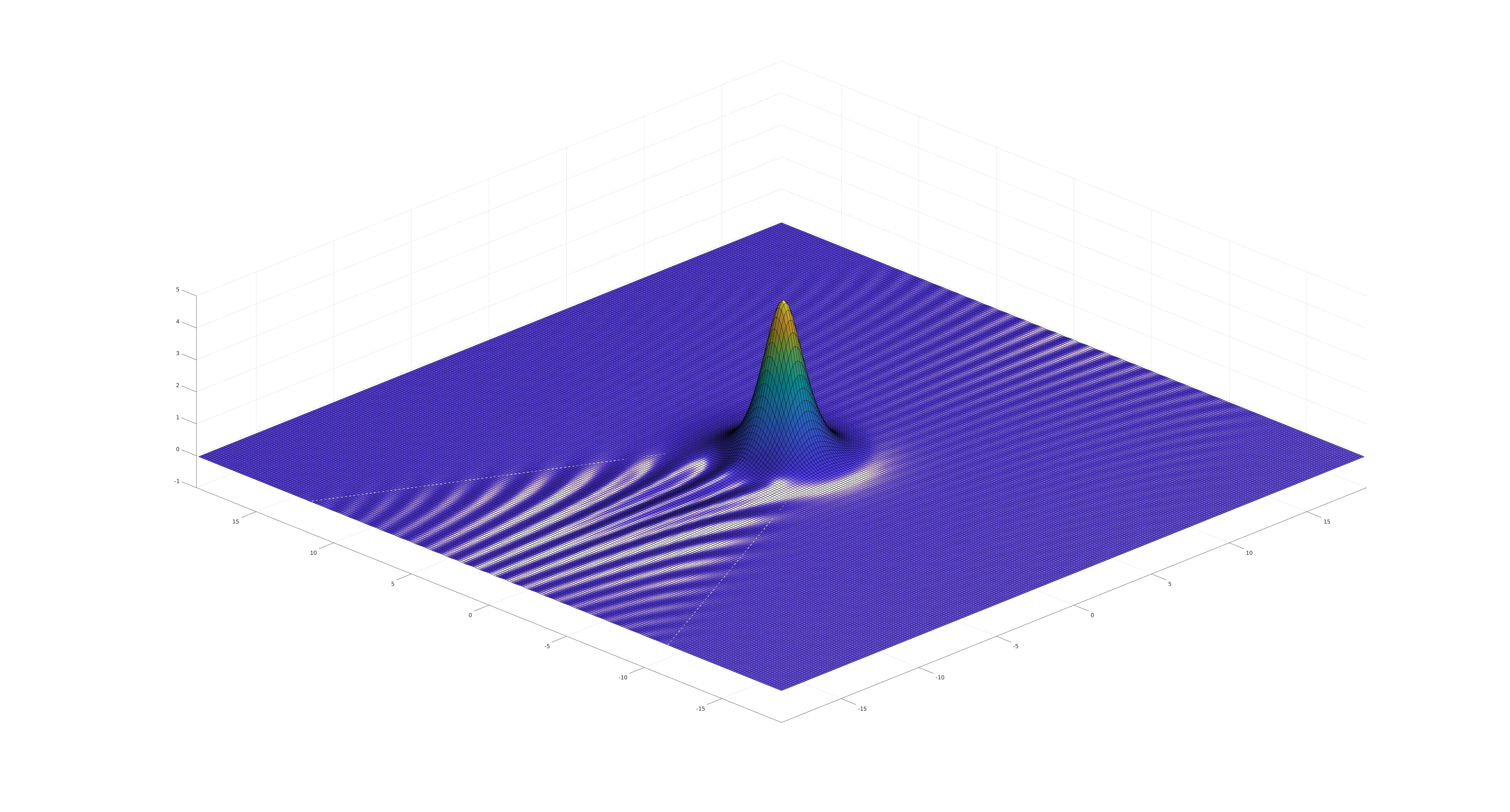}
\caption{Outgoing radiation in the asymmetric deformation of the initially perturbed soliton $u_0=Q + e^{-(x^2+y^2+4z^2)}$ at $t = 0.12, 0.24, 0.36, 0.48$. White lines: the expected wedge of the radiation front \eqref{E:angle} forming a total angle of $\pi/3$.}
\label{F:Qasym-radiation}
\end{figure}

Finally, we discuss the formation of radiation in this example. In 
Fig.~\ref{F:Qasym-radiation} we show the development of outgoing 
radiation in time (again projected onto $z=0$). The radiation escapes to the left conic wedge of the moving rightward solution at an angle of $30^0$ with the negative $x$-axis (for a total opening of $60^0$). This is in agreement with the asymptotic stability result in \cite{FHRY2020}.  

\subsection{Soliton Resolution} 
We next study initial data with monotone decay, vanishing at infinity and either a single maximum or a continuous interval of the same maximum. 
The goal here is to investigate how the ZK evolution of such data resolves into solitons and radiation. 
Before we consider non-soliton-like initial data, we mention that if 
$u_0=\lambda \, Q$ with either $\lambda \gg 1$ or $\lambda \ll 1$, 
then we typically observe that the $L^\infty$ norm saturates at a 
certain level (much faster for larger $\lambda$ and significantly 
slower for smaller $\lambda$). It is plausible to think that such 
data always asymptotically approach a rescaled and shifted soliton (since the solitons can be of an arbitrary size in the ZK equation). The time of how fast the solution resolves into the soliton and radiation is inversely proportional to the initial amplitude of the data.

We next consider initial data, which is different from a soliton, its multiple or a small perturbation. The decay rates that we consider are either faster than the exponential decay of the soliton (for example, a Gaussian) or a polynomial decay. Namely, we examine the following cases of initial data: 

\begin{itemize}
\item[(a)] Gaussian, 
\item[(b)] flattened Gaussian,
\item[(c)] wall-type, 
\item[(d)] fast algebraic decay, referred here to as `super-Lorenzian'.
\end{itemize}

\underline{\bf Case (a): Gaussian initial data.}
\smallskip

\begin{figure}[!htb]
\includegraphics[width=0.51\hsize]{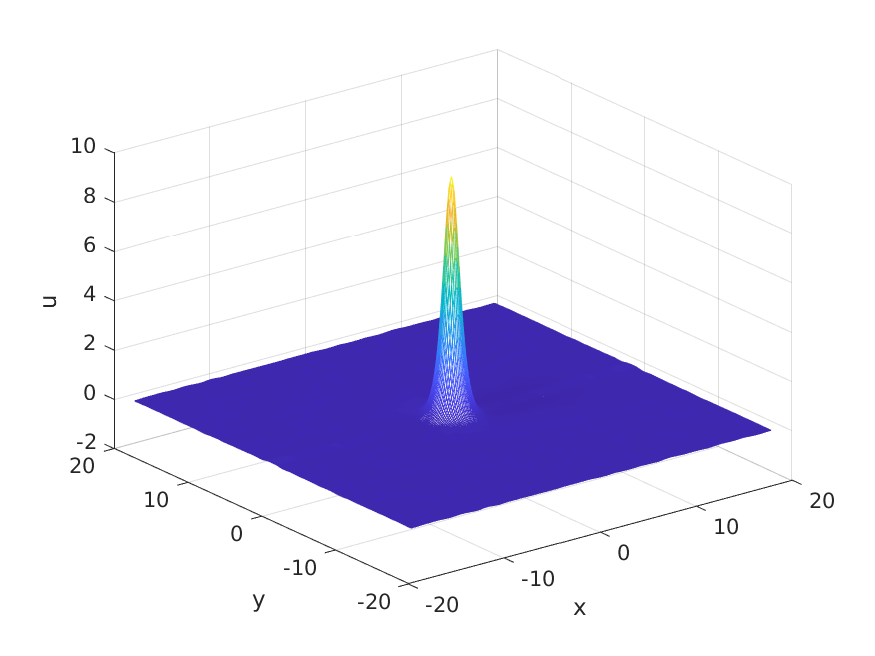}
\includegraphics[width=0.46\hsize]{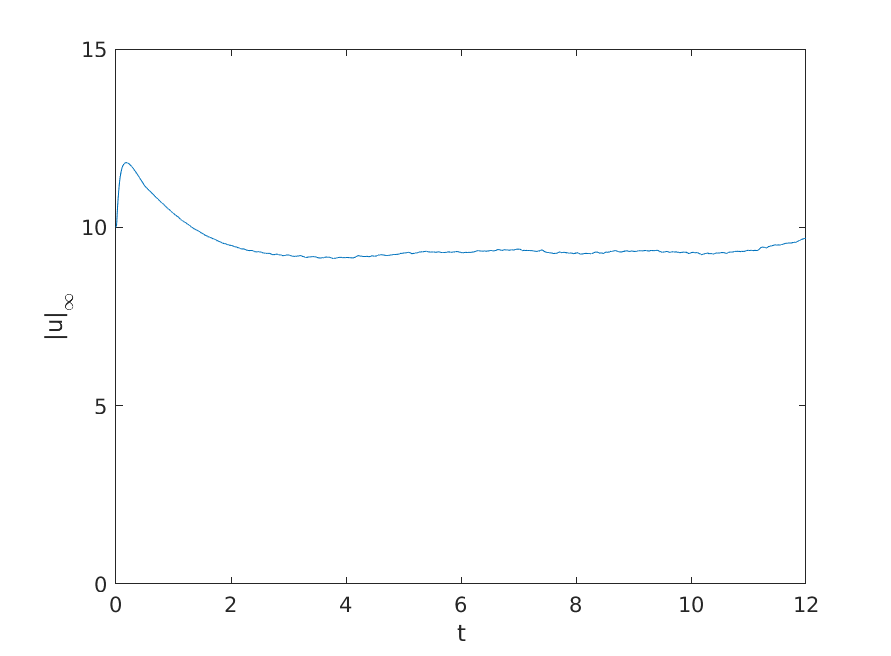}
\includegraphics[width=0.49\hsize]{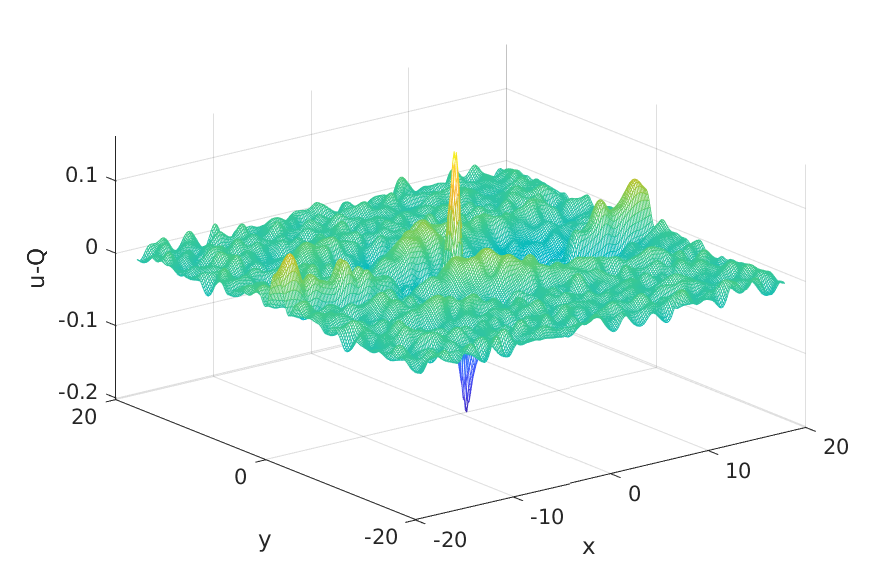}
\includegraphics[width=0.49\hsize]{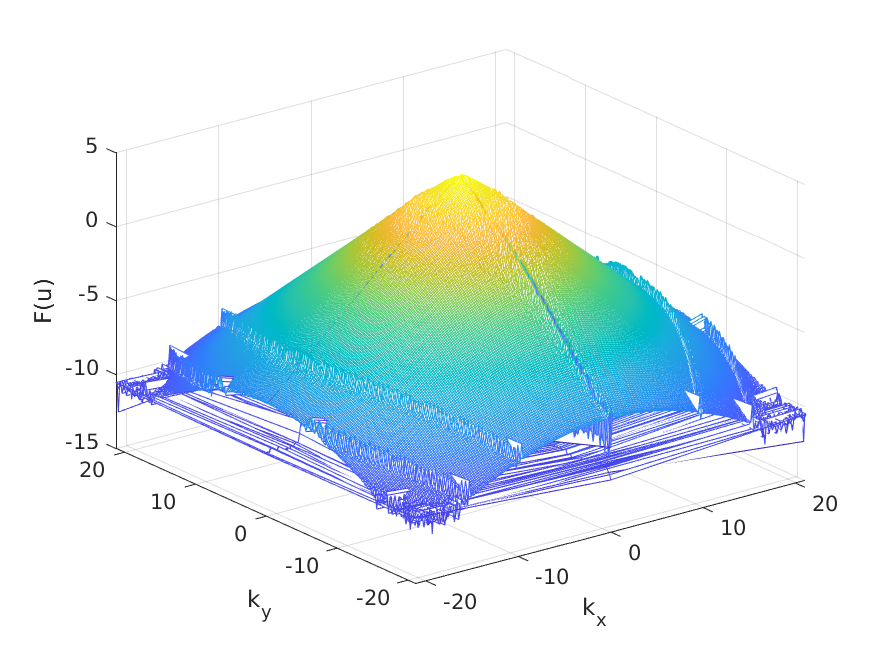}
\caption{ZK solution with Gaussian initial data $u_0=10 \, e^{-(x^2+y^2+z^2)}$:
solution (projected onto $z=0$) at $t = 12$ (top left), 
the $L^{\infty}$ norm depending on time (top right), the 
difference between the solution and a rescaled soliton $Q_c$ (bottom left), 
the Fourier coefficients at $t=12$ (bottom right).}
\label{Gauss_final}
\end{figure}

\begin{figure}[!htb]
\includegraphics[width=0.49\hsize]{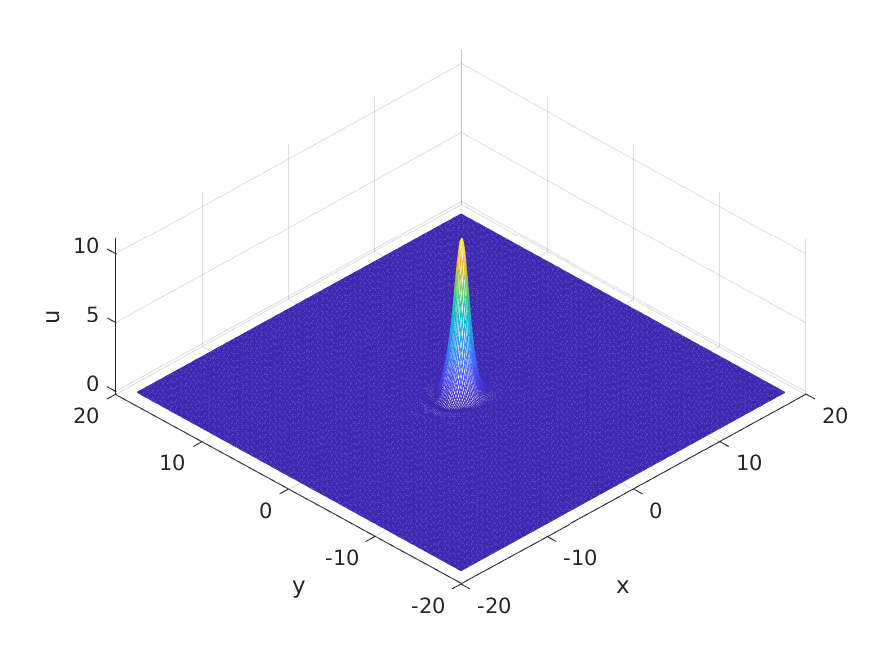}
\includegraphics[width=0.49\hsize]{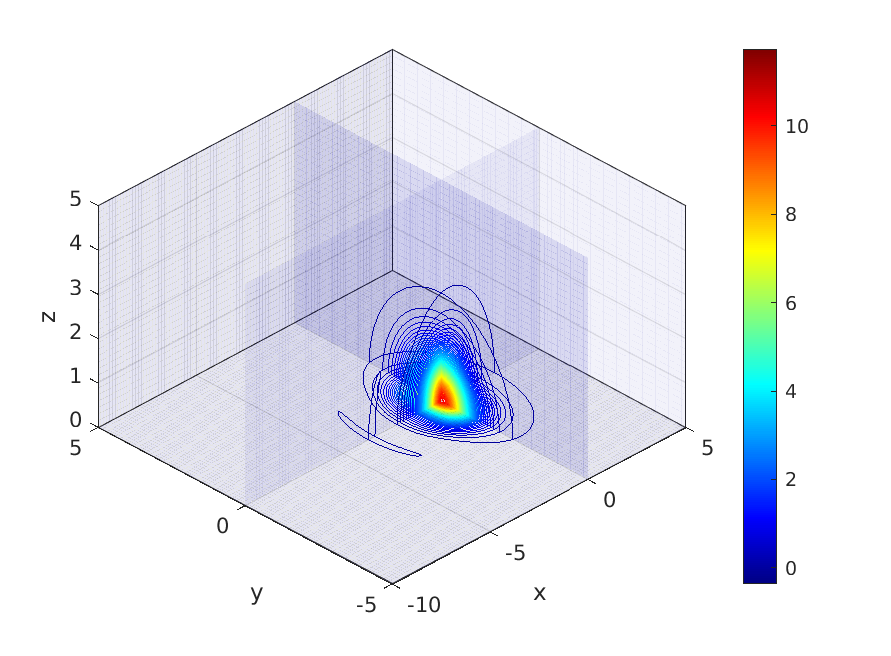}\\    
\includegraphics[width=0.49\hsize]{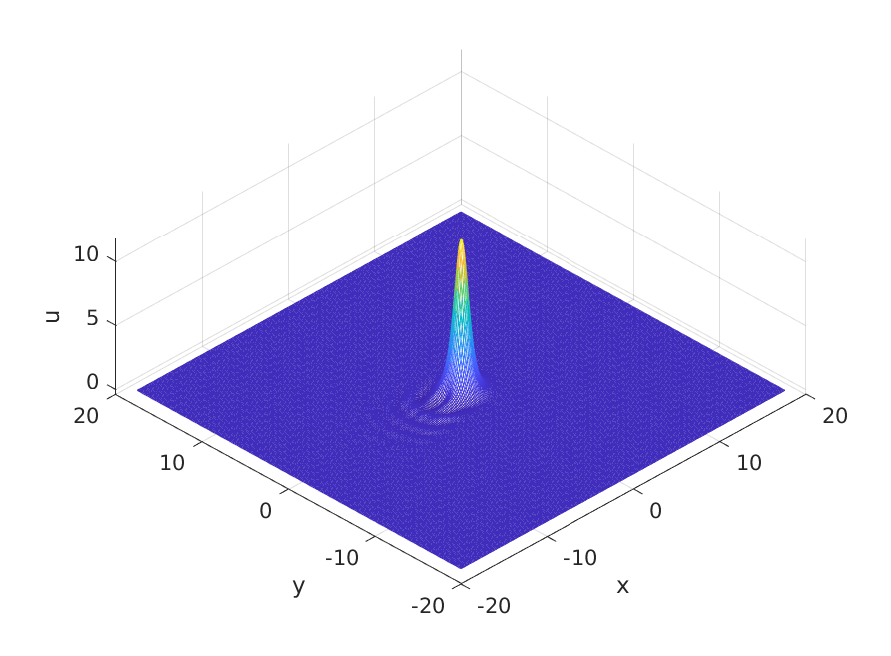}
\includegraphics[width=0.49\hsize]{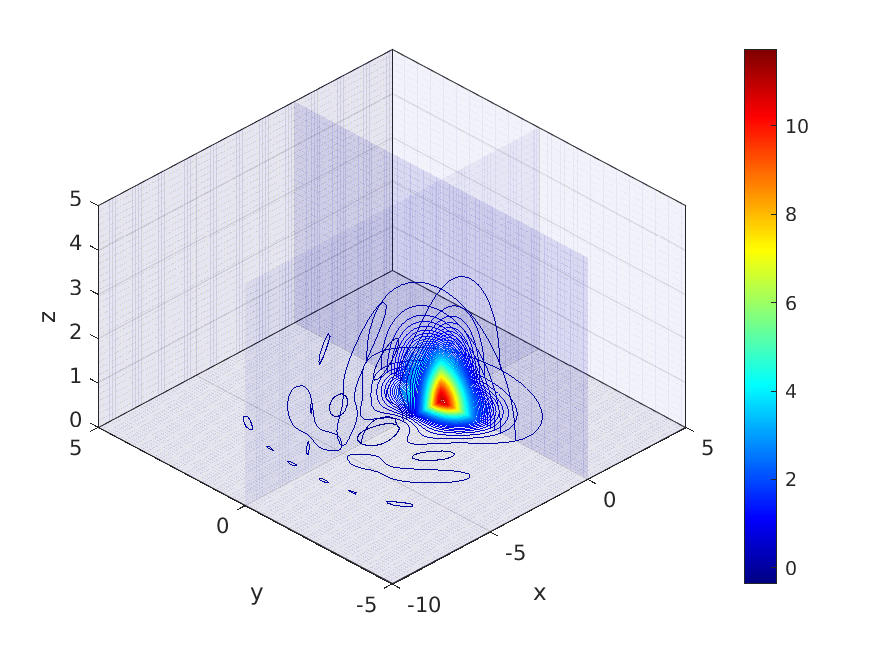}\\  
\includegraphics[width=0.49\hsize]{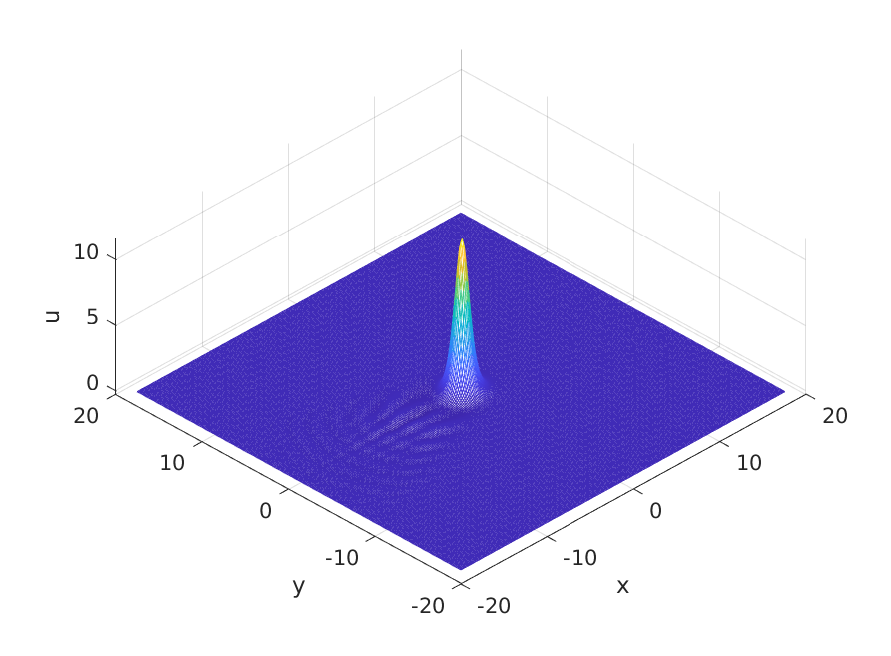}
\includegraphics[width=0.49\hsize]{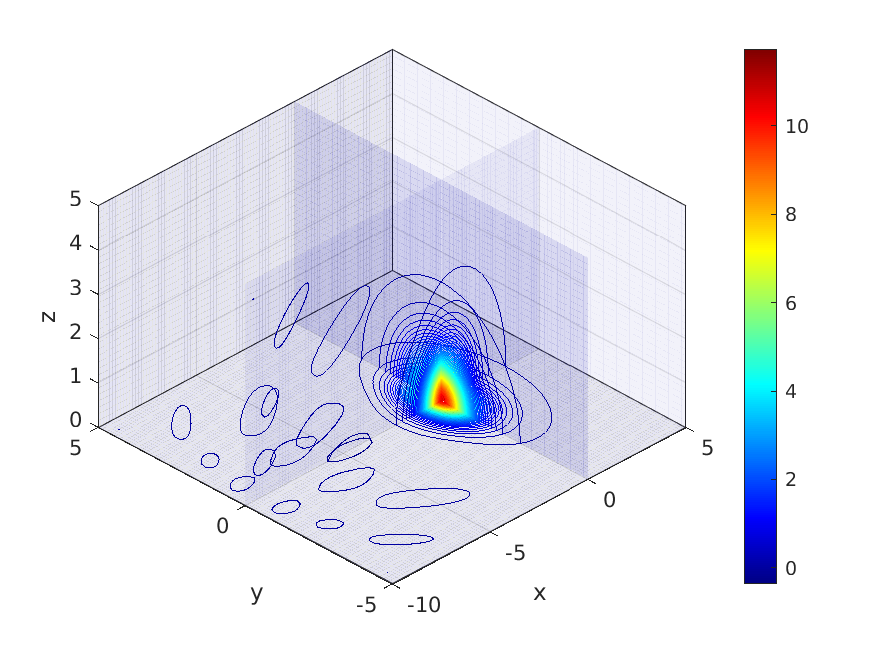}\\  
\caption{Snapshots of the ZK solution with Gaussian 
initial data $u_0 = 10 \, e^{-(x^2+y^2+z^2)}$ at $t = 0.05$, $0.15$, $0.35$. Left: two dimensional projections onto $z=0$. Right: 
3D isocurves on the slices of the coordinate planes.}
\label{Gaussslices}
\end{figure}

\begin{figure}[!htb]
\includegraphics[width=0.49\hsize]{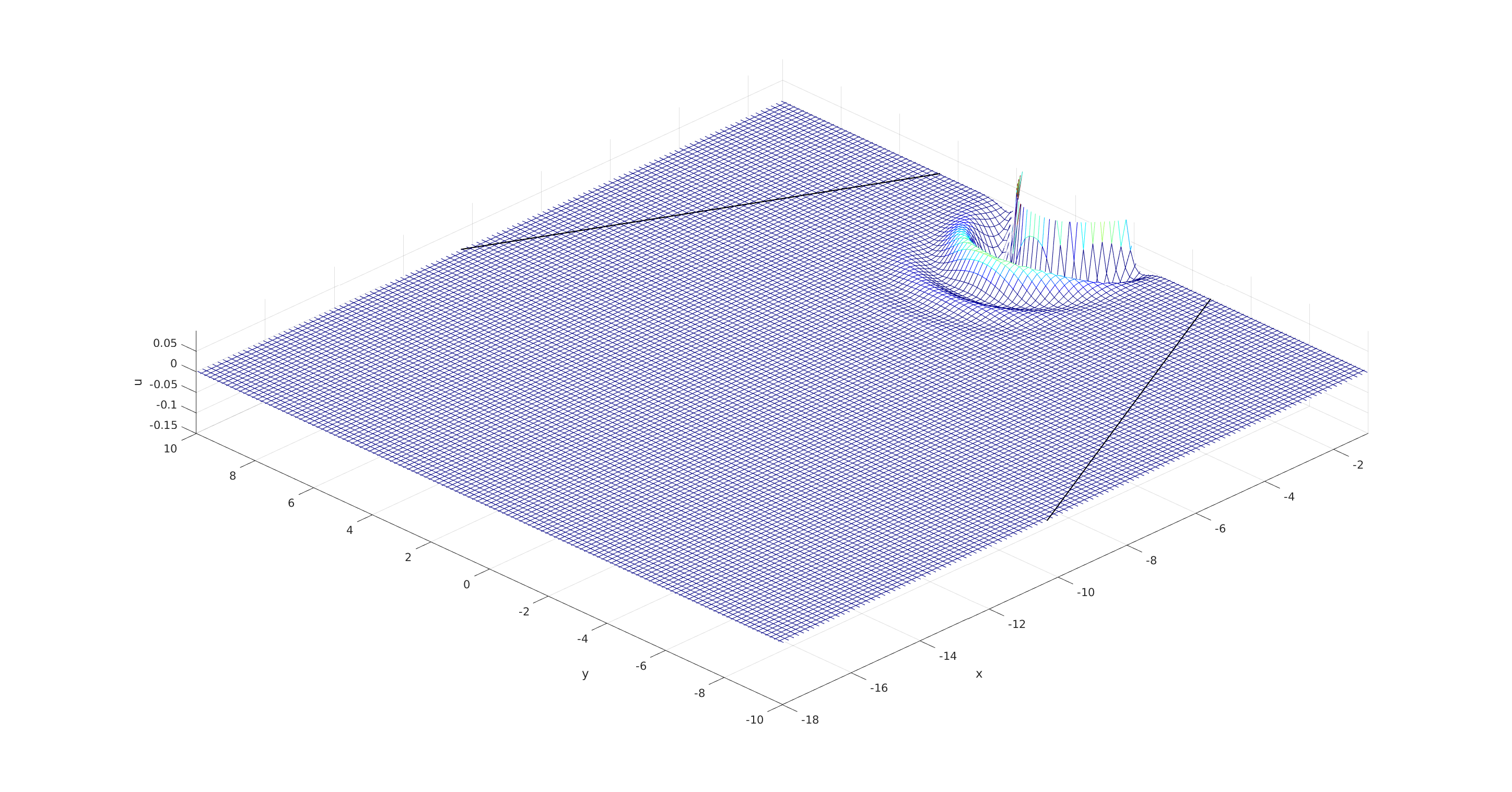}
\includegraphics[width=0.49\hsize]{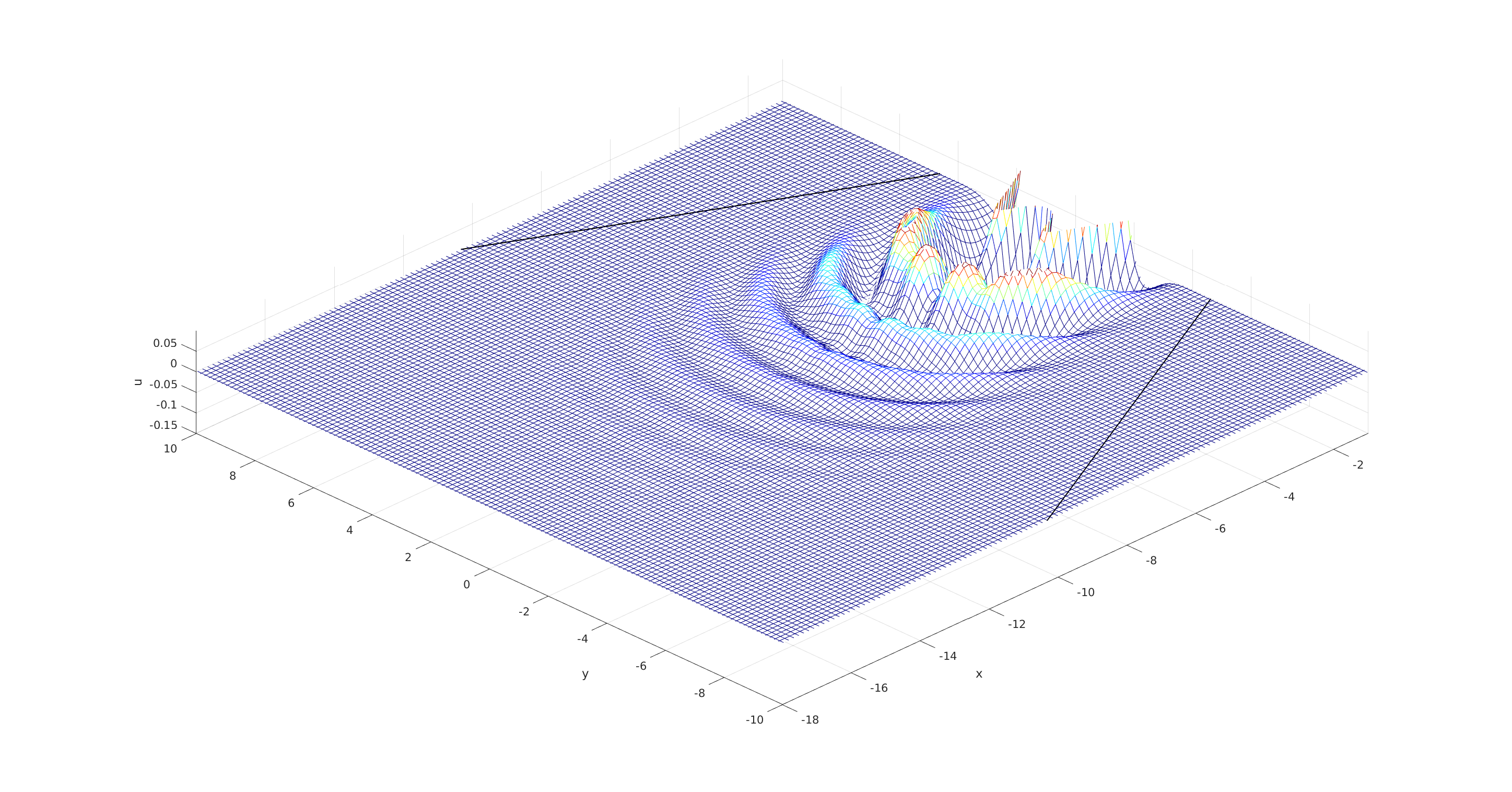} \\ 
\includegraphics[width=0.49\hsize]{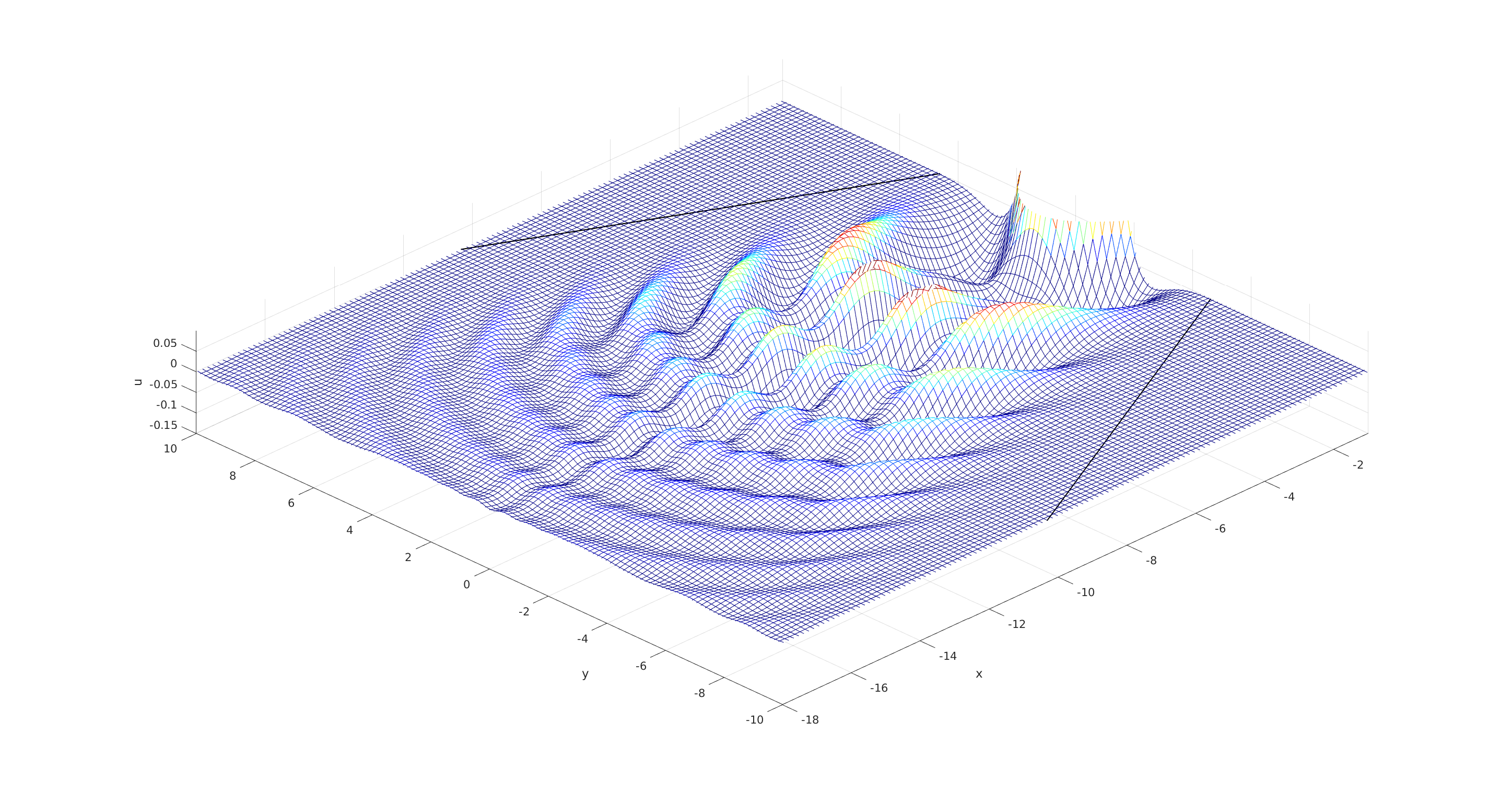}
\includegraphics[width=0.49\hsize]{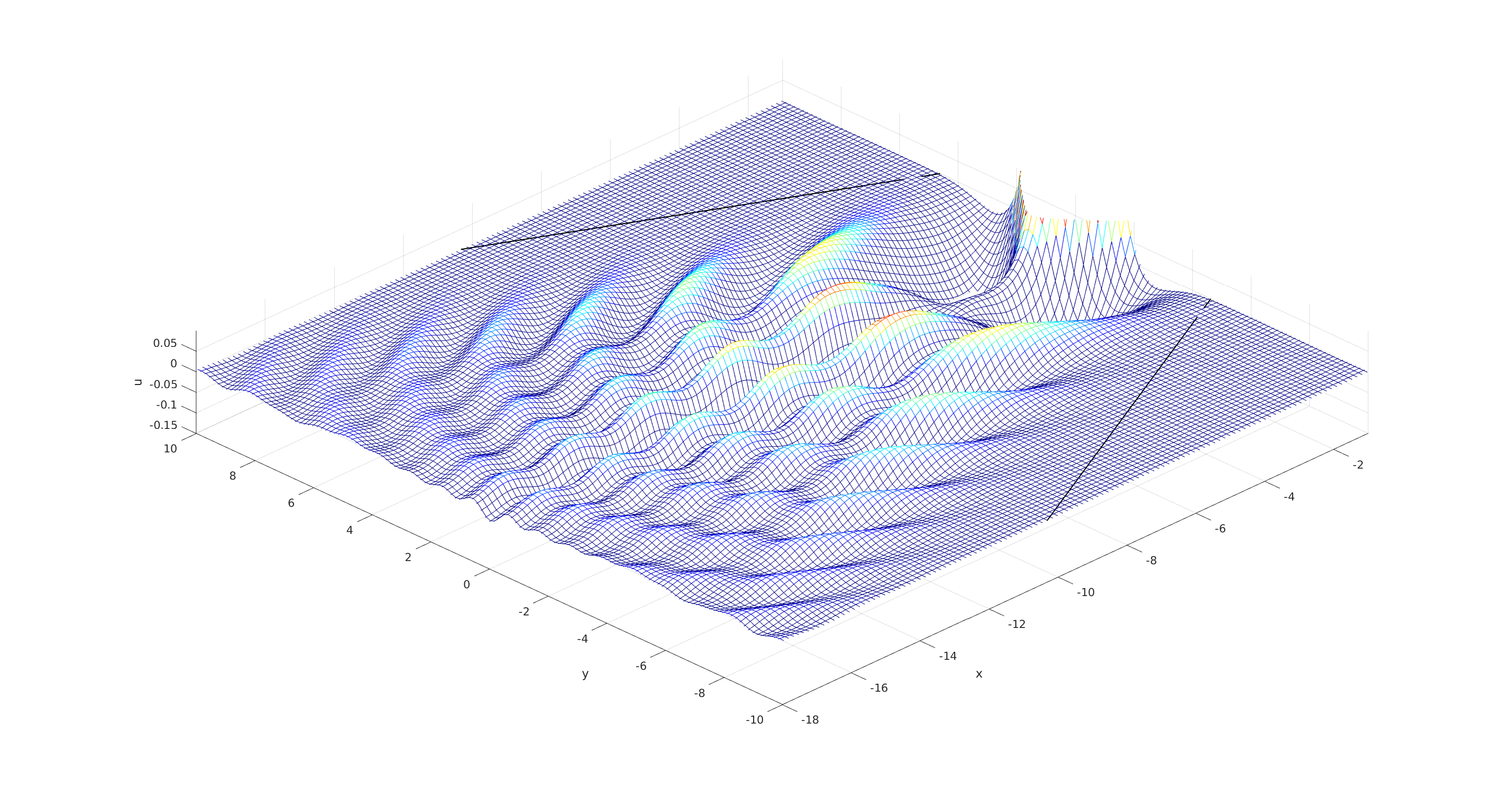}   
\caption{Detail view of the radiation developed in the ZK solution with Gaussian 
initial data $u_0 = 10 \, e^{-(x^2+y^2+z^2)}$ at $t = 0.05$, $0.15$, $0.35$, $0.5$. Two dimensional projections onto $z=0$ with black lines indicating the $\pi/3$ total opening angle of the radiation cone.}
\label{Gauss-radiation}
\end{figure}

We start here with initial data of the form
\begin{equation}\label{ID:Gauss}
u(x,y,z,0)= A \, e^{-(x^2+y^2+z^2)}, ~~ A \gg 1,
\end{equation}
with sufficiently large $A$ (so that the solution would not all disperse into radiation and could form  solitons). 
In Fig. \ref{Gauss_final} - \ref{Gauss-radiation} we discuss the solution with initial amplitude of the Gaussian $A=10$. 

As the solution evolves with time, a soliton appears to emerge. 
A snapshop of the solution at $t = 12$ is shown on the top left of 
Fig.~\ref{Gauss_final}, on the top right we track the $L^\infty$ norm, which stabilizes around time $t = 3$, approximately the time when the soliton $Q_c$ forms. To check the shape of the solution against the rescaled soliton $Q_c$, we show 
the difference between the solution and $Q_c$
on the bottom left of Fig. \ref{Gauss_final}, noting that the resulting profile differs from a soliton by less than 2\%. This is roughly the size of the radiation, which cannot 
escape to infinity (since we work on a periodic domain). Here, we use a co-moving frame \eqref{ZKcom} with the speed $v_{x}=2$, which is slightly slower 
than the velocity of the resulting soliton. The Fourier coefficients 
of the solution at the terminal computational time $t=12$ are plotted 
on the bottom right of Fig. \ref{Gauss_final} and indicate the good numerical resolution of the solution.

To study further the soliton resolution in this example, we consider 
the initial stages of the time evolution for these initial data more closely and see the formation of outgoing radiation.  

The initial stages of the ZK time evolution are shown in Fig. \ref{Gaussslices}. On the right there are the isocurves of the solution in three coordinate planes, and as time increases one can clearly see that the radiation is propagating in the negative direction of the $x$-axis, widening out in the $y$ and $z$ directions into a cone-type region. 

Since the solution is symmetric in the $y$ and $z$ coordinates, we suppress for a moment the $z$-coordinate, and plot how the radiation develops in Fig. \ref{Gauss-radiation} for different times ($t=0.05$, $0.15$, $0.35$ and $0.5$). The angle of the radiation is $30^0$ from the negative $x$-axis (or the total opening angle is $60^0$) as shown by the black lines, which corresponds to the cone $C$ in \eqref{E:cone}-\eqref{E:angle}, i.e., $60^0$ to the $y$ or $z$-axis.  

\smallskip

\underline{\bf Case (b): flattened Gaussian initial data.}
\smallskip

We next examine a flattened (in the $y$ and $z$ directions) Gaussian of the form 
\begin{equation}
u(x,y,z,0)= A \, e^ {-(x^2+0.05 \rho^2)}, ~~\rho^2 = y^2 +z^2, ~~ A >1.
\end{equation}

\begin{figure}[!htb]
\includegraphics[width=0.49\hsize]{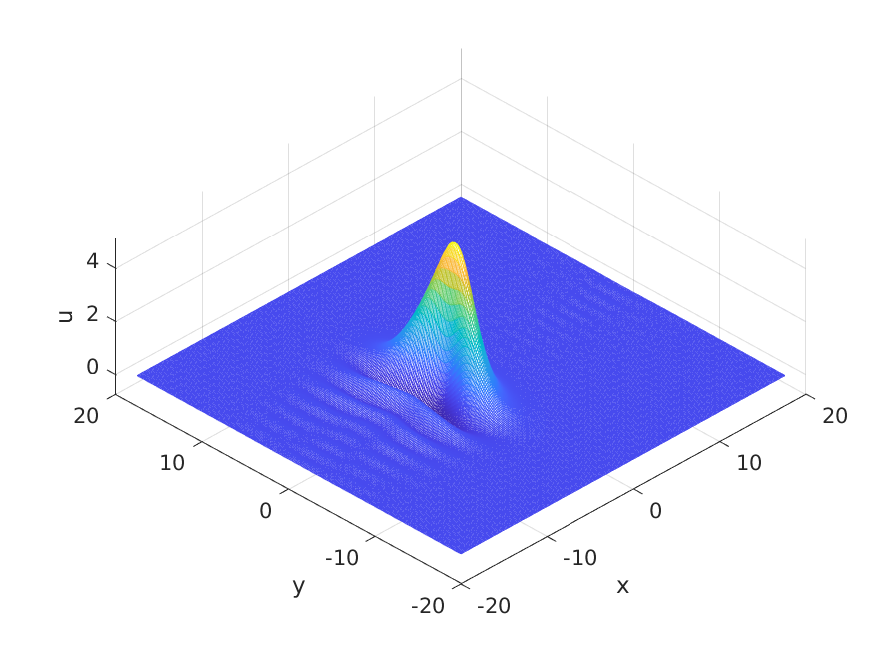}
\includegraphics[width=0.49\hsize]{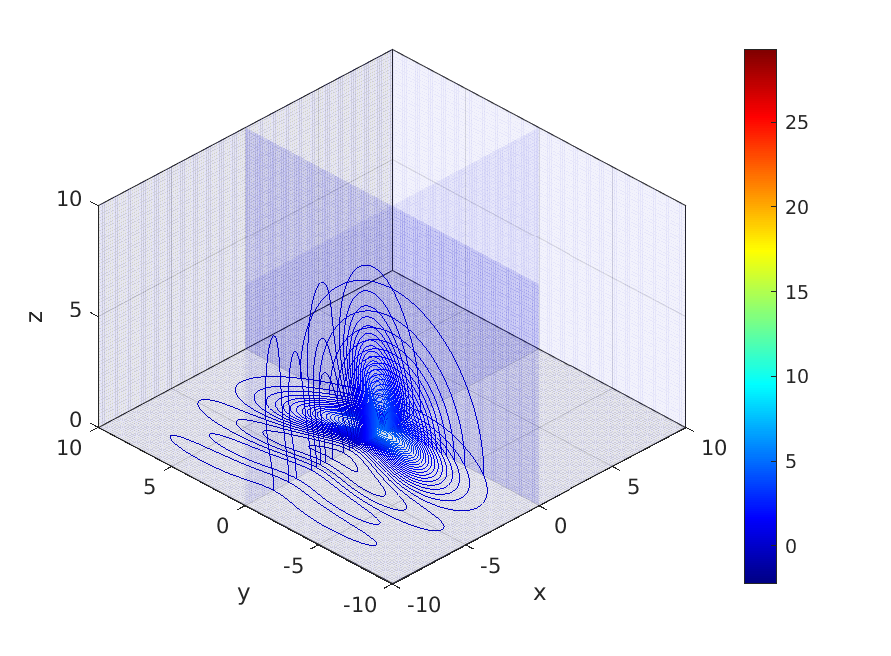}\\      
\includegraphics[width=0.49\hsize]{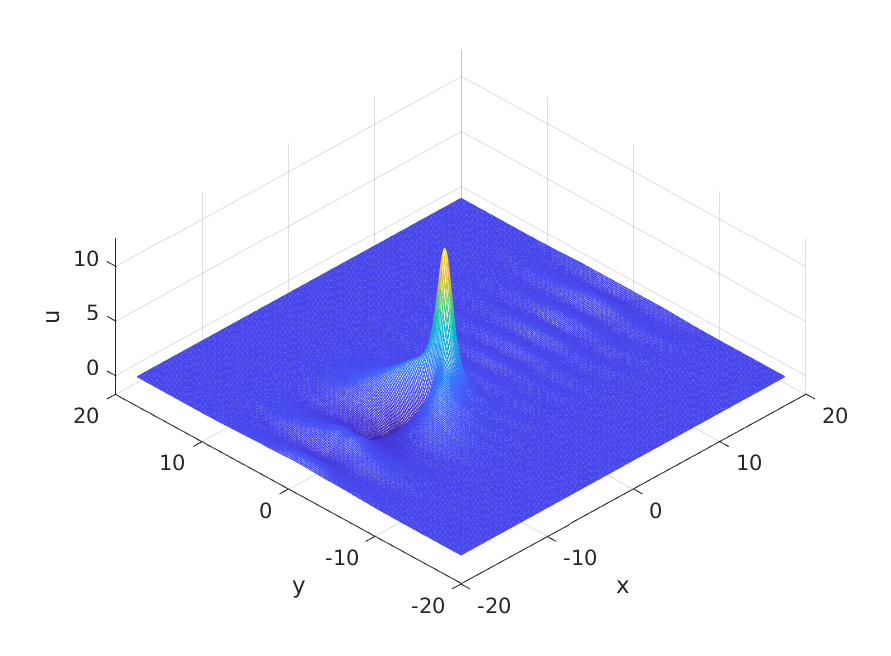}
\includegraphics[width=0.49\hsize]{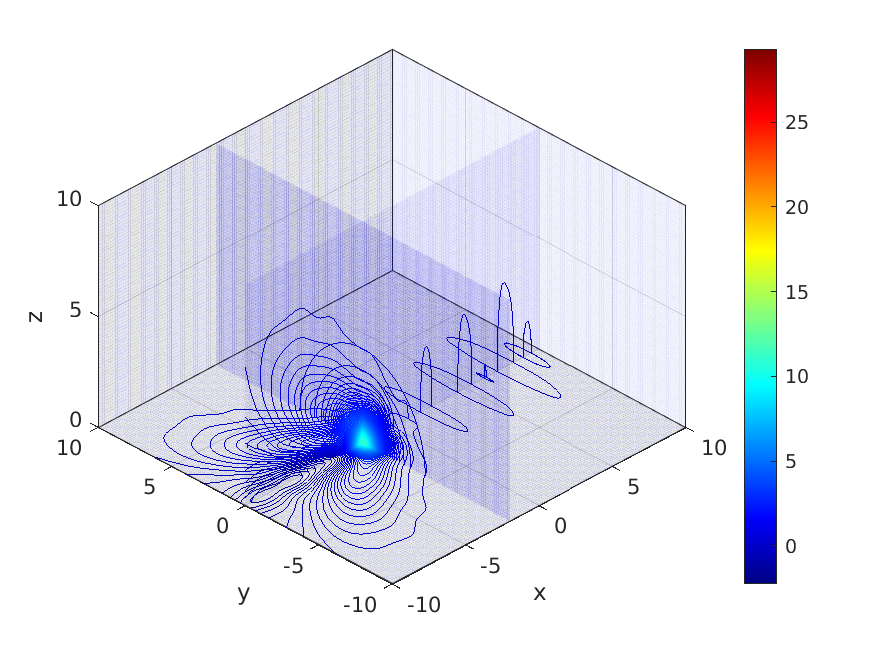}\\
\includegraphics[width=0.49\hsize]{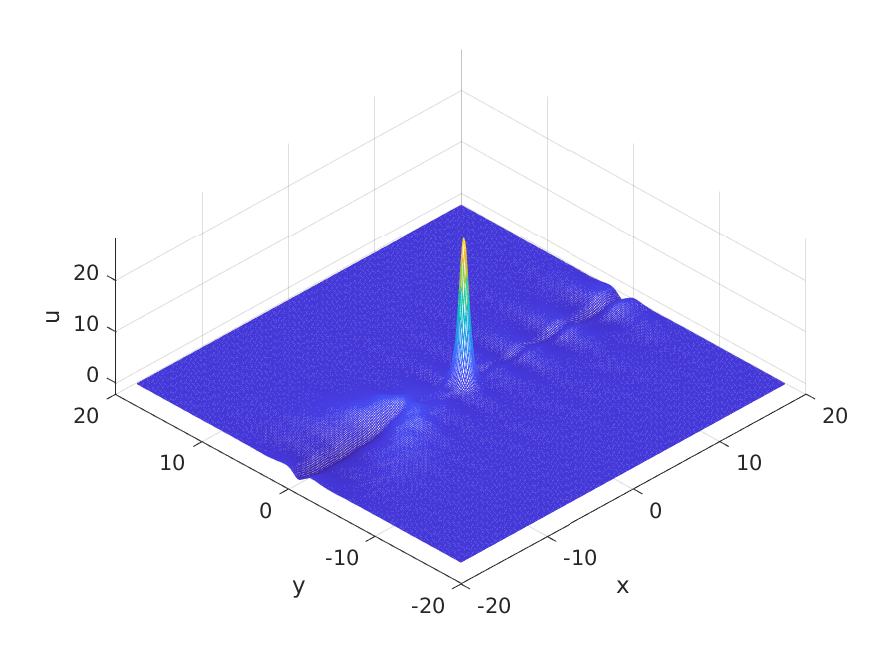}
\includegraphics[width=0.49\hsize]{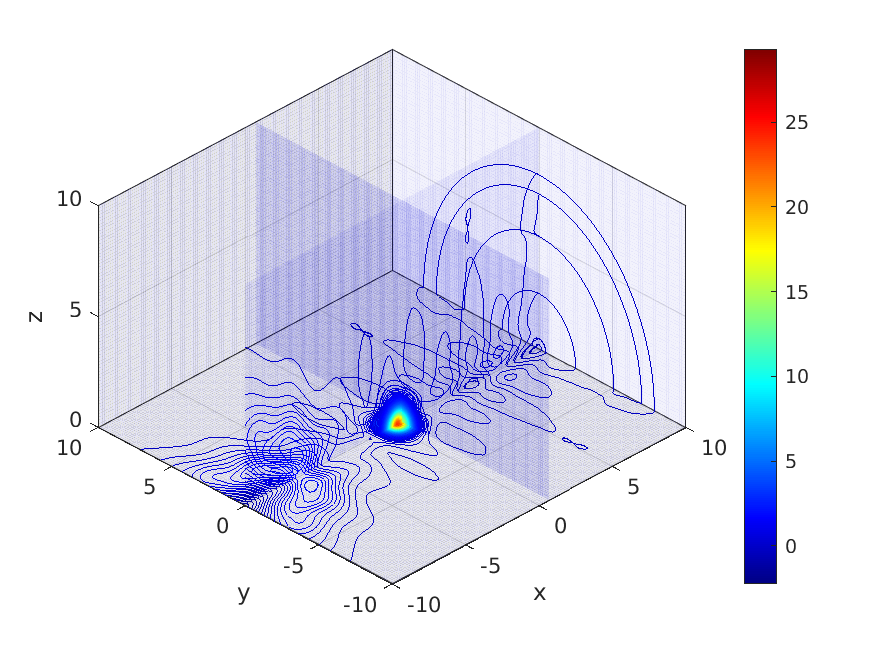}    
\caption{Snapshots of the ZK solution with flattened Gaussian initial data $u_0 = 5 \, e^{-(x^2+0.05 \rho^2)}$ at  $t = 0.4$, $1.6$, $2.8$. Projections onto the plane $z=0$ (left), the 3D isocurves on the slices of the coordinate planes (right).}
\label{ZK_FlatGauss_snaps}
\end{figure}

We take $A=5$ and plot the early stages (at times $t=0.4$, $1.6$, $2.8$) of the corresponding ZK evolution in Fig.~\ref{ZK_FlatGauss_snaps}. Similar to the previous case, one can observe that the radiation is outgoing in the negative $x$-direction (and due to the periodic domain, reappearing on the far right of the $x$-axis, which does not interfere with the soliton formation until much later times). The isocurves on the right of Fig. \ref{ZK_FlatGauss_snaps} show that the radiation does go out at some angle around the negative $x$-axis, though it is more challenging to give a precise description of the radiation region.  

\begin{figure}[!htb]
\includegraphics[width=0.51\hsize]{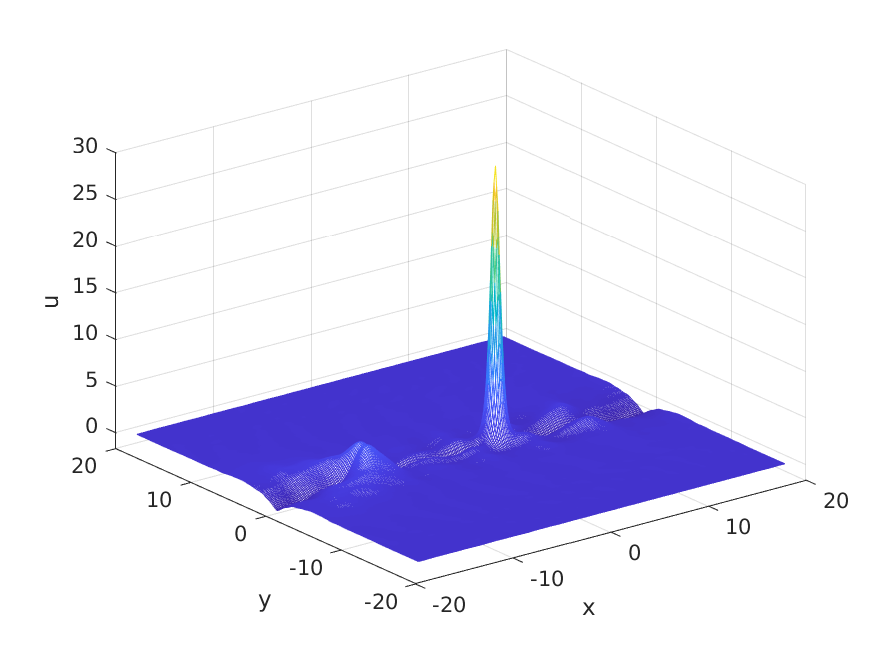}
\includegraphics[width=0.46\hsize]{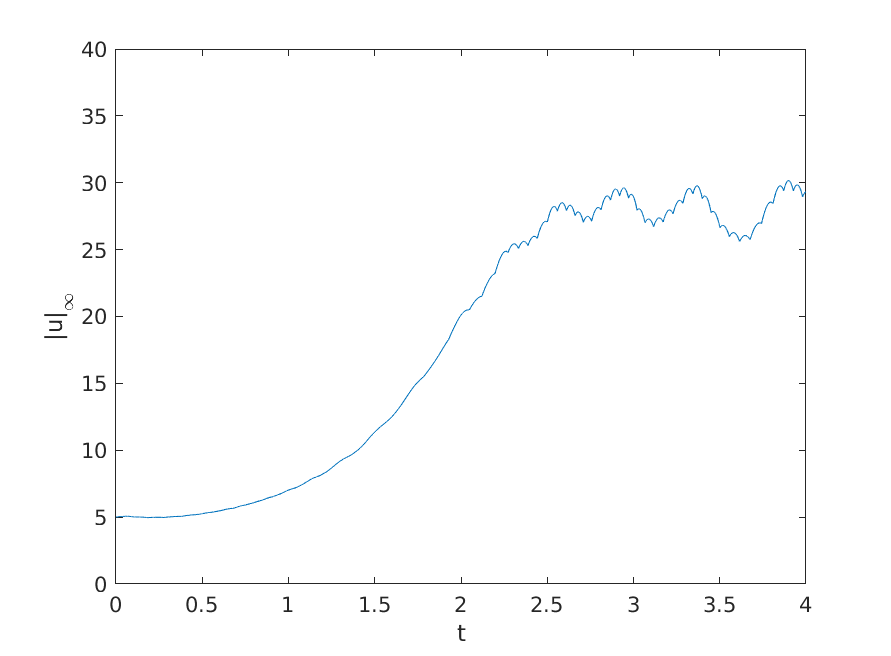}
\includegraphics[width=0.49\hsize]{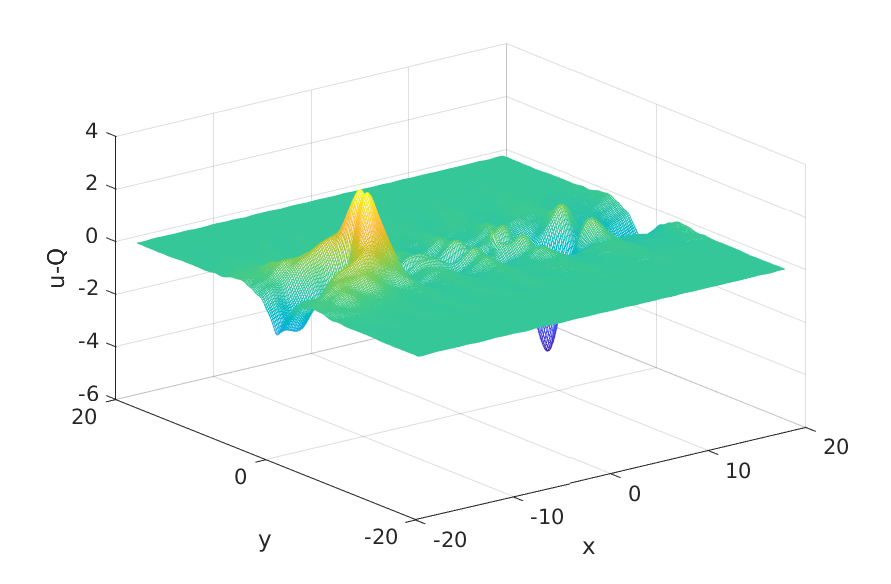}
\includegraphics[width=0.49\hsize]{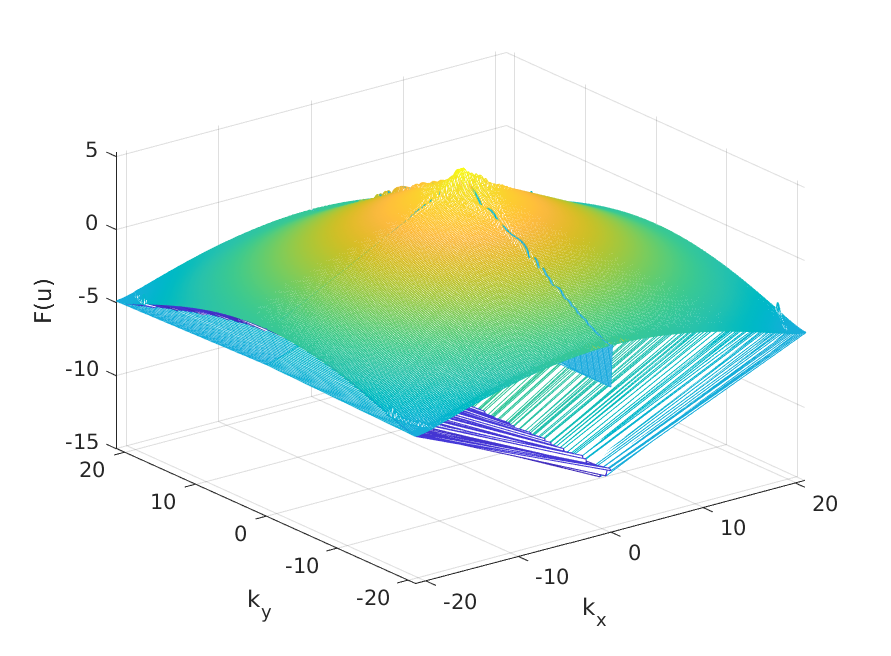}
\caption{ZK solution with flattened Gaussian initial data  
$u_0=5\, e^ {-(x^2+0.05 \rho^2)}$ at $t=4$ (top left), 
the $L^{\infty}$ norm (top right), 
the difference of the solution (at $t=4$) and a rescaled soliton $Q_c$ (bottom left), 
the Fourier coefficients of the solution at $t=4$ (bottom right). In all plots the presented solution is projected onto $z=0$.}
\label{Flat_final}
\end{figure}

The ZK solution of this flattened Gaussian data (with $A=5$) at the time $t=4$ is plotted on the top left of Fig.~\ref{Flat_final}. Tracking the $L^\infty$ norm of this solution on the top right of Fig. \ref{Flat_final} shows that the solution grows for some time (till about $t \approx 2.5$) until it reaches the height of the appropriate rescaled soliton while shedding some radiation.  
The oscillations in the $L^\infty$ height are due to the periodic domain setting as well as the reappearence of the radiation on the right. The difference between this solution and a rescaled  soliton $Q_c$ on the bottom left of the same figure once more indicates that the solution asymptotically approaches the soliton as its final state. 
The Fourier coefficients of the solution at $t=4$
are given on the bottom right of the same figure, showing that the solution is numerically well 
resolved. 
 
\newpage

\underline{\bf Case (c): Wall-type initial data.}
\bigskip

Instead of single maximum initial data, we consider a wall-like setup 
with the maximum value spread out continuously over an interval (for 
instance $-a\leq y+z \leq a$), while still having a fast decay and vanishing at infinity. For example, we take 

\begin{equation}\label{wallcond}
u(x,y,z,0) = 
\begin{cases}
 A \, e^{-x^{2}} & |y+z|\leq a \\
 A \, e^{-(x^{2}+(y+z-a)^{8})} & y+z > a \\
 A \, e^{-(x^{2}+(y+z+a)^{8})}& y+z<-a.
\end{cases}
\end{equation}
\medskip

\begin{figure}[!htb]
\includegraphics[width=0.49\hsize]{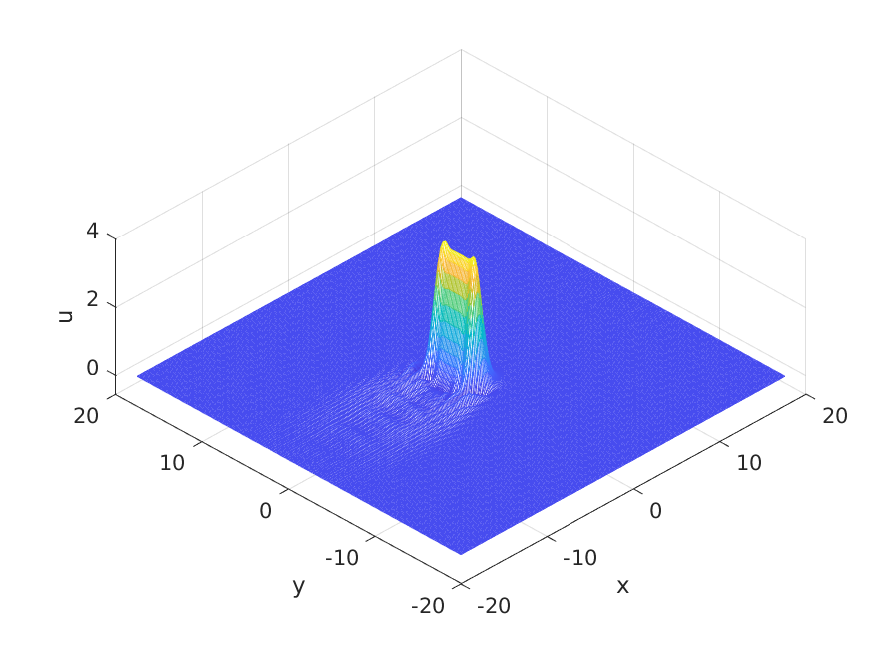}
  \includegraphics[width=0.49\hsize]{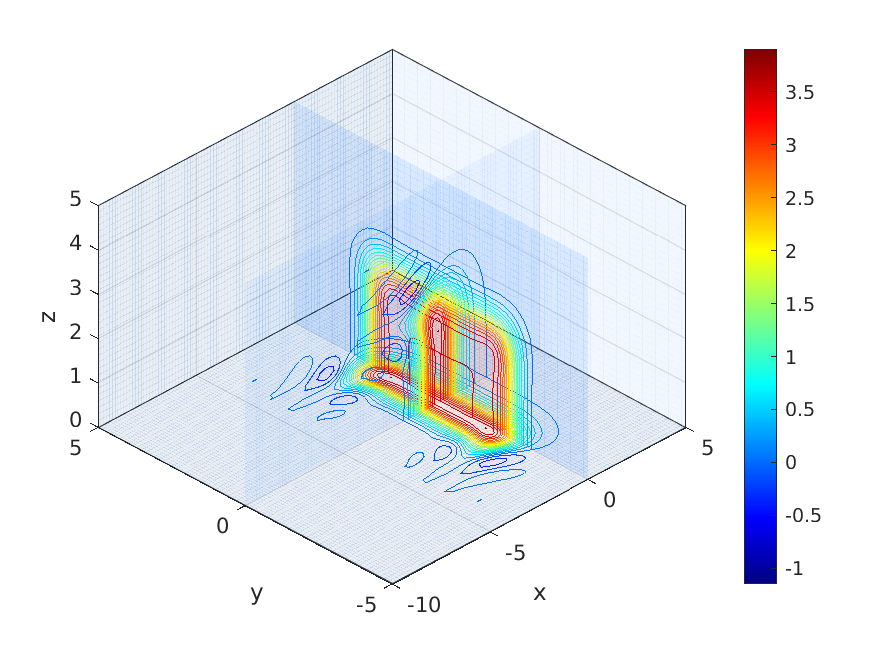}\\      
\includegraphics[width=0.49\hsize]{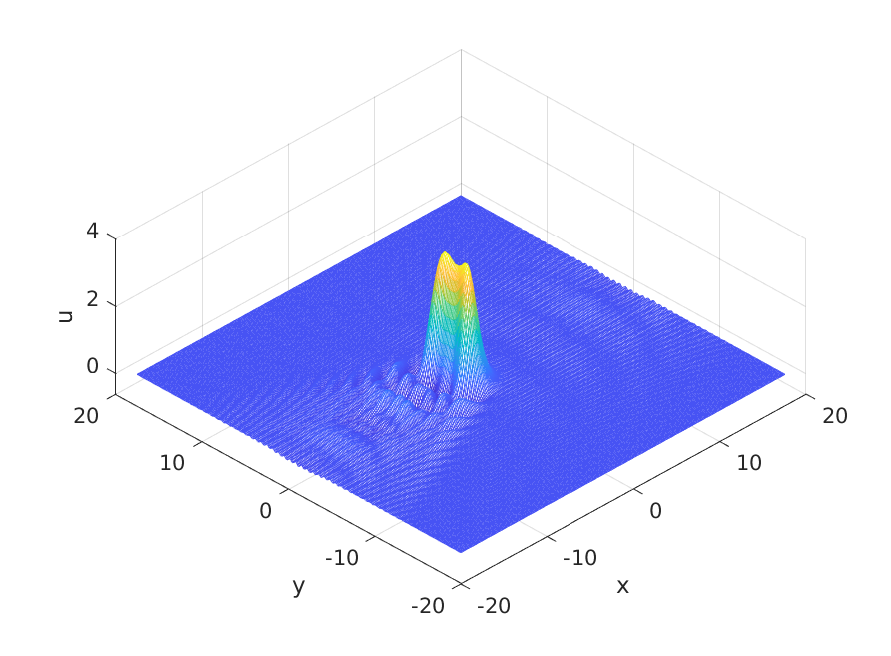}
  \includegraphics[width=0.49\hsize]{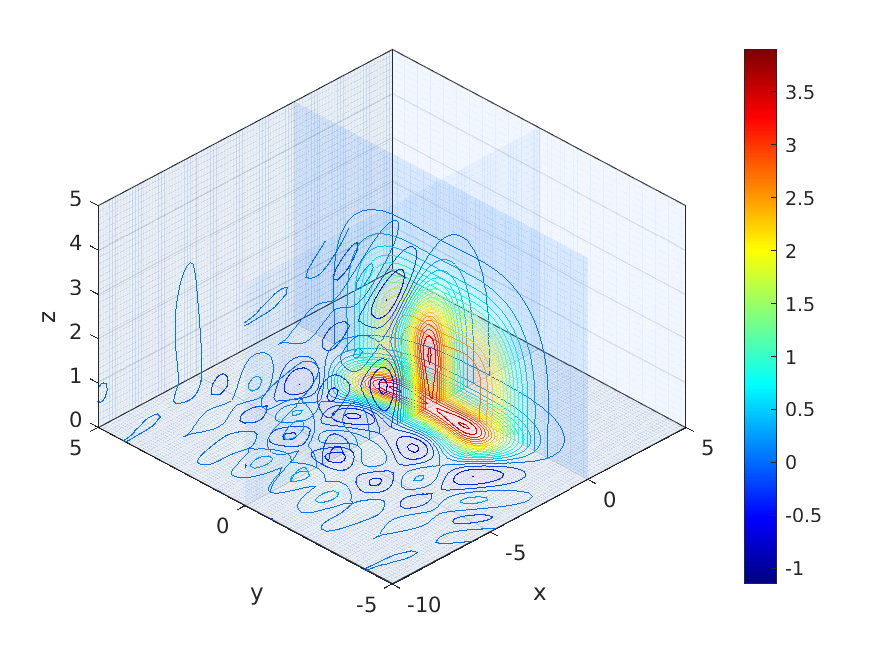}\\
  \includegraphics[width=0.49\hsize]{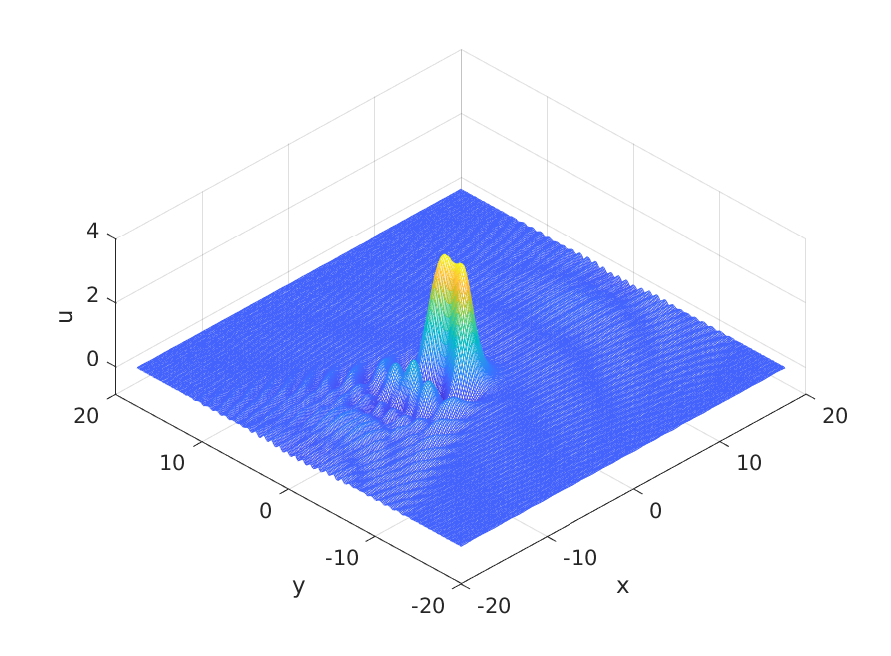}
  \includegraphics[width=0.49\hsize]{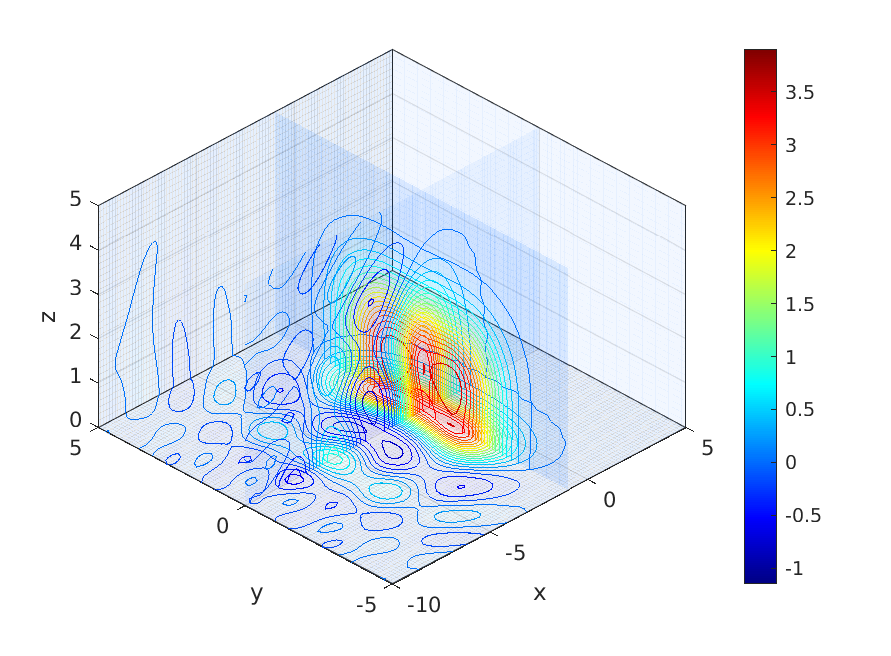}\\ %
\caption{Snapshots of the ZK solution with the wall-type initial data  
\eqref{wallcond} at $t=0.05$, $0.2$, $0.35$. Projections onto the plane $z=0$ (left), 
the 3D contour plots on the slices of the coordinate planes (right).}
\label{ZK_wall_rad}
\end{figure}

Our goal here is to study the time evolution of initial data that are not 
single peaked, but have the same maximum along a certain set, for example, an interval or a curve. In the wall-type data above \eqref{wallcond} we have the maximum elongated in the $y$ and $z$ directions, and in 
the $x$-direction the maximum is still localized in a single point (at $x=0$). The ZK evolution for such data with $A=3.6$, $a=1.5$ at $t=0.05$, $0.2$, $0.35$ is shown in Fig.~\ref{ZK_wall_rad}.   

\newpage

\begin{figure}[!htb]
\includegraphics[width=0.51\hsize]{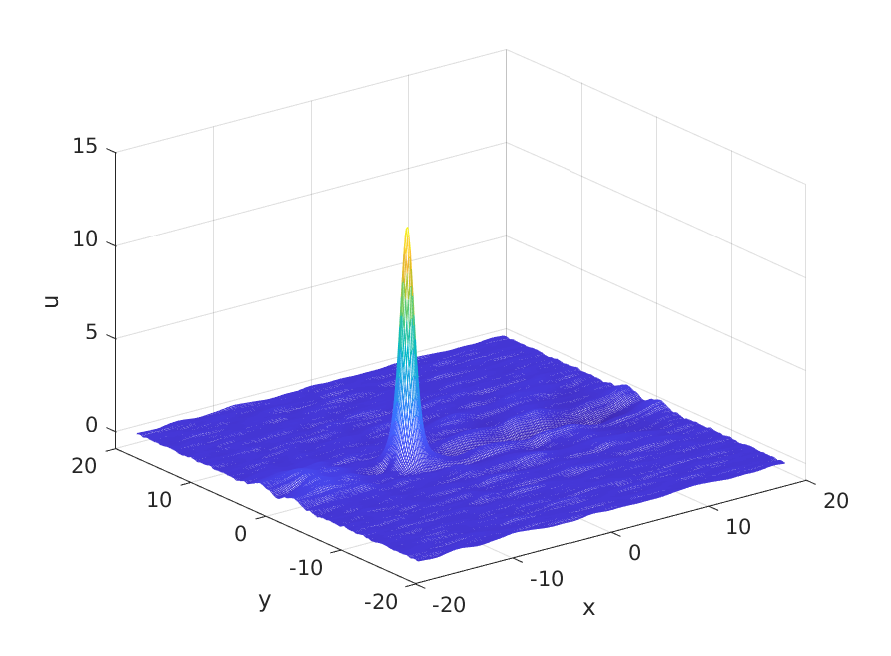}
\includegraphics[width=0.46\hsize]{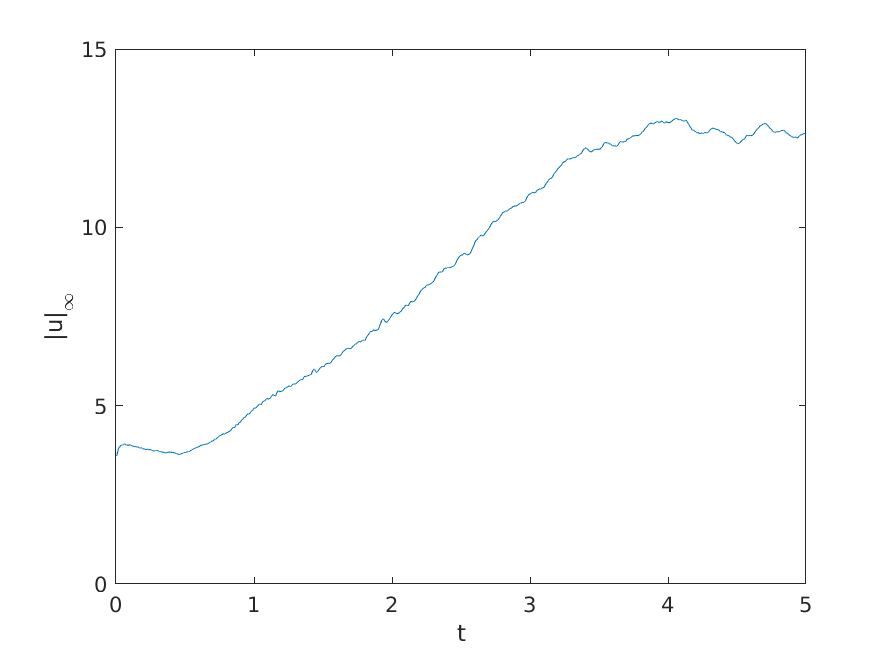}
\includegraphics[width=0.49\hsize]{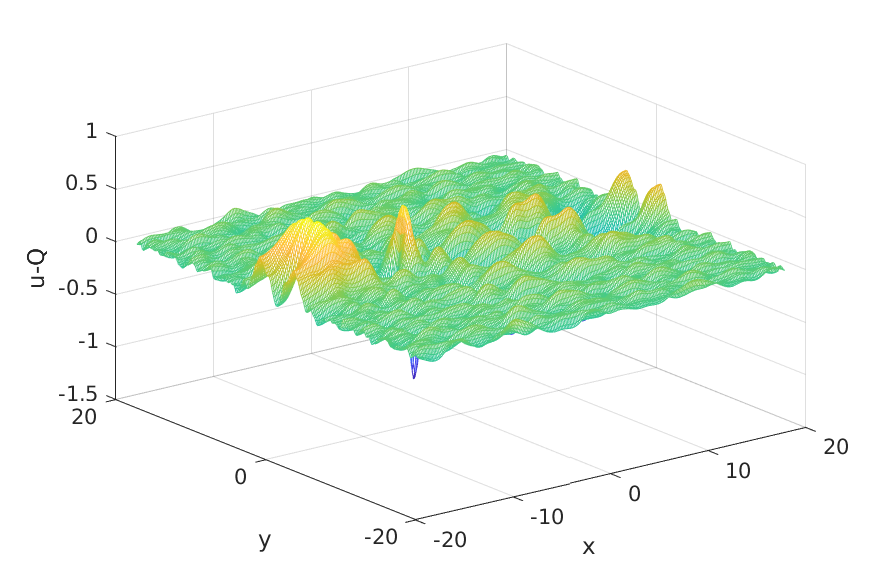}
\includegraphics[width=0.49\hsize]{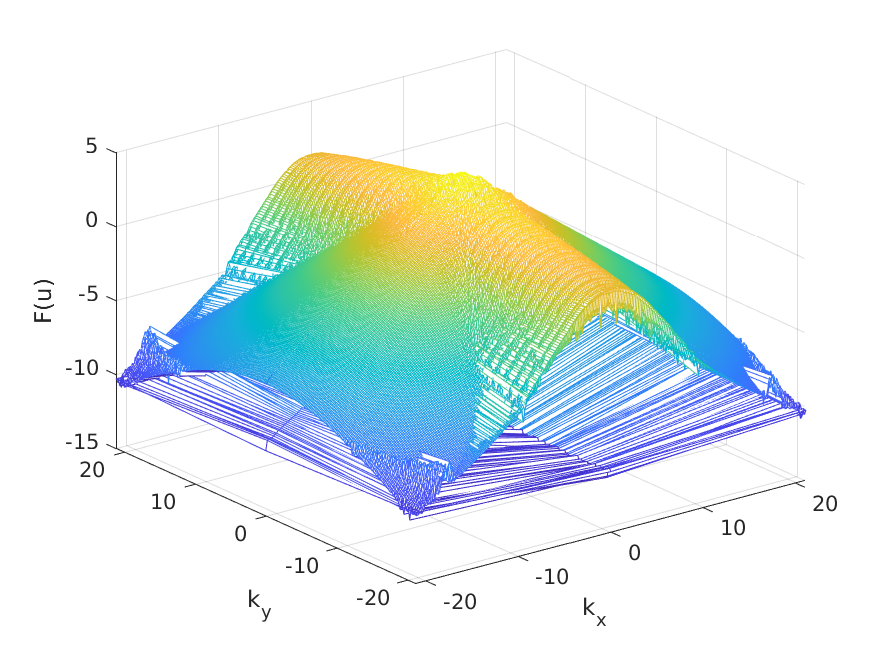}
\caption{Soliton resolution and control parameters for the wall-type initial data \eqref{wallcond}. The profile of the solution at $t=5$ (top left), time dependence of the $L^\infty$ norm (top right), the difference of the solution at $t=5$ and a rescaled soliton (bottom left), the Fourier coefficients at $t=5$ (bottom right).}
\label{ZK_wall}
\end{figure}

Continuing tracking the time evolution of this solution, we observe that eventually it forms one single peak, which then moves along the positive $x$-axis. The corresponding solution at time $t=5$ is plotted on the top left of Fig.~\ref{ZK_wall}. 

The $L^{\infty}$ norm of the solution is given on the top right of Fig.~\ref{ZK_wall}, it shows that it is growing and then stabilizes around time $t=4$. 
The difference of the solution at time $t=5$ and a rescaled soliton $Q_c$ is shown on the bottom left of the same figure: it does indicate that the single peak bump, into which the solution evolved thus far, is very close to a rescaled soliton, thus, showing that the solution is asymptotically approaching as $t\to\infty$ a rescaled (and shifted) soliton. This, furthermore, confirms the soliton resolution conjecture for initial data with non-single peaked maximum.   
We note that since the initial condition is not very smooth (at the edges), the spatial resolution suffers, as indicated by the Fourier coefficients of the 
solution at the final time on the bottom right of the same figure. Nevertheless, this gives a positive confirmation to the soliton resolution conjecture. 

We also comment on the radiation cone in this example: 
we observe that initially the angle $\theta$ of the radiation front with respect to the $x$-axis is such that $\tan \theta = 4$, and as the wall narrows into a soliton (for example, compare the top left and bottom left plots in Fig.~\ref{ZK_wall_rad}), the angle gets closer to $\pi/6$, which is easier to see here than in the case of the flattened Gaussian example.    
\smallskip

\clearpage

\underline{\bf Case (d): Algebraic decay initial data.}
\smallskip

So far we have been studying data with exponential decay, either the 
same rate as the decay of a soliton, or even faster such as a Gaussian. We now consider polynomial decay data. 
Currently, we are able to consider sufficiently fast polynomial decay (such as $|x|^{-20}$), that is, an algebraic decay 

\begin{figure}[!htb]
\includegraphics[width=0.49\hsize]{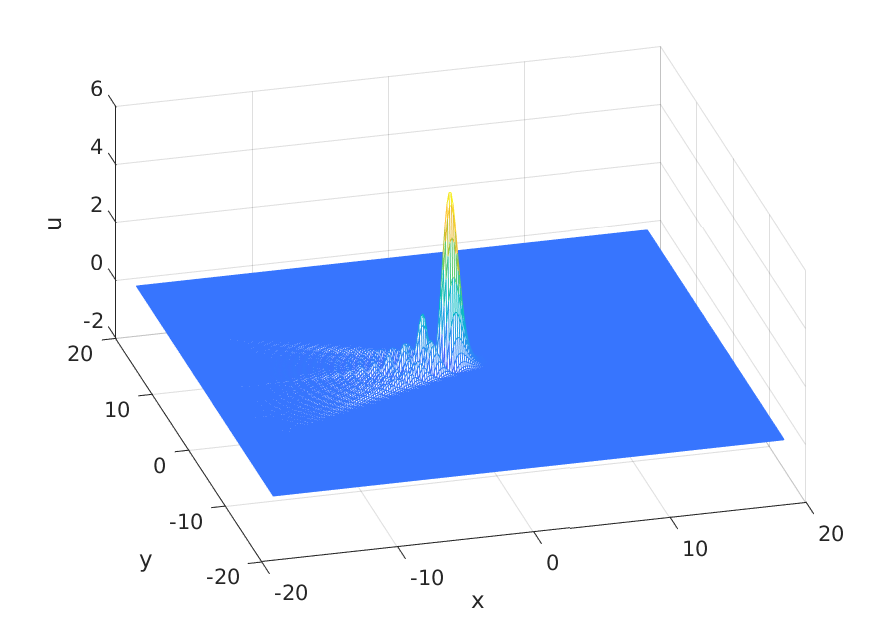}
  \includegraphics[width=0.49\hsize]{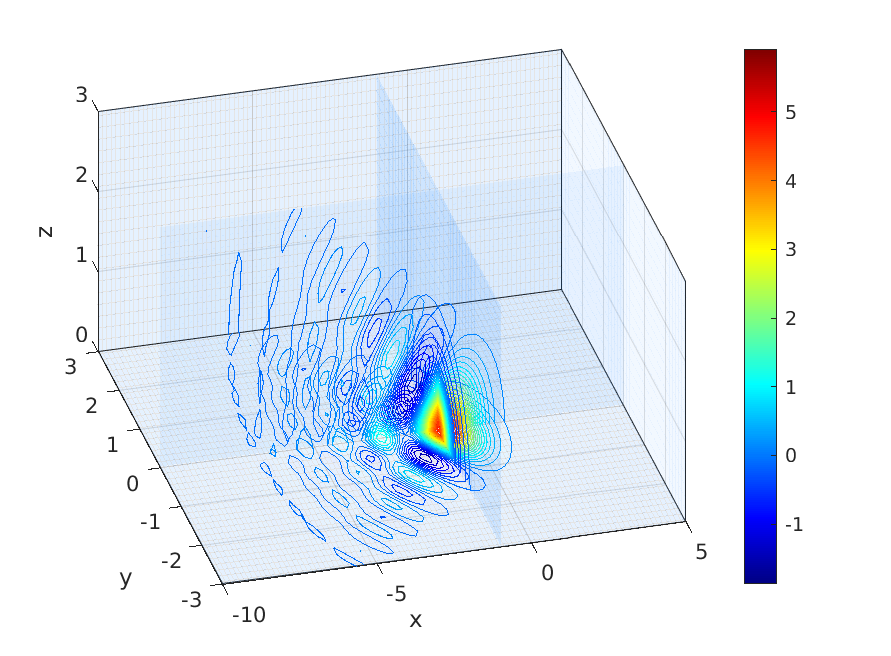}\\      
\includegraphics[width=0.49\hsize]{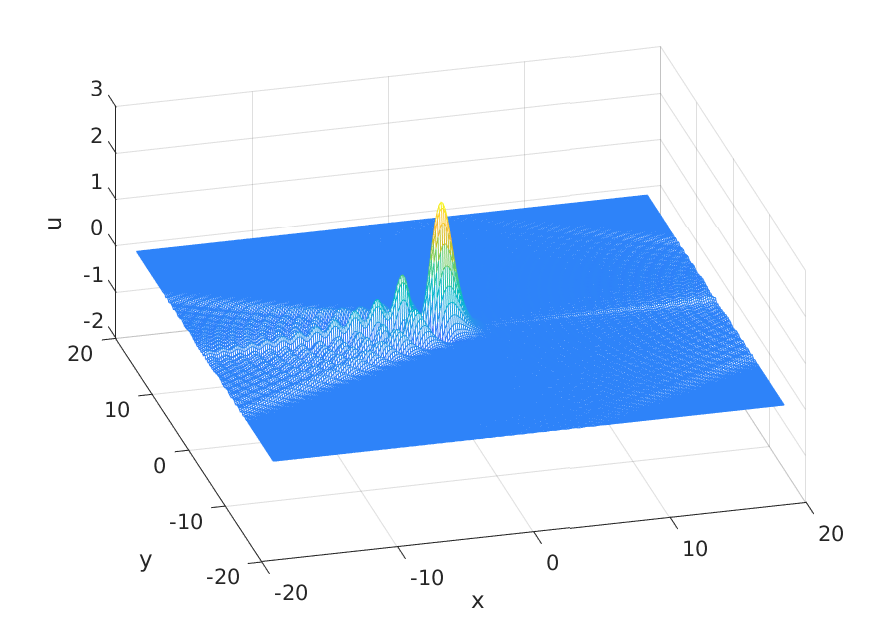}
  \includegraphics[width=0.49\hsize]{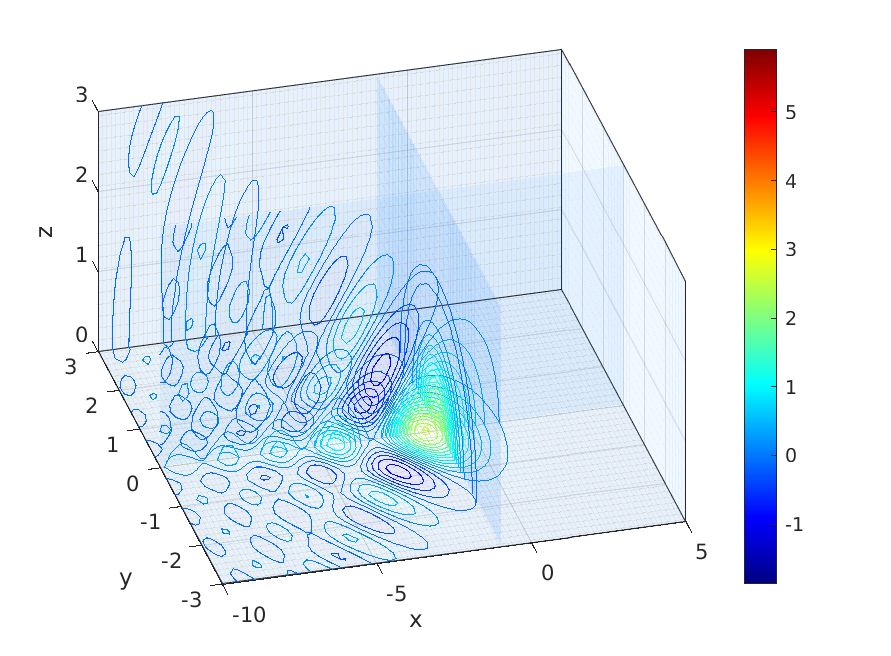}\\ 
\includegraphics[width=0.49\hsize]{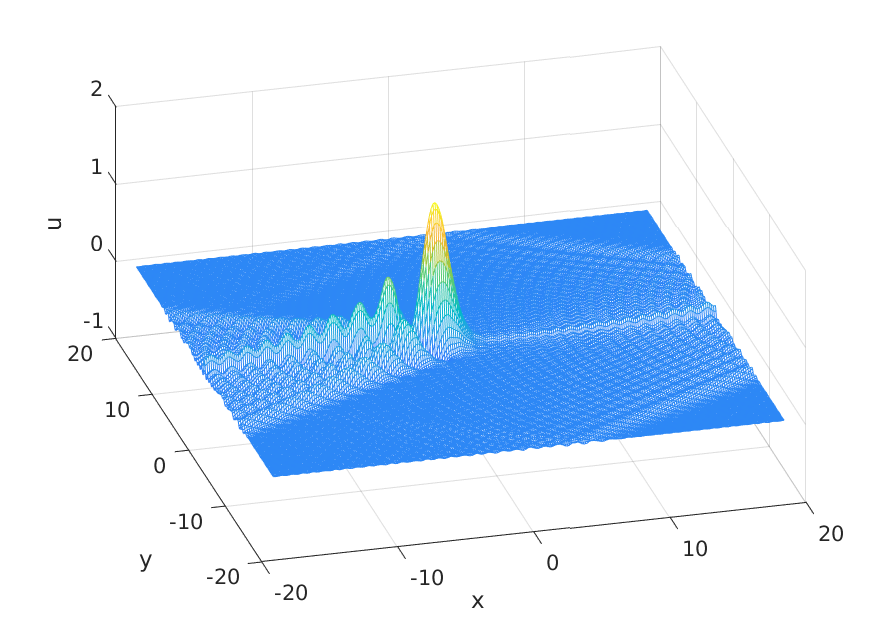}
  \includegraphics[width=0.49\hsize]{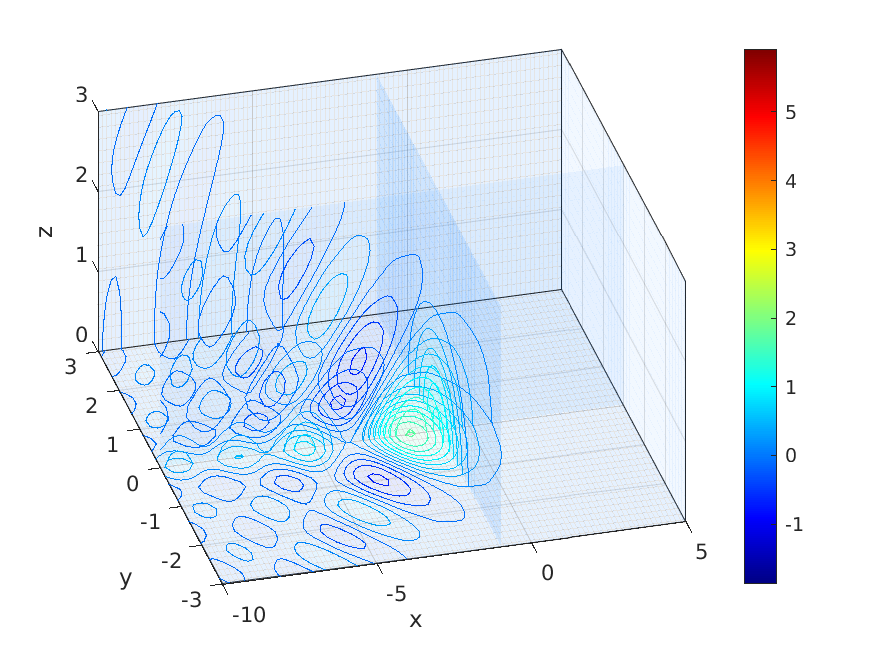}\\ 
\caption{Snapshots of the ZK solution with $u_0 =  {20}/{(1+ x^2 + y^2 +z^2)^{10}}$ at $t = 0.05$, $0.15$, $0.25$. 
Left: 2D projections onto $z=0$ (note the changing scale for $u$).
Right: 3D isocurves on the slices of the coordinate planes.}
\label{ZK_SuperLorentz_snaps}
\end{figure}

\noindent toward infinity is rapid enough so that Fourier methods could be used efficiently.
We first take initial data with the algebraic decay of the form
\begin{equation}\label{E:polynomial}
u(x,y,z,0) =  \frac{A}{(1+ x^2 + y^2 +z^2)^{10}}, ~~A \gg1.
\end{equation}

The corresponding solution, despite the significant mass, if we take $A$ large (e.g. $A=10,20$), appears to disperse without forming a soliton. The snapshots of the solution with $A=20$ at different times  are shown in Fig. \ref{ZK_SuperLorentz_snaps}.

\begin{figure}[!htb]
\centering
\includegraphics[width=0.6\hsize]{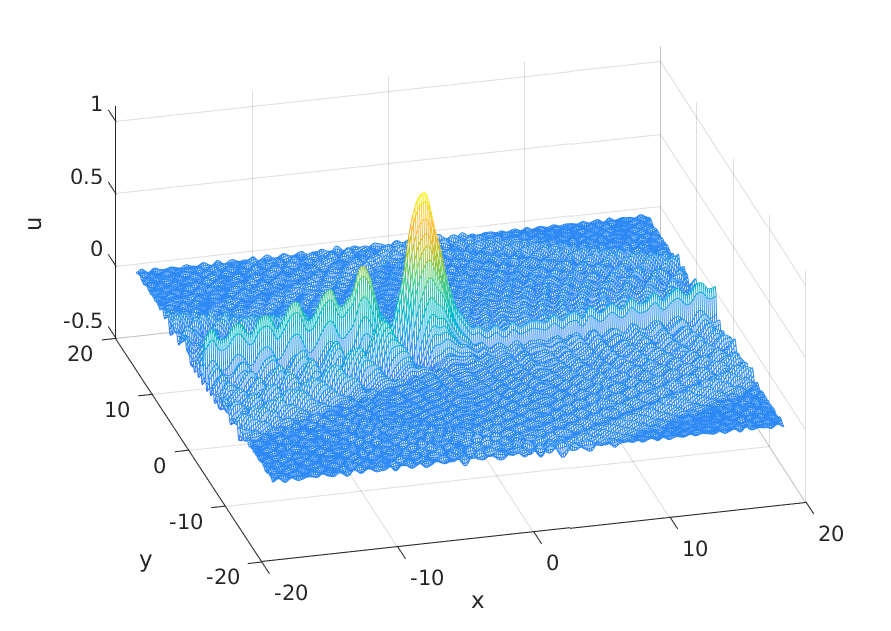}\\
 \includegraphics[width=0.49\hsize]{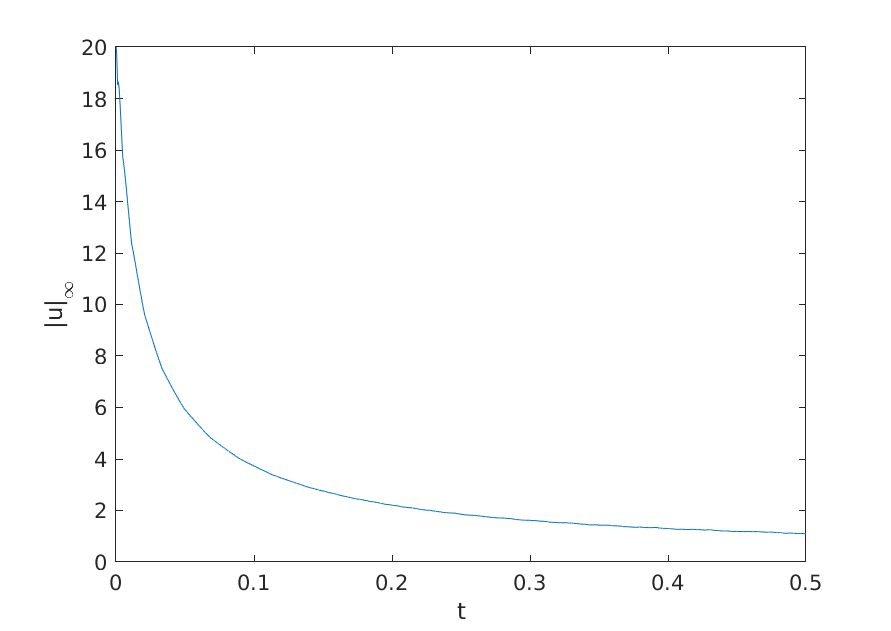}
 \includegraphics[width=0.49\hsize]{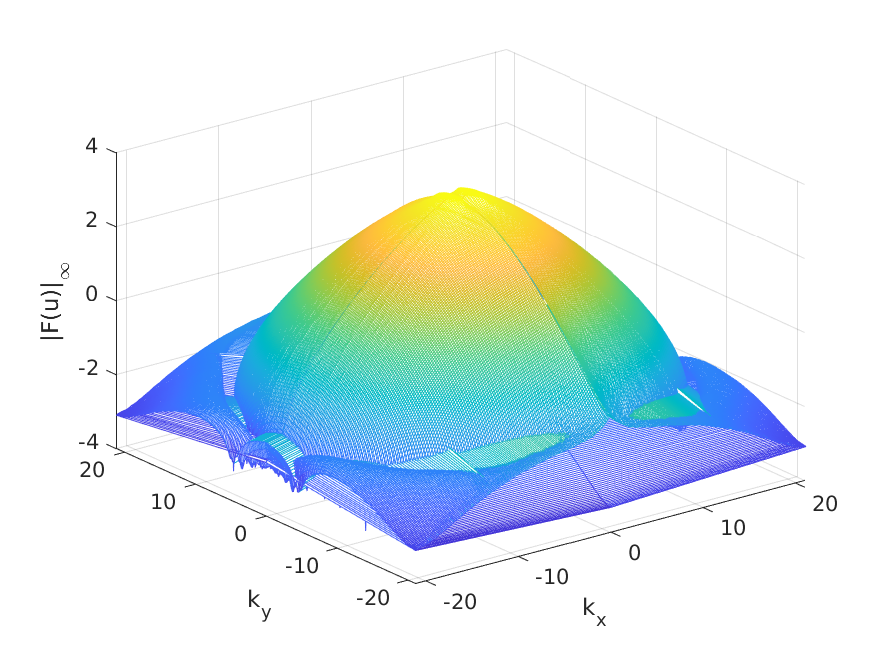}
\caption{ZK solution with super-Lorentzian $u_0 = {20}/{(1+ x^2 + y^2 
+z^2)^{10}}$: the projection onto the plane $z=0$ at $t=0.5$ (top); the $L^{\infty}$ norm 
depending on time (bottom left), the Fourier coefficients at $t=0.5$ 
depending on $k_{x}$ and $k_{y}$ on the right (the dependence on $k_{z}$ is suppressed due to the symmetry in $y$ and $z$). }
\label{ZKSuperLorentz}
\end{figure}

We continue to track the solution even though the radiation keeps 
reappearing on the right as long as it does not influence the 
solution too much, till $t=0.5$.  The snapshot of that time is given 
on the top of Fig.~\ref{ZKSuperLorentz}.

On the bottom left plot we track the $L^{\infty}$ norm in time, which appears to be monotonously decreasing, thus, it is plausible to think that this solution 
appears to be purely radiative and no solitons would be formed. The Fourier coefficients of the solution at $t=0.5$ decay considerably less than those for the exponentially localized initial data, but still indicate a numerical resolution to the order of the plotting accuracy. 
We also check the relative energy conservation, which is on the order of $10^{-7}$, indicating that the resolution in time does not cause any concerns.

Qualitatively an identical result is obtained for the initial condition with faster algebraic decay:
\begin{equation}\label{E:SL-2}
u(x,y,z,0) =  \frac{20}{(1+ x^2 + y^2 +z^2)^{20}},
\end{equation}
with all of the solution dispersing into the radiation. 
We note that since the above initial data do not have exponential decay, 
the spatial resolution with Fourier methods is not on the same order as in the exponentially decaying examples. Nevertheless, it still provides six orders of magnitude precision in space. Temporal resolution is of the same order. 
We also note that as we compute on a finite grid, the exponential 
decay should become indistinguishable from the algebraic decay at the order of 
approximately $N/2$. 
The tested super-Lorentzian initial data \eqref{E:SL-2} with the decay rate $~|x|^{-40}$ still exhibits ``not sufficiently compactified" behavior, and therefore, we conclude that the initial data of the forms \eqref{E:polynomial} and \eqref{E:SL-2} of algebraic decay data do not form a soliton.

\begin{figure}[!htb]
\includegraphics[width=0.51\hsize]{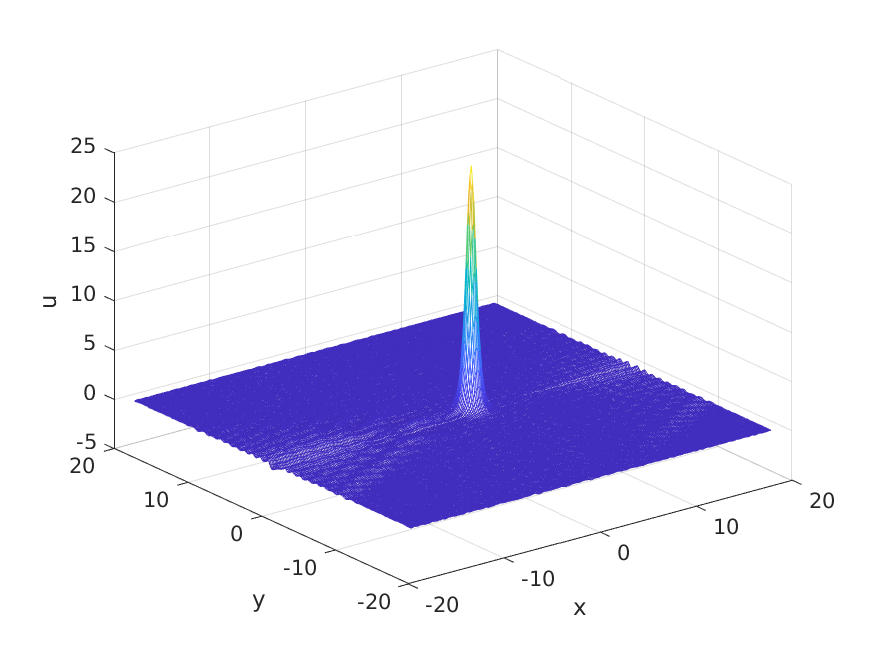}
\includegraphics[width=0.46\hsize]{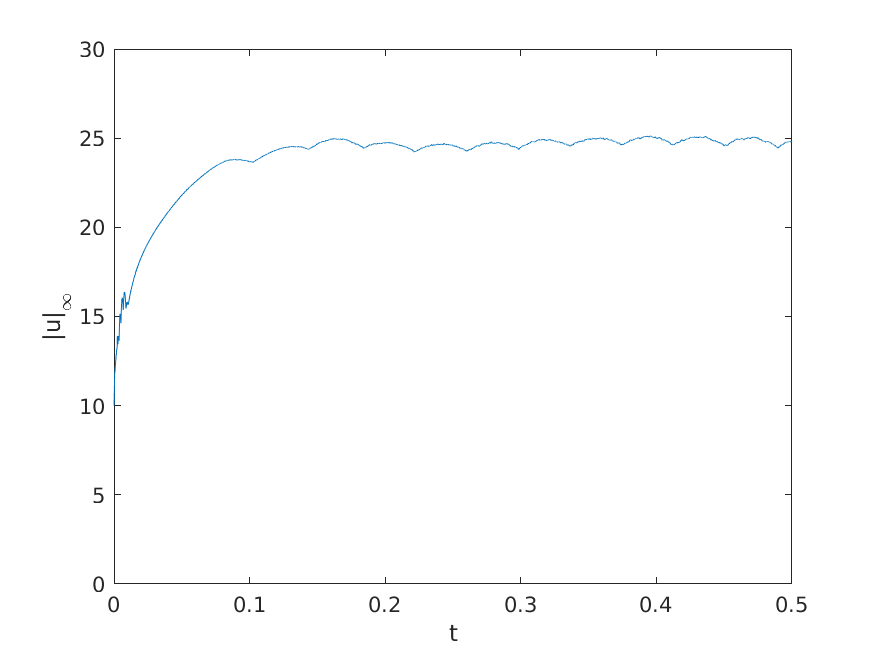}
\includegraphics[width=0.49\hsize]{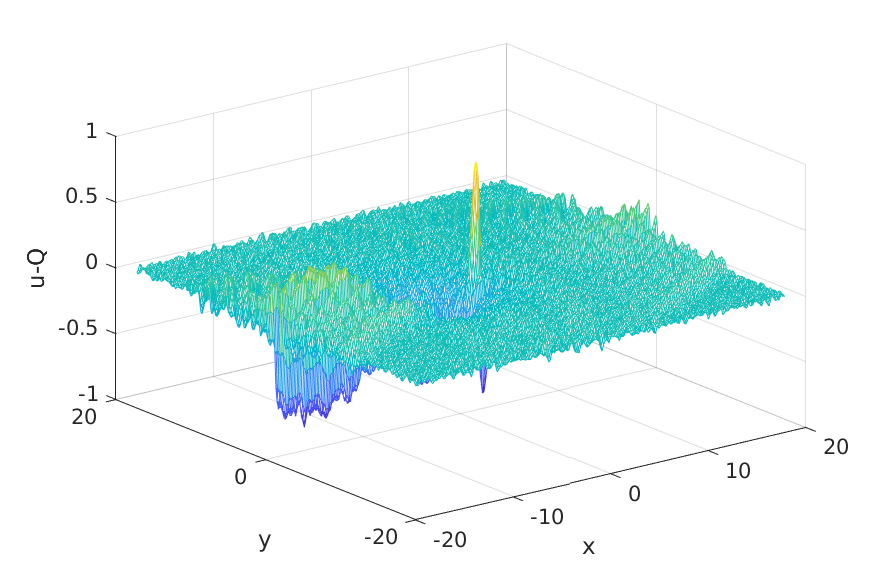}
\includegraphics[width=0.49\hsize]{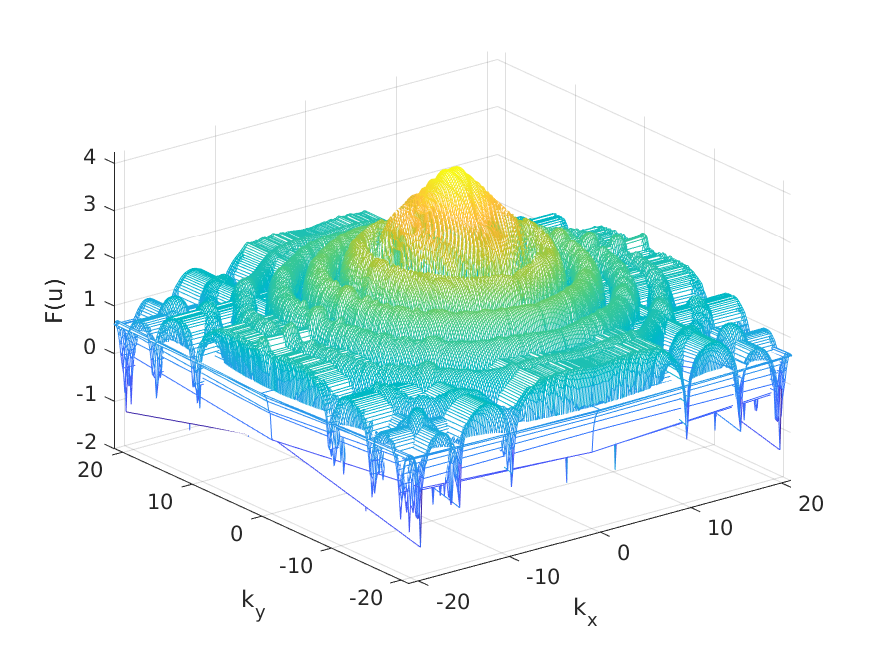}
\caption{ZK solution for the more flat polynomial initial data \eqref{E:polynomial_flat} with $A=10$: solution (projected onto $z=0$) at $t=0.5$ (top left), time dependence of the $L^\infty$ norm (top right), the difference of the solution at $t=0.5$ and a rescaled soliton (bottom left), the Fourier coefficients at $t=0.5$ (bottom right).}
\label{ZK-SL2}
\end{figure}

We next consider a more flat at the origin initial data with 
a similar decay rate 
\begin{equation}\label{E:polynomial_flat}
u(x,y,z,0) =  \frac{A}{1+ (x^2 + y^2 +z^2)^{10}}, ~~A \gg1.
\end{equation}

We plot the ZK solution with the initial condition \eqref{E:polynomial_flat} and $A=10$ at $t=0.5$ in the top left of Fig.~\ref{ZK-SL2} and the time dependence of the $L^\infty$ norm in the top right of the same figure. One can observe that around $t=0.15$ the $L^\infty$ norms starts to stabilize. The difference between the solution and a rescaled and shifted soliton $Q_c$ is obtained in the bottom left of Fig.~\ref{ZK-SL2}. The difference does confirm that this solution forms a soliton. We conclude that it is not only the decay rate that affects the asymptotic behavior of the solution but also the shape (flatness at the origin in this case) of the initial data.


\section{Soliton interaction}\label{S:Interaction}
In this section we study the interaction of solitons. Since the ZK 
equation is not completely integrable, there are no {\it exact} 
multi-soliton solutions known, as, for instance, in the KdV or mKdV equations. 
However, since the soliton solutions are rapidly decaying as $|(x,y,z)| \to \infty$,  
one can study their interactions by considering initial data, that are composed as the 
sum of displaced solitons. Because of the finite numerical precision, 
this gives multi-soliton initial data within the available accuracy. 
Indeed, since solitons have an exponential decay, their contribution to the locations that are far away from their joint center of mass is zero within the numerical precision.
Therefore, we consider initial data by superimposing two one-soliton solutions that are sufficiently well separated. 

We study the following three scenarios:
\begin{enumerate}
\item[(a)] 
Head-on collision: a classical quasi one-dimensional interaction.
\item[(b)] 
Twin solitons: the interaction of two identical solitons that 
are separated only in $y$ or in $z$ direction.
\item[(c)] 
Off-set solitons: the interaction of two solitons that are separated in $x$ and $y$ (or $z$) directions.
\end{enumerate}

We show that in some cases the solitons can interact essentially elastically, although 
with radiation, as the equation is not integrable, and in other cases the interaction is strong: the solitons can merge into a single one (also with the outgoing radiation). We mention that the mass and energy conservation that we obtain in the examples below are on the order of $10^{-8}$ or smaller. 
\smallskip

\underline{\bf Case (a): Head-on collision.}
\medskip

We start with considering the initial data with two localized and sufficiently  
separated solitons: one of them being a large soliton $Q_c$, $c>1$, shifted away from the origin in the negative $x$-direction, and another one is the soliton $Q$ itself centered at the origin: 
\begin{equation}\label{E:Qc-Q}
u(x,y,z,0) = Q_c(x - a,y,z) + Q(x,y,z), ~~c > 1. 
\end{equation}
We show an example with $c=2$ centered at $a=-10$, that is, the 
larger soliton is behind the smaller soliton on the $x$-axis, and we 
expect the larger one to travel twice as fast (at least in the beginning) than the smaller one. 
In our simulations we actually solve 
\begin{equation}\label{E:Qc}
-cQ + Q_{xx}+Q_{yy}+Q_{zz} + Q^2=0
\end{equation}
with $c=2$ to obtain the appropriate soliton $Q_2$. As an alternative, we 
could also use the scaling property \eqref{dilation} to get the soliton 
with $c=2$.

\begin{figure}[!htb]
\includegraphics[width=0.49\hsize]{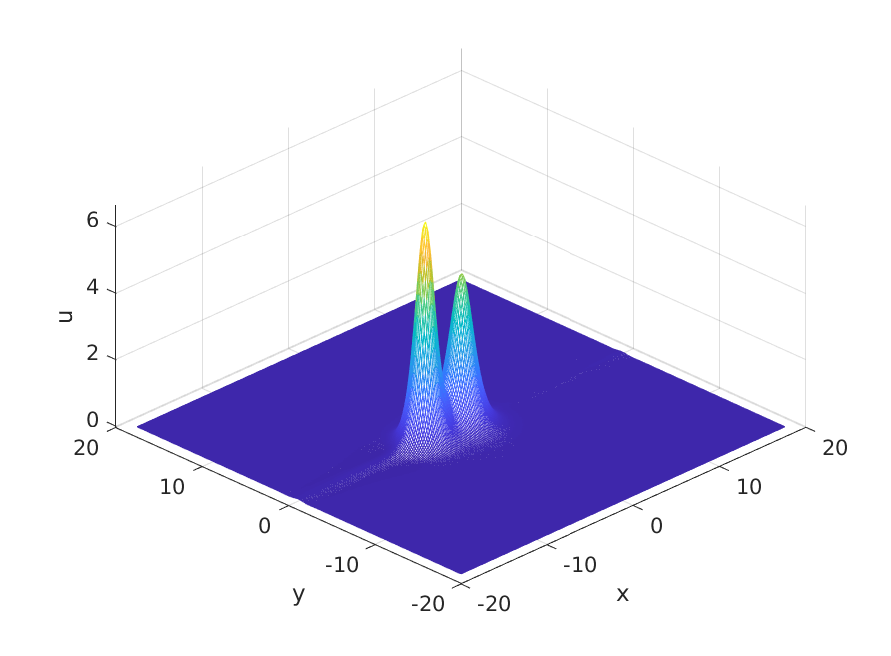}
\includegraphics[width=0.49\hsize]{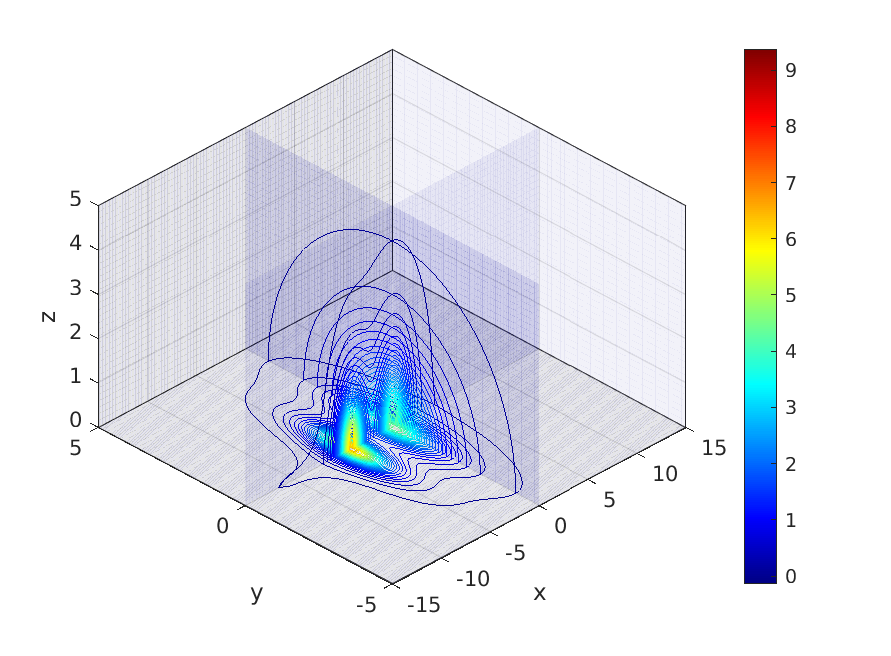}\\
\includegraphics[width=0.49\hsize]{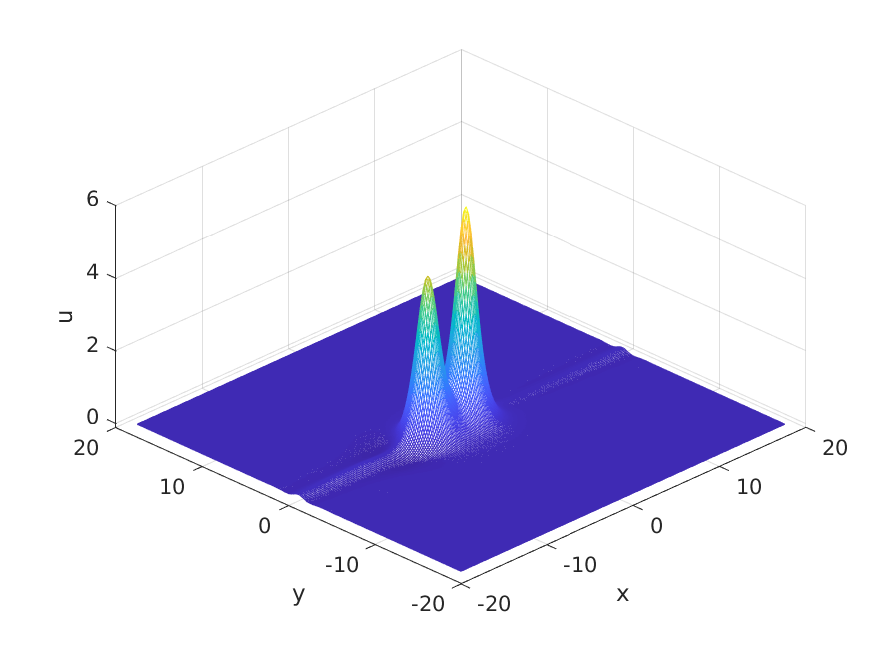}
 \includegraphics[width=0.49\hsize]{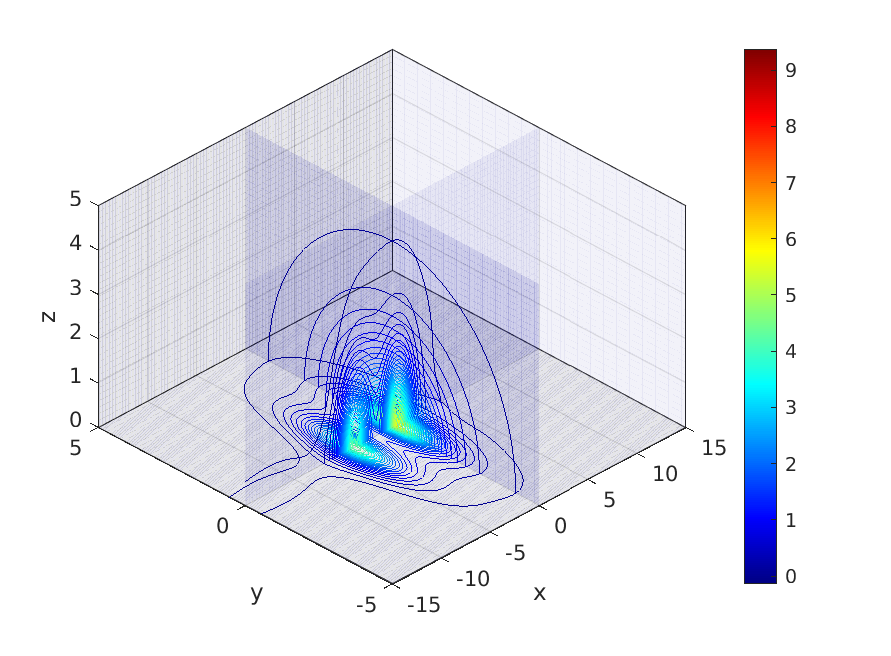}\\ 
 \includegraphics[width=0.49\hsize]{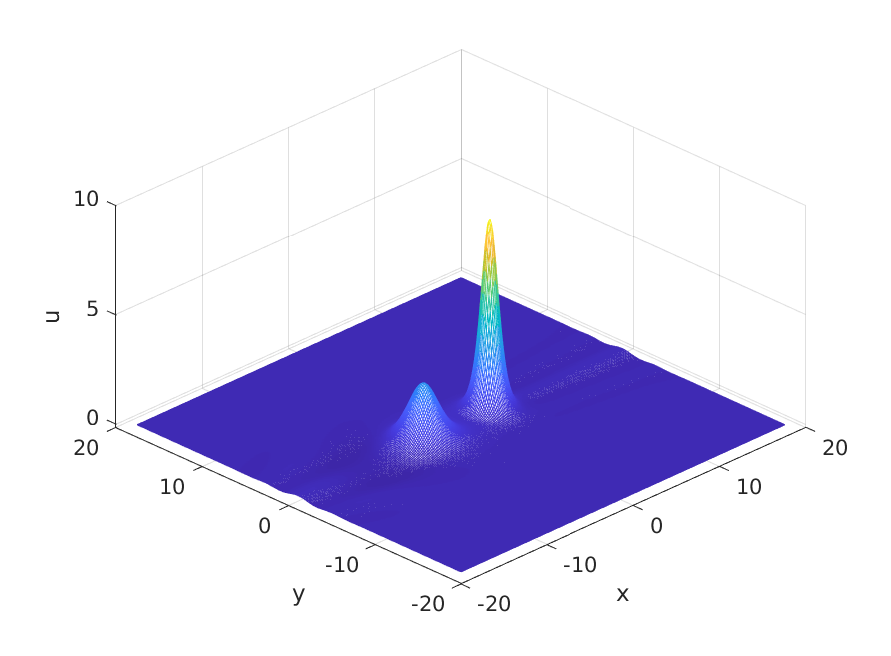}
 \includegraphics[width=0.49\hsize]{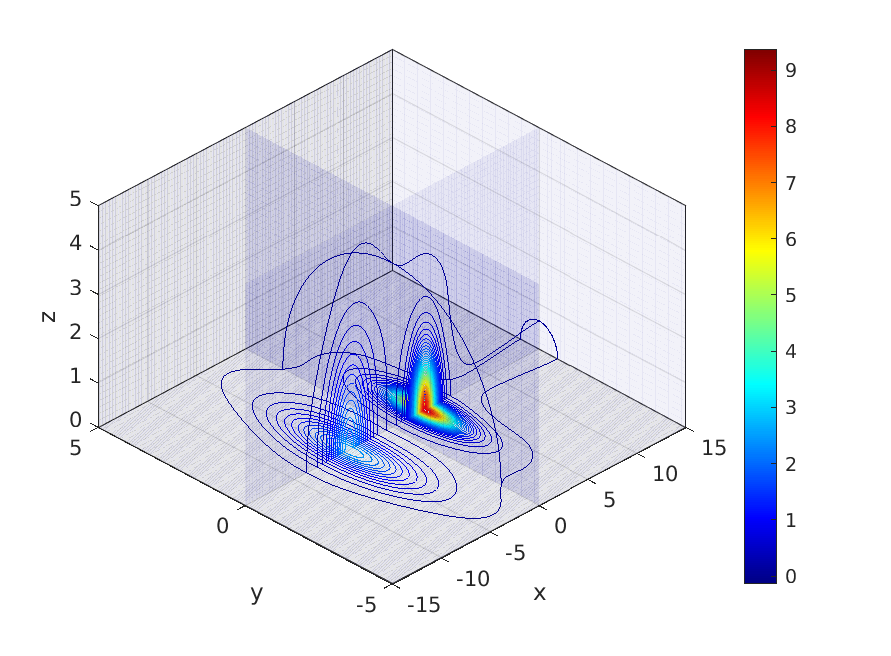}\\%
\caption{Head-on soliton interaction for initial data 
$u_0= Q_{2}(x+10,y,z)+Q(x,y,z)$ at times $t = 6.0$ (before the interaction), $t=7.5, 10.5$ (after the interaction), in a co-moving frame 
with the soliton that is initially at the origin (i.e., $v_x=1$). 
Left: interaction in 2D (projected onto the plane $z=0$). Right: the corresponding  
3D  isocurves on the slices of the coordinate planes.}
\label{Chase_sol_snaps_contour}
\end{figure}

It can be seen in Fig.~\ref{Chase_sol_snaps_contour} that the faster soliton will hit the 
slower soliton around \( t=10 \) (note 
that we are still in a co-moving frame with \( v_{x}=1 \)). 
The collision, 
while not really elastic, does preserve the number of solitons, it also exchanges their features, though not exactly.  
The smaller soliton becomes even smaller and (in a co-moving frame) shifts back (thus, the speed is slower than  $c=1$), the faster soliton grows, and thus, moves faster forward (in the positive $x$-direction), some radiation is emitted, outgoing as before in a cone-type region around the negative $x$-axis. We refer to this interaction as quasi-elastic.

\begin{figure}[!htb]
 \includegraphics[width=0.51\hsize]{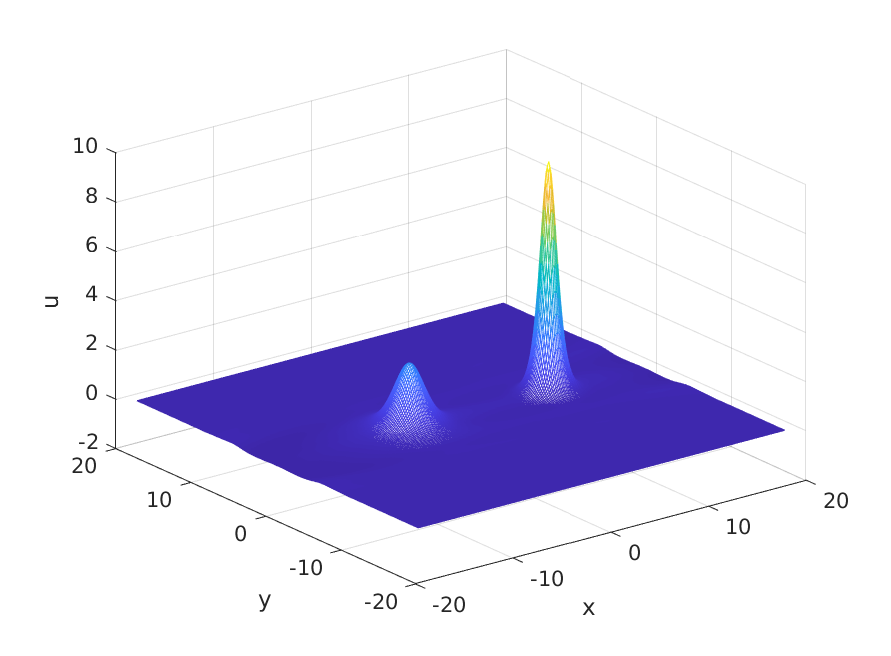}
 \includegraphics[width=0.46\hsize]{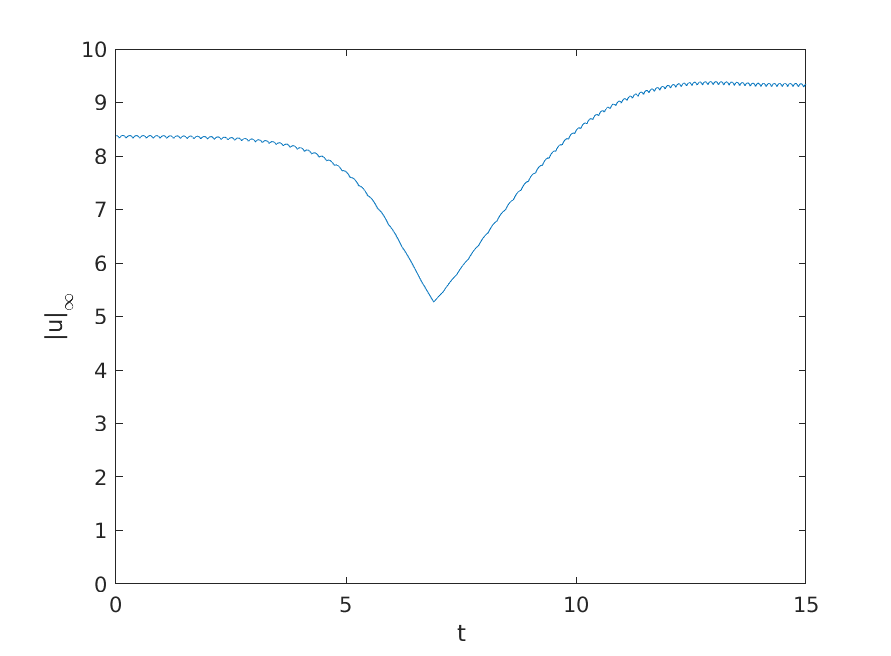}
 \includegraphics[width=0.49\hsize]{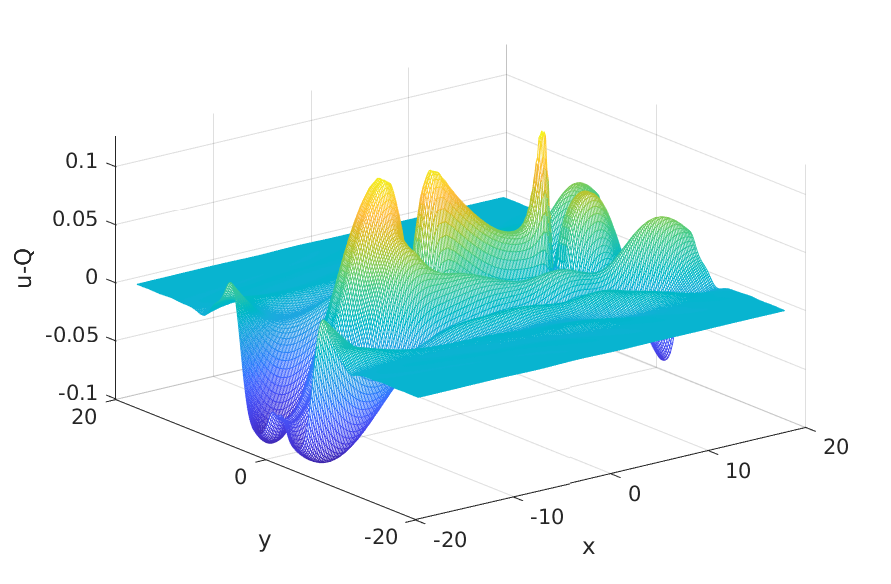}
 \includegraphics[width=0.49\hsize]{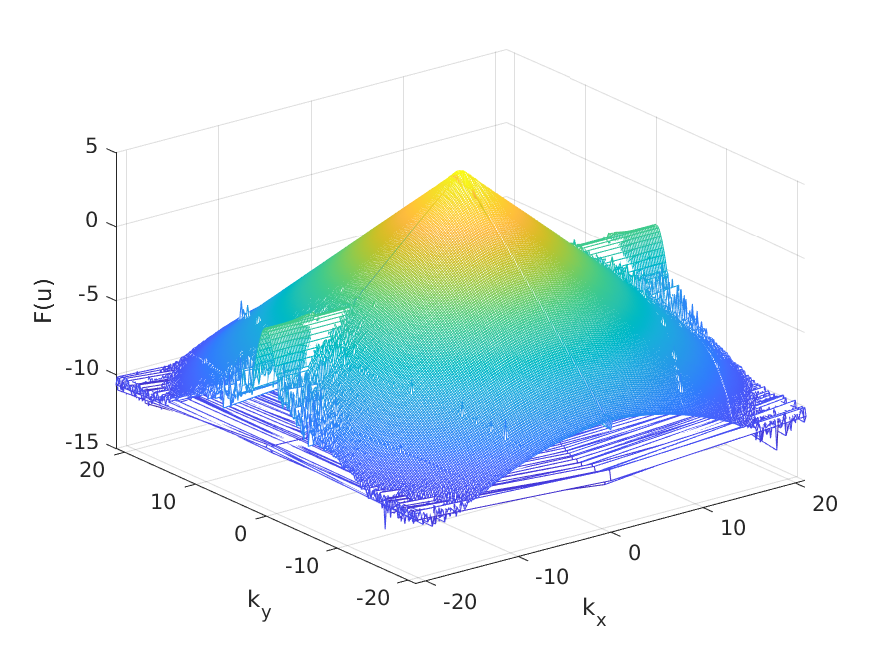}
\caption{ZK solution for initial data $u_0= Q_{2}(x+10,y,z)+Q(x,y,z)$: 
solution (projected onto $z=0$) at $t=15$ (top left), 
time dependence of the $L^{\infty}$ norm (top right),  
the difference of the solution with two rescaled solitons at $t=15$
(bottom left), 
the Fourier coefficients at $t=15$, projected onto $k_{z}=0$ (bottom right).}
\label{Chase_sol}
\end{figure}

We continue the simulation till $t=15$ and plot the snapshot of that solution in the top left of Fig.~\ref{Chase_sol}, again as the projection onto the $z=0$ plane. 
The top right of the same figure shows the $L^{\infty}$ norm. The appearance of the larger soliton can be seen in the time evolution of the $L^\infty$ norm around $t=12$. We find the scaling parameter $c$ for that soliton from \eqref{E:c} using $t=15$. 
For the second soliton we find the local maximum of the $L^\infty$ norm, where the first soliton is not present.  
The difference of the solution at $t=15$ and the two rescaled and appropriately shifted solitons 
is shown on the bottom left of Fig.~\ref{Chase_sol} and is less than $10^{-1}$. 
The latter suggests that the asymptotic solution on $\mathbb{R}^{3}$ is a superposition of two solitons (recall that we approximate this setting 
on the three-dimensional torus, 
where the radiation cannot escape to infinity), thus, confirming the soliton resolution at the (numerical) final state.  
The Fourier coefficients on the bottom right of Fig.~\ref{Chase_sol} indicate that 
the solution is numerically well resolved. 

\smallskip


\underline{\bf Case (b): {Twin Solitons.}}

\begin{figure}[!htb]
\includegraphics[width=0.49\hsize]{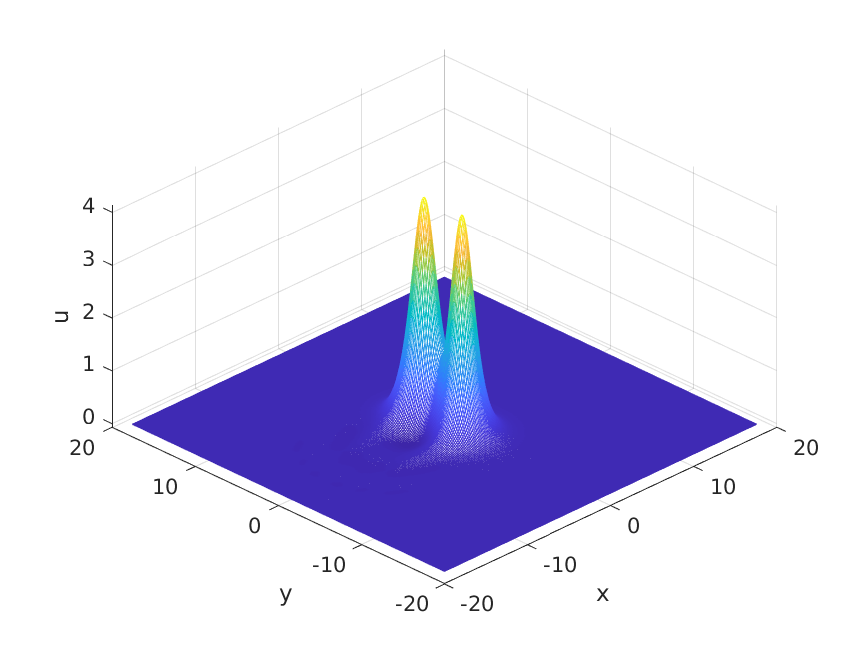}
 \includegraphics[width=0.49\hsize]{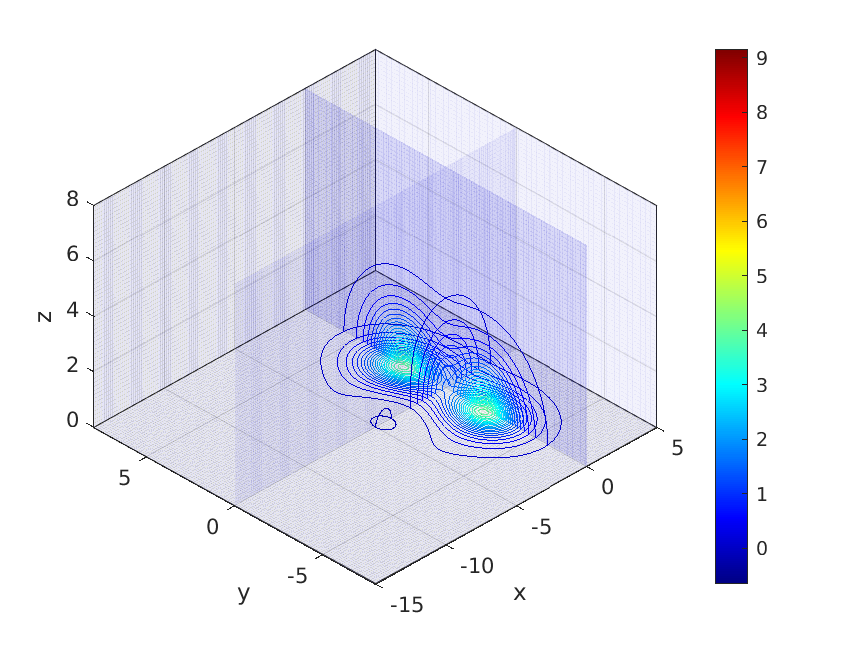}\\
\includegraphics[width=0.49\hsize]{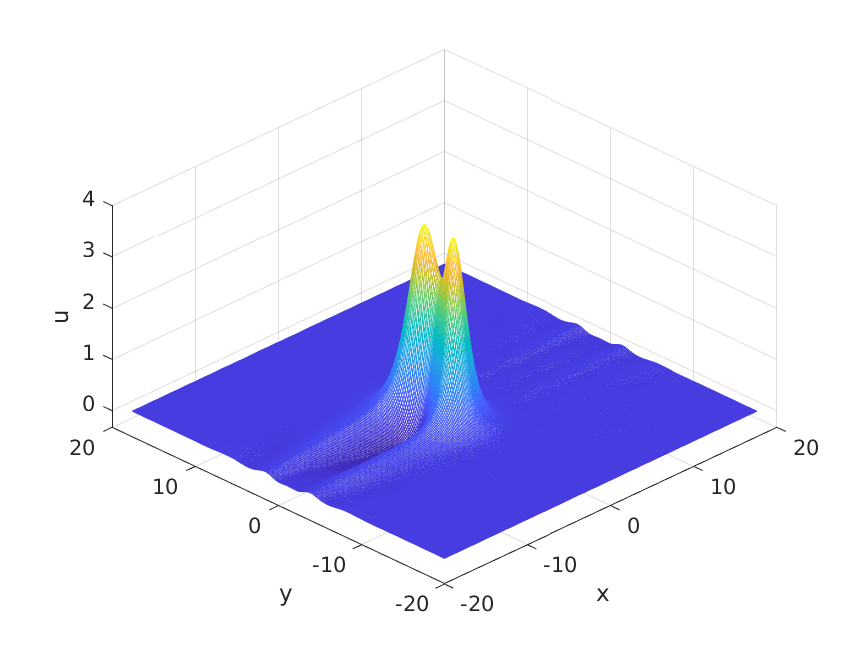}
 \includegraphics[width=0.49\hsize]{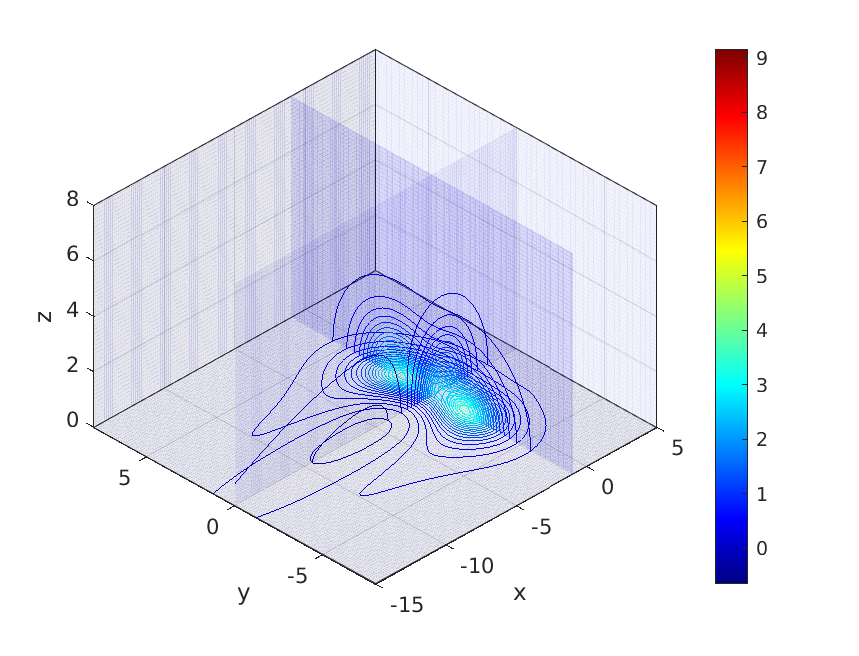}\\ 
 \includegraphics[width=0.49\hsize]{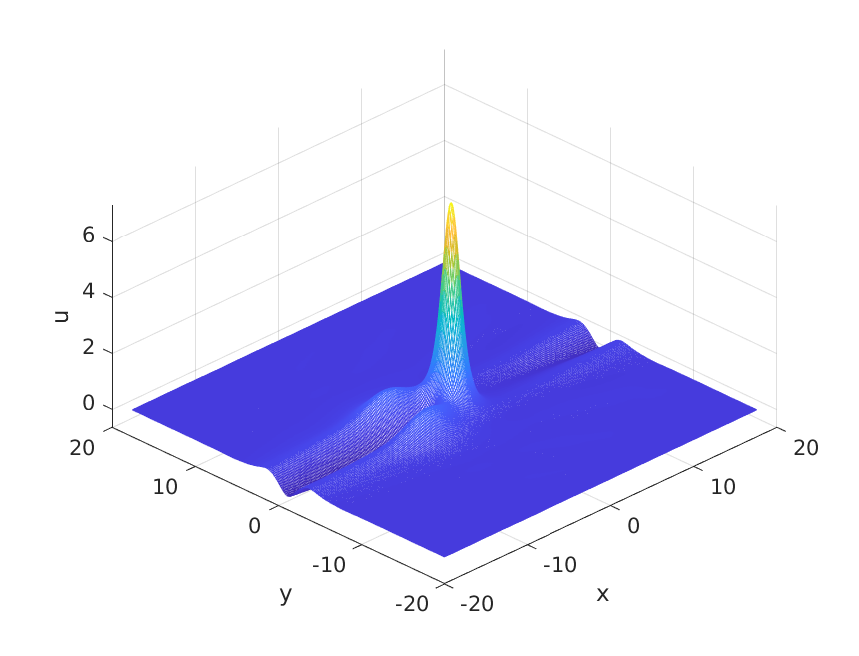}
 \includegraphics[width=0.49\hsize]{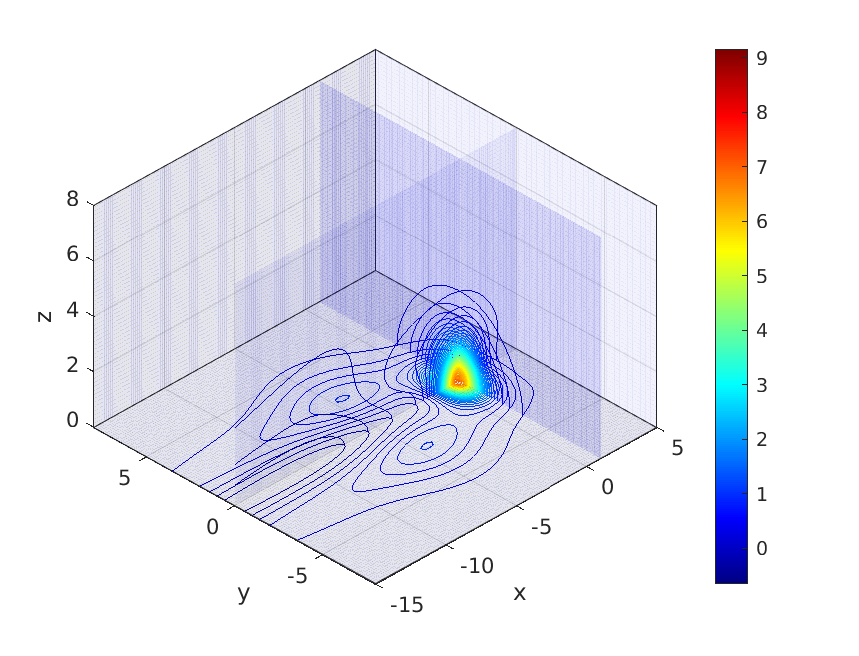}\\%
\caption{Snapshots of the ZK solution of strong soliton interactrion with $u_0 = 
Q(x,y-a,z)+Q(x,y+a,z)$, $a=\pi \, L/8$ at $t = 1.5, 6, 10.5$. 
Left: 2D projections (onto $z=0$ plane). Right: 3D contour plots on the slices of the coordinate planes.}
\label{ZK_doubleSol_slices}
\end{figure}

The next example that we consider has initial data of two solitons next to each other shifted symmetrically in one of the non-leading axis, either in the $y$-axis or in the $z$-axis, for instance, 
\begin{equation}\label{E:twin}
u(x,y,z,0) = Q(x,y-a,z)+Q(x,y+a,z), ~~a > 0.
\end{equation}

In Fig.~\ref{ZK_doubleSol_slices} we show the snapshots of the solution with   
initial data \eqref{E:twin} and $\displaystyle a=\frac{\pi \, L}{8}$, noting that the periods are $2\pi L_{x}$, $2\pi L_{y}$, $2\pi L_{z}$, respectively in each direction, and $L_{x}=L_{y}=L_{z}=L$. 
The two nearby solitons clearly merge into a single soliton of larger velocity.
We refer to this interaction as a strong insteraction, since the number of solitons has changed.

\begin{figure}[!htb]
 \includegraphics[width=0.51\hsize]{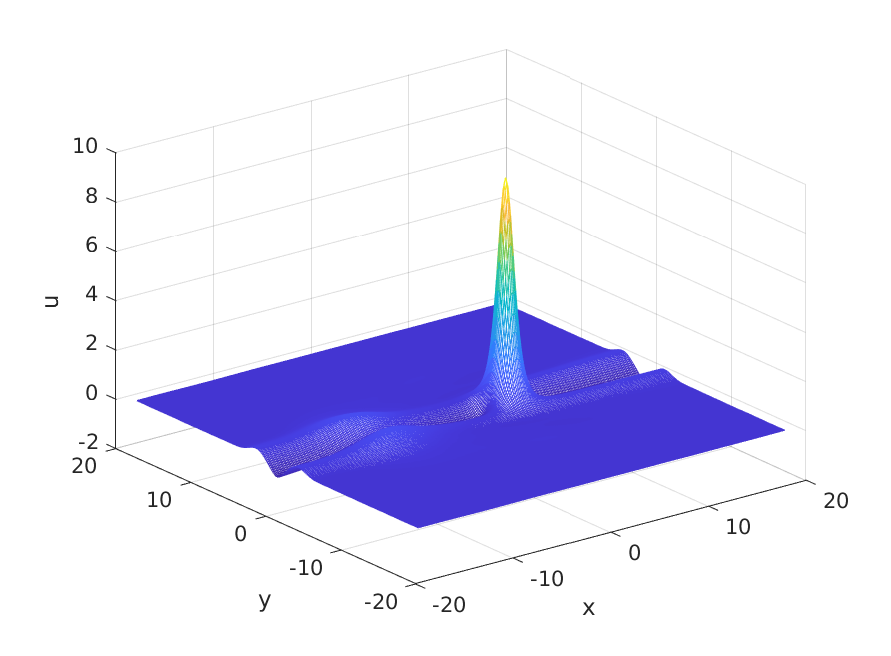}
 \includegraphics[width=0.46\hsize]{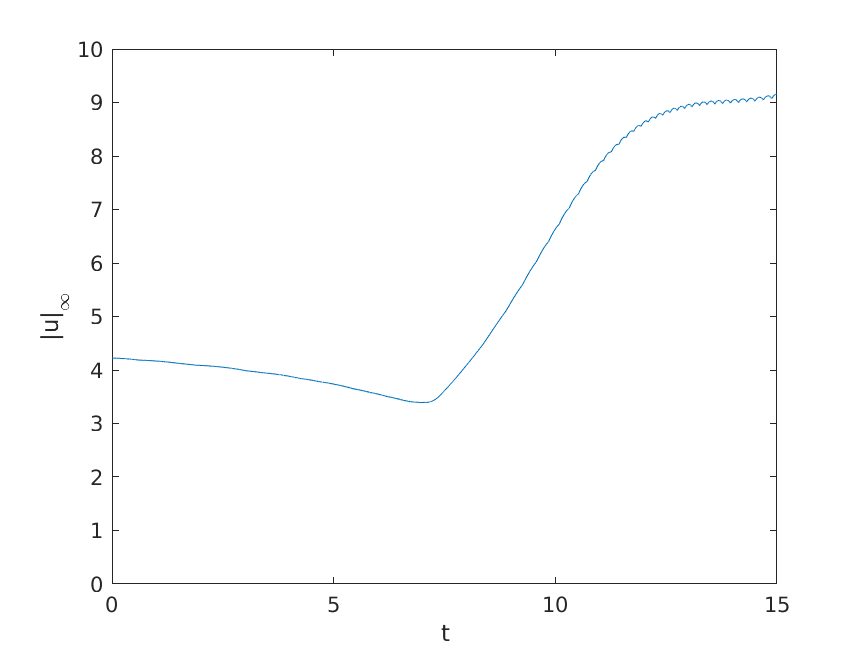}
 \includegraphics[width=0.49\hsize]{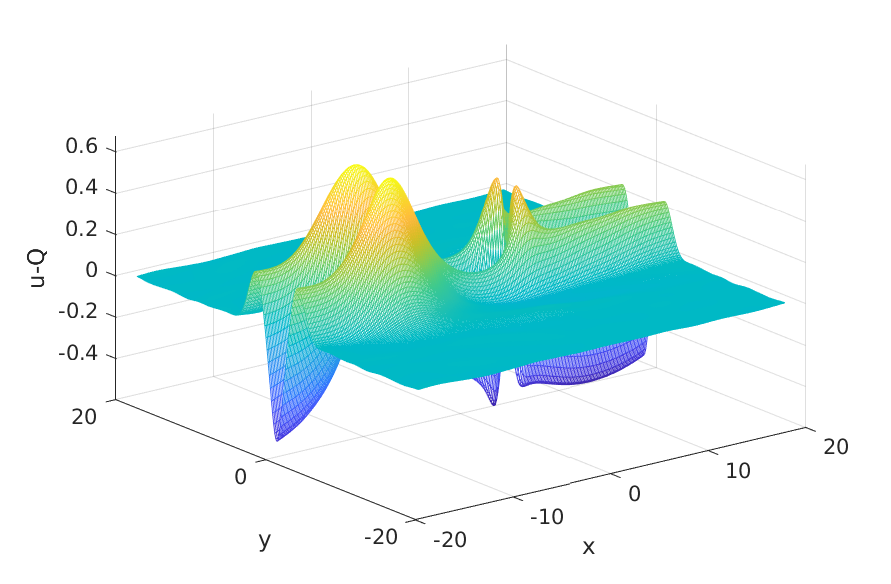}
 \includegraphics[width=0.49\hsize]{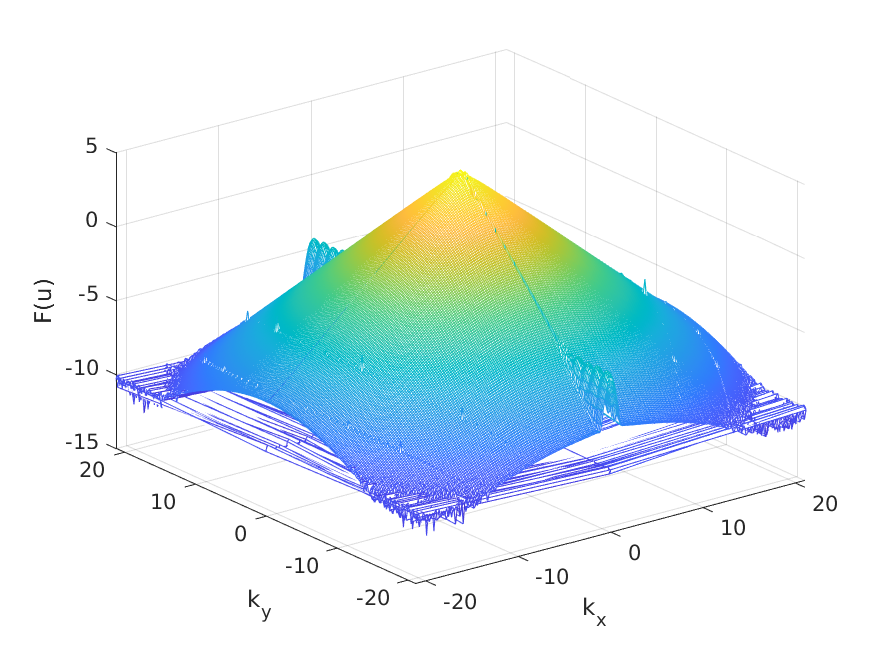}
\caption{Soliton resolution for $u_0 = 
Q(x,y-a)+Q(x,y+a)$, $a=\pi \, L/8$. The solution (projection onto $z=0$) for $t=15$ (top left), 
the $L^{\infty}$ norm (top right), the difference between the resulting profile and a rescaled soliton $Q_c$ (bottom left), the Fourier coefficients at $t=15$ (bottom right).} 
\label{double_sol}
\end{figure}

Continuing tracking the solution up to $t=15$, we show a projection onto the $z$-plane of the solution at that time on the top left of Fig. \ref{double_sol}. 
The difference with a rescaled soliton $Q_c$ is on the bottom left of the same figure; this suggests that the final state in this interaction is indeed a soliton plus 
radiation. 

This result is also confirmed by the time dependence of the $L^{\infty}$ norm of the 
solution on the top right of Fig.~\ref{double_sol} (from which we 
determined the scaling parameter $c$ for the rescaled soliton $Q_c$).  The Fourier 
coefficients on the bottom right of the same figure indicate that the 
solution is numerically well resolved. 

\underline{\bf {Case (c): Off-set solitons.}}

\begin{figure}[!htb]
\includegraphics[width=0.49\hsize]{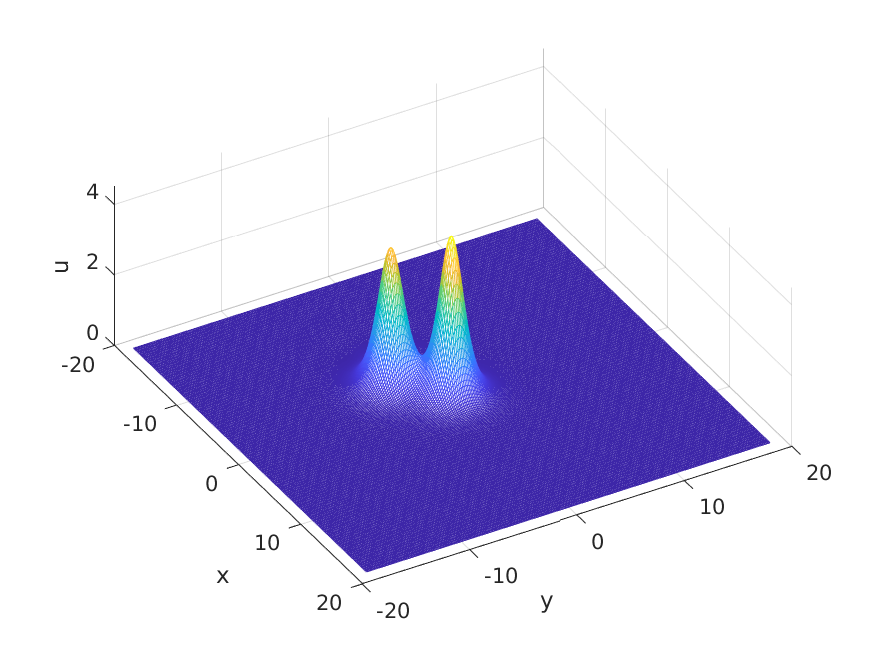}
  \includegraphics[width=0.49\hsize]{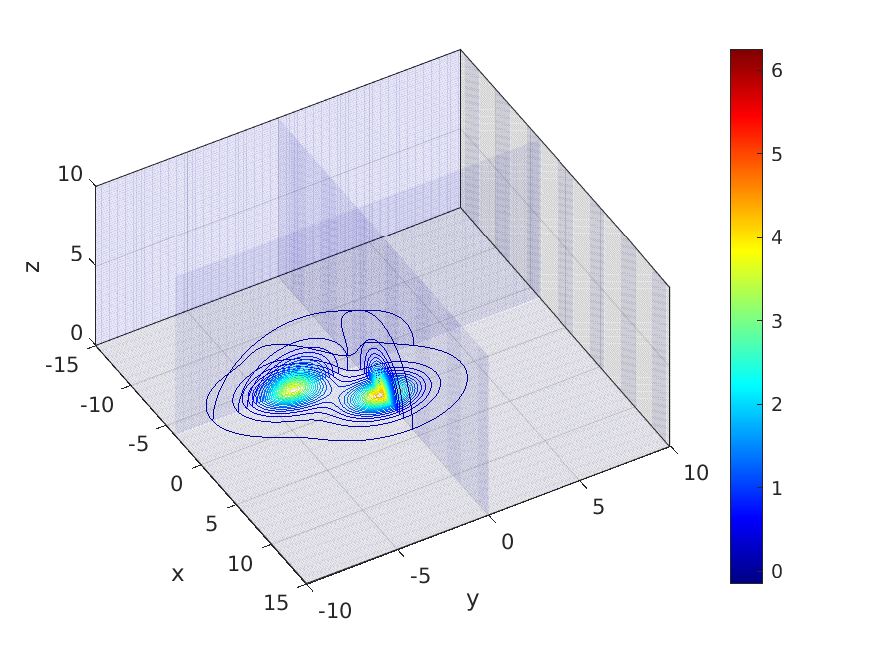}\\      
  \includegraphics[width=0.49\hsize]{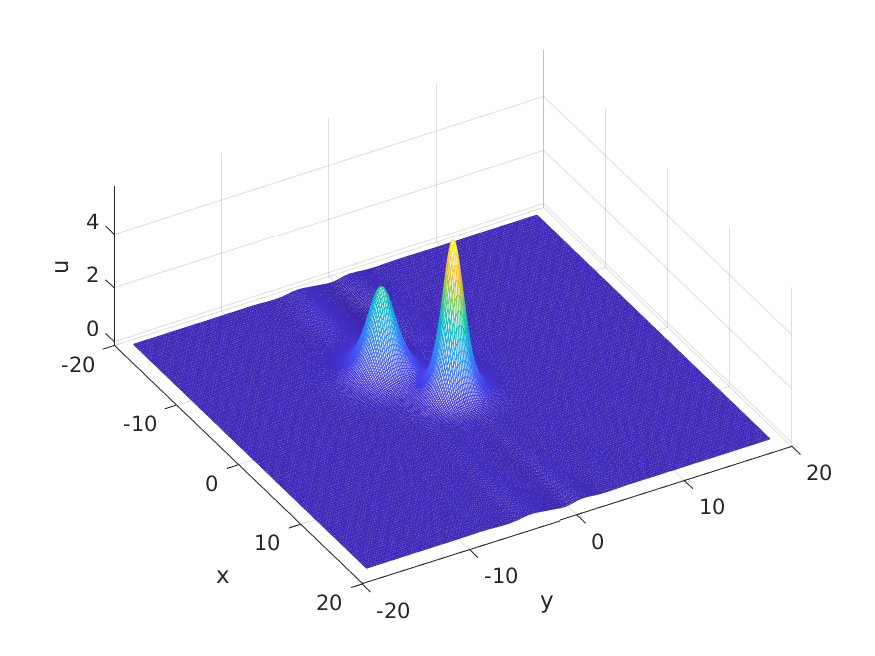}
  \includegraphics[width=0.49\hsize]{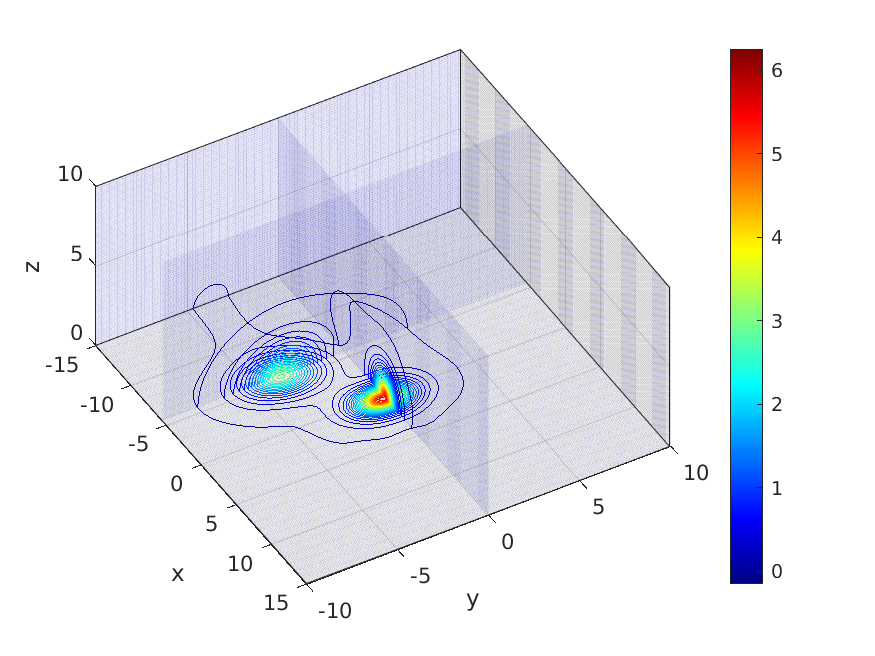}\\
\includegraphics[width=0.49\hsize]{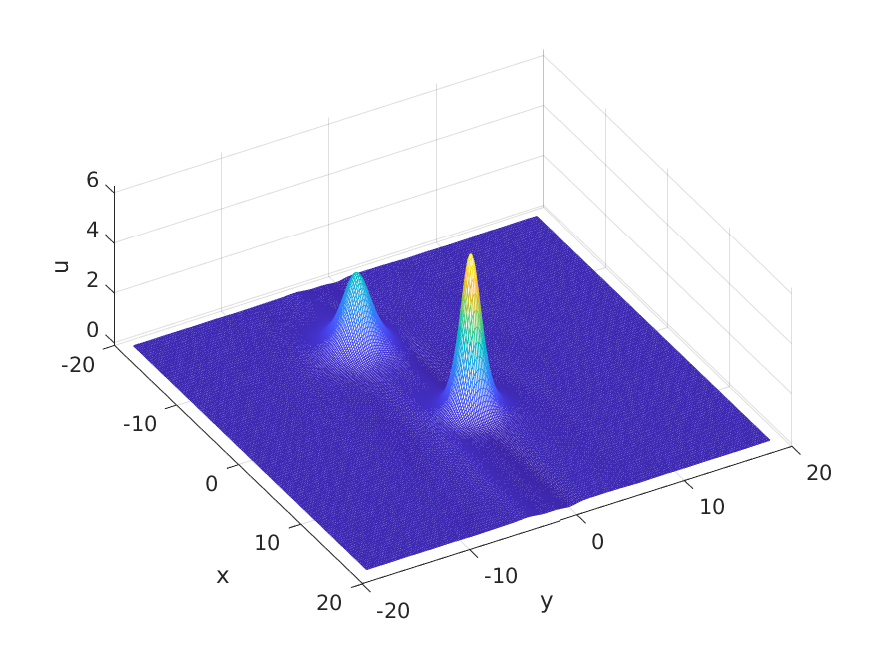}
\includegraphics[width=0.49\hsize]{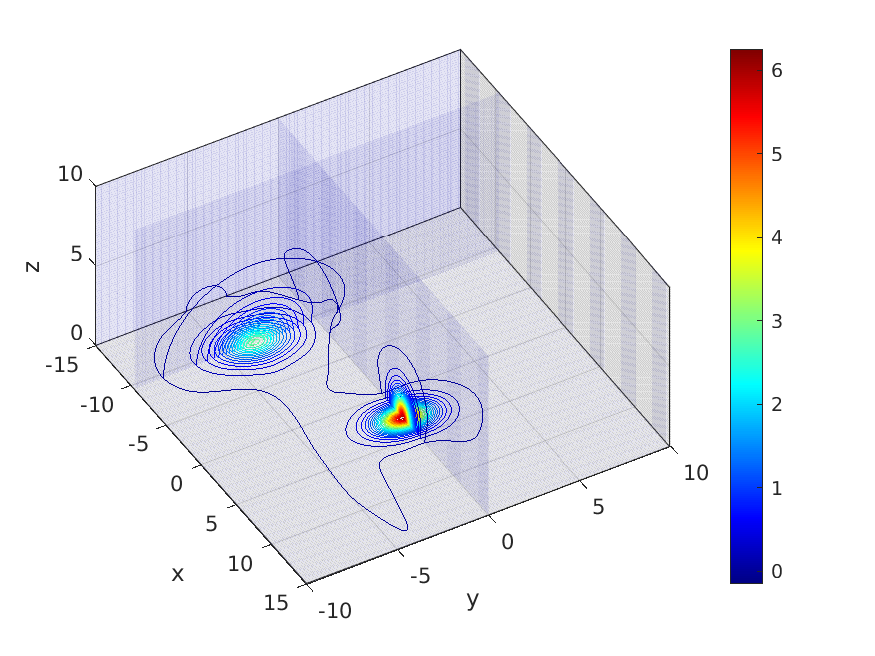}\\
\caption{Snapshots of ZK solution with off-set initial data 
$u_0(x,y,z) = Q(x,y,z,0) + Q(x+a, y+a, z, 0)$, $a = 3/8\pi$, at $t = 1.5, 6, 15$. 
Left: 2D projections onto the $z$-plane. Right: 3D contour plots on 
the coordinate plane slices.}
\label{offset_slices}
\end{figure}

Finally, we consider two identical solitons displaced in both $x$ and $y$ (or in $x$ and $z$) directions to break the symmetry in the interaction. 
Specifically, we study the initial data of the form  
\begin{equation}
u(x,y,z,0) = Q(x,y,z) + Q(x+a, y+a, z), ~~a>0. 
\label{offset_init}
\end{equation}

As an example, we take $a = \dfrac38 \,\pi$. We observe a quasi-elastic 
soliton interaction: some mass is transferred to the front soliton, the number of solitons remain the same, though their characteristics have changed: the back soliton is slower 
than the original $Q$ and the fast soliton is faster than $Q$ after the interaction. 
One can notice an increasing separation of them 
due to their different speeds.

\begin{figure}[!htb]
 \includegraphics[width=0.51\hsize]{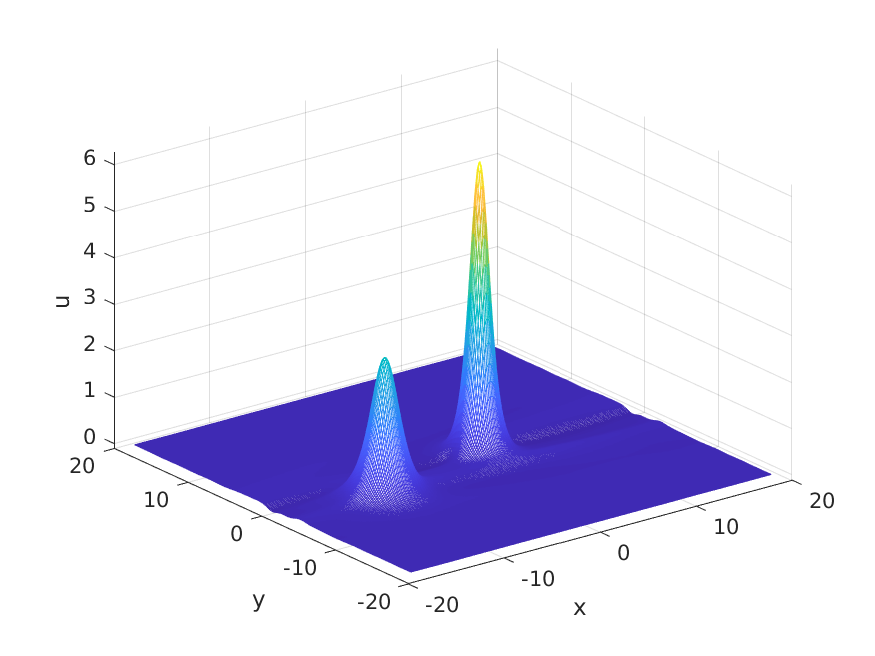}
 \includegraphics[width=0.46\hsize]{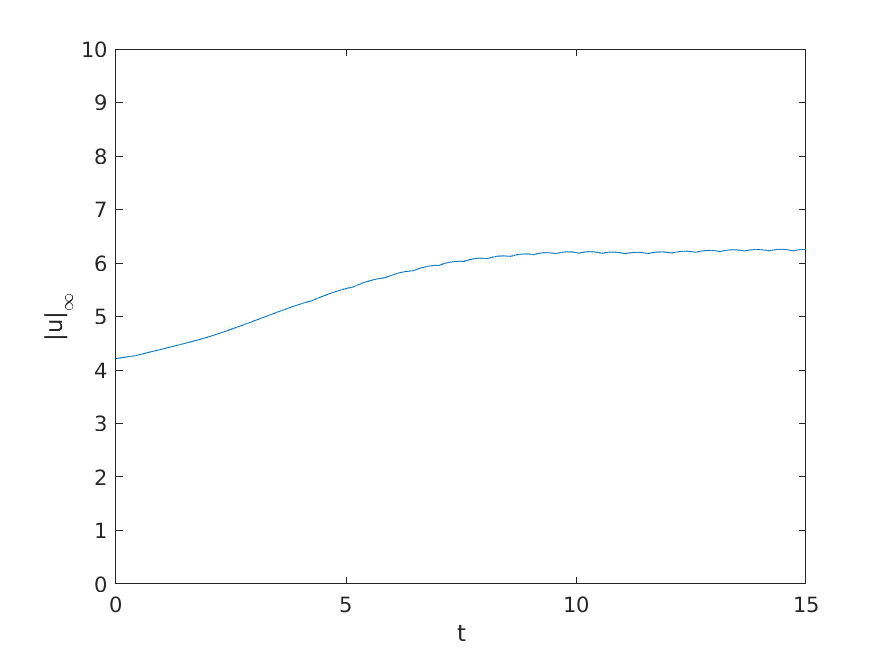}
 \includegraphics[width=0.49\hsize]{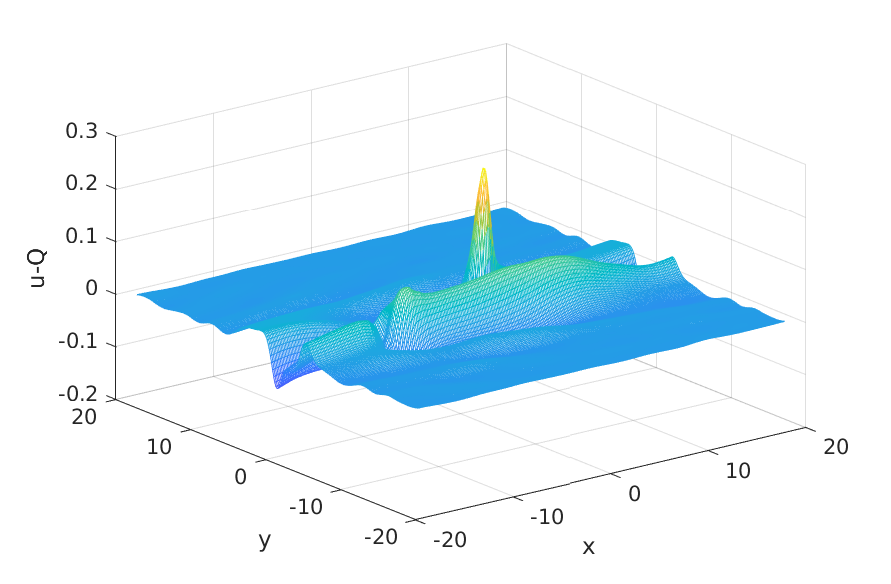}
 \includegraphics[width=0.49\hsize]{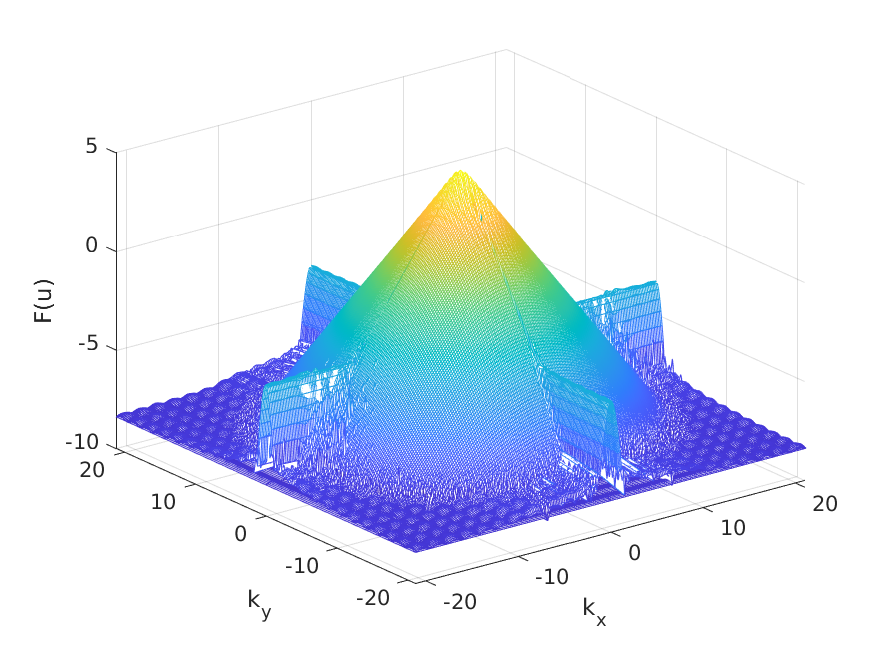}
\caption{Soliton resolution for the off-set initial data \eqref{offset_init}, $a=\frac38 \pi$: final profile (at $t=15$), projected onto $z=0$ (top left),  
the time dependence of the $L^{\infty}$ norm (top right),
the difference between the leading profile and the rescaled soliton (bottom left), 
the Fourier coefficients at $t=15$ (bottom right).}
\label{offset_sol}
\end{figure}

The projection of the solution at the `final' time $t=15$ onto the $z=0$ plane 
is shown on the top left of Fig.~\ref{offset_sol}.  
The difference of the solution at $t=15$ and the leading soliton is on the bottom left of the same figure; it indicates that the final state can again be interpreted as a 
superposition of solitons. 
This result is also confirmed by the $L^{\infty}$ norm of the 
solution on the top right of Fig.~\ref{offset_sol}. The Fourier 
coefficients on the bottom right of the same figure indicate that the 
solution is numerically well resolved.

We conclude that the soliton stability, soliton resolution hold for the 3D ZK equation as well as the interaction of solitons in this equation can be quasi-elastic or strong.


\end{document}